\documentclass[a4paper,10pt]{article}
\usepackage{latexsym}
\usepackage{amssymb}
\usepackage{amsfonts}
\usepackage{epsfig}
\usepackage{amsmath}
\usepackage{amscd}
\usepackage[active]{srcltx}
\usepackage{psfrag}
\usepackage{array}
\usepackage{tikz}
\usepackage[all]{xy}
\usepackage{xcolor}% http://ctan.org/pkg/xcolor
\usepackage{colortbl}% http://ctan.org/pkg/colortbl
\usepackage{graphicx}% http://ctan.org/pkg/graphicx
\usepackage{color}
\usepackage{amsmath}
\usepackage{multirow,bigdelim}
\usepackage{enumerate}
\usepackage{booktabs}
\usepackage{authblk}
\usepackage{verbatim}
\usepackage{lscape}

\definecolor{pink}{rgb}{1,0,1}

\newtheorem{defn0}{Definition}[section]
\newtheorem{prop0}[defn0]{Proposition}
\newtheorem{thm0}[defn0]{Theorem}
\newtheorem{lemma0}[defn0]{Lemma}
\newtheorem{corollary0}[defn0]{Corollary}
\newtheorem{example0}[defn0]{Example}
\newtheorem{conjecture0}[defn0]{Conjecture}
\newtheorem{notation0}[defn0]{Notation}
\newtheorem{remark0}[defn0]{Remark}
\newtheorem{assumption0}[defn0]{Claw $d$-tree hypothesis}

\newenvironment{defn}{\begin{defn0} \rm}{\end{defn0}}
\newenvironment{prop}{\begin{prop0}}{\end{prop0}}
\newenvironment{thm}{\begin{thm0}}{\end{thm0}}
\newenvironment{lem}{\begin{lemma0}}{\end{lemma0}}
\newenvironment{cor}{\begin{corollary0}}{\end{corollary0}}
\newenvironment{example}{\begin{example0}\rm}{\end{example0}}
\newenvironment{rem}{\begin{remark0}\rm}{\end{remark0}}

\newenvironment{notation}{\begin{notation0}\rm}{\end{notation0}}
\newenvironment{assumption}{\begin{assumption0}\rm}{\end{assumption0}}

\newcommand{\N}[1]{N_{#1}}
\newcommand{\FF}[2]{\mathcal{F}_{#1}(#2)}
\newcommand{\xx}{\mathtt{X}}

\newcommand{\yy}{\mathtt{Y}}
\newcommand{\zz}{\mathtt{Z}}
\renewcommand{\ss}{\mathtt{S}}
\newcommand{\vv}{\mathtt{V}}
\newcommand{\ww}{\mathtt{W}}
\renewcommand{\aa}{\mathtt{A}}
\newcommand{\cc}{\mathtt{C}}
\renewcommand{\gg}{\mathtt{G}}
\renewcommand{\tt}{\mathtt{T}}
\newcommand{\CC}{\mathbb{C}}
\newcommand{\J}[1]{\mathbf{J}_{#1}}
\newcommand{\MM}{\mathcal{M}}

\newcommand{\baa}{\underline{\aa}}
\newcommand{\bcc}{\underline{\cc}}
\newcommand{\bgg}{\underline{\gg}}
\newcommand{\btt}{\underline{\tt}}
\newcommand{\bxx}{\underline{\xx}}
\newcommand{\byy}{\underline{\yy}}
\newcommand{\lra}{\longrightarrow}
\newcommand{\m}{\mathbf{m}}
\renewcommand{\o}{\otimes}
\newcommand{\op}[1]{\otimes^{#1}W}
\newcommand{\opG}[1]{(\op{#1})^G}

\newcommand{\codim}{\mathrm{codim}\,}
\newcommand{\rk}{\mathrm{rank}\,}
\newcommand{\Tf}{\mathrm{Tflat}\,}
\newcommand{\Hom}{\mathrm{Hom}\,}
\newcommand{\Par}{\mathrm{Par}\,}
\renewcommand{\Im}{\mathrm{Im}\,}

\newcommand{\Ba}{\mathcal{B}}
\newcommand{\eq}{\mathfrak{eq}}
\newcommand{\eeq}{\mathrm{eq}}
\newenvironment{proof}{\noindent {\textsc{Proof.}}}{$\square$ \vspace{3mm}}

\newif\ifprivate
\privatetrue
\def\???{\ifprivate {\bf {???}} \marginpar{{\Huge {\bf ?}}}
\else \fi}

\def\p{\mathbb{P}}

\allowdisplaybreaks

\begin{document}

\title{Complete intersection for equivariant models}
\author[1]{M. Casanellas \thanks{marta.casanellas@upc.edu}}
\author[1]{J. Fern\'andez-S\'anchez \thanks{jesus.fernandez.sanchez@upc.edu}}
\author[2]{M. Micha{\l}ek \thanks{wajcha@berkeley.edu}} 
\affil[1]{Departament de Matem\`atiques, Universitat Polit\`ecnica de Catalunya}
\affil[2]{Freie Universit{\"a}t, Berlin, Germany; UC Berkeley, USA; Polish Academy of Sciences, Warsaw}

% \author{M. Casanellas, J. Fern{\'a}ndez-S{\'a}nchez \thanks{The first two authors were partially supported by Spanish government MTM2012-38122-C03-01/FEDER
% and Generalitat de Catalunya 2014SGR634.}, M. Micha{\l}ek \thanks{The third autor was supported by National Science Center grant SONATA UMO-2012/05/D/ST1/01063.} \\ 
% \small{Departament de Matem\`atiques, Universitat Polit\`ecnica de Catalunya,}\\
% %\small{Av. Diagonal 647, 08028-Barcelona, Spain.} \\ \small{e-mail: marta.casanellas@upc.edu, jesus.fernandez.sanchez@upc.edu} \\
% \small{Freie Universit{\"a}t, Berlin, Germany; UC Berkeley, USA; Polish Academy of Sciences, Warsaw}

\maketitle
% \blfootnote{The first two authors were partially supported by Spanish government MTM2012-38122-C03-01/FEDER
% and Generalitat de Catalunya 2014SGR634. The third autor was supported by National Science Center grant SONATA UMO-2012/05/D/ST1/01063.}
\begin{abstract}
Phylogenetic varieties related to equivariant substitution models have been studied largely in the last years. One of the main objectives has been finding a set of generators of the ideal of these varieties, but this has not yet been achieved in some cases (for example, for the general Markov model this involves the open ``salmon conjecture'', see \cite{allman_salmon}) and it is not clear how to use all generators in practice. % should it be useful for biologists
%\red{mm:should it $\rightarrow$ it should. Also maybe we can make it a little softer like 'not clear how to apply/use all generators/ideal in biology/practice'}.
Motivated by applications in biology, we tackle the problem from another point of view. The elements of the ideal that could be useful for applications in phylogenetics only need to describe the variety around certain \emph{points of no evolution} (see \cite{casfer2008}).
%\red{mm: abstract is quite long. Maybe the bracket could be change just to a citation?}. We call them \red{mm: maybe insted of new sentence just '- generic ..'}
 %\red{mm: In this paper, we manage to $\rightarrow$ We}
 We produce a collection of explicit equations that describe the variety on a Zariski open neighborhood of these points (see Theorem \ref{invariantsCI}). Namely, for any tree $T$ on any number of leaves (and any degrees at the interior nodes) and for any equivariant model on any set of states $\kappa$, we compute the codimension of the corresponding phylogenetic variety. We prove that this variety is smooth at general points of no evolution, and provide an algorithm to produce a complete intersection that describes the variety around these points.
 %\red{The last sentence is quite long. Also we do the first induction step only for most well-know models.}
%We study phylogenetic invariants of equivariant phylogenetic models, including General Markov Model. We provide an explicit description of a Zariski open %neighbourhood of the no evolution point in the model as a complete intersection. In other words, we provide a minimal possible number of explicitly constructed %phylogenetic invariants, that allow to describe the model at biologically meaningful points. Our work is based on previous inductive constructions of %phylogenetic invariants \cite{}, however is motivated mostly by applications, as the number of phylogenetic invariants we construct is much lower.
\end{abstract}

\section{Introduction}

In the recent years there has been a huge amount of work done on phylogenetic varieties -- we advise the reader to consult e.g.~\cite{allman2003, AllmanRhodeschapter4, MR2314109, eriksson2004, MR2811326, sturmfels2004} and references therein. These algebraic varieties contain the set of joint distributions at the leaves of a tree evolving under a  Markov model of molecular evolution. From the biological point of view, these varieties are interesting because they provide new tools of non-parametric inference of phylogenetic trees. At present, the algebraic/geometric framework of phylogenetic varieties has allowed proving the identifiability of parameters of certain evolutionary models widely used by biologists \cite{allman2006a,allman2008}, proposing new methods of model selection \cite{Kedzierska2012}, and producing new phylogenetic reconstruction methods \cite{eriksson2005,casfer2015}.

From the mathematical point of view, there has been a great effort in finding a whole description of the ideal of these phylogenetic varieties \cite{allman2004,sturmfels2004, draisma2008, MR3016420}. Still, for some models, many questions remain open for trees on an arbitrary number of leaves $n$. Indeed, if one is interested in using these algebraic tools with real data, one would need a small set of generators of the ideal (rather than a description of the whole ideal); it would also be desirable to know the degree at which the ideal is generated \cite{MR3336833, MR3092693, sturmfels2004}; and as the codimension of the variety is exponential in $n$, it is necessary to distinguish between generators that only account for the underlying evolutionary model, those that account for the tree topology, and those that could be useful for inferring the numerical parameters (see \cite{AllmanRhodeschapter4} for a good introduction to this topic).

For instance, the authors of \cite{sturmfels2004} and \cite{AllmanRhodeschapter4} raised the question whether knowing complete intersection containing the phylogenetic variety and of the same dimension would be enough. Eminently, for biological applications it is only relevant to know the description of the variety around the points that make sense \emph{biologically} speaking. If these points are smooth, then a complete intersection can define the variety on a neighborhood of these points.

This is the approach that was considered in \cite{casfer2008, hagedorn2000a, MR3159089} for the particular model of Kimura 3-parameter and is the same goal that we pursued in \cite{casfermich2015} for abelian group-based models. In the present paper we address this problem for the more general class of \emph{G-equivariant models}, which contains more general algebraic models of interest to biologists
%\tp{jf: ``\textbf{all} the algebraic models of interest to biologist'' sounds a little excessive to me; I'd suggest ``most of the algebraic models of interest...'' or remove ``all''}
(for example the strand-symmetric and the general Markov models). We give an explicit algorithm to construct a complete intersection that describes the variety on a dense open subset around a generic \emph{point of no evolution}. Points of no evolution represent molecular sequences that remain invariant from the common ancestor to the leaves of the tree. Points in the phylogenetic variety that arise from biologically meaningful parameters are supposed to be near these points of no evolution (otherwise phylogenetic inference could not be made), and therefore it is important to study the variety around these points. Also, as we describe the variety at a dense open subset containing these points, we cover most (actually, all but possibly a subset of smaller dimension) of the biologically meaningful points of the variety.
Also, in the same papers mentioned above, it is argued that a complete intersection can contain points in other irreducible components that can mislead the results in practice. However, the complete intersection we give contains a regular sequence of the edge invariants, which are known to be phylogenetic invariants \cite{casfer2010}. Therefore, the other irreducible components of the complete intersection do not contain other phylogenetic varieties.
%Thus, all points in the complete intersection o Theorem that areclose to a no evolution point and that once it is assumed that a point $p$ lies in one of the phylogenetic varieties and is a point close to a no evolution point, the complete inte$p$ cannot lie in any of the other irreducible components of
%\red{I am no expert here, but I guess that points that are biologically interesting should be 'near' points of no evolution - as changes in DNA appear with small probabilities. This would be a good motivation to study neighbourhood of those. Otherwise a description on any open set would also 'cover most of the biologically meaningful points' and would be just as good.}

We prove first that these points are non-singular and therefore the variety can be described locally at these points by the smallest possible number of equations, the codimension of the variety.
A system of generators of the local complete intersection can be explicitly computed. The degree of these generators is low and depends on a local description of claw trees related to the interior nodes of the tree and of the multiplicities of the permutation representation of the group $G \subset \mathfrak{S}_{\kappa}$. For example, for the biologically interesting models mentioned above, the complete intersection we provide has generators of degree at most 13. One should contrast this to the generators of the complete intersection given in \cite{steel1993} for the Kimura 3-parameter model, which had exponential degree in the number of leaves.   Our approach is also useful in case one wants to use differential geometry for this variety (for example to compute the distance of a point to this variety, \cite{draisma2015}).

%\red{Well... for 3K we know that degree 4 is enough to describe the whole variety set theoretically for any number of leaves and any valency, but we can keep it as an argument.}.

The
 %\red{mm: maybe we could remove 'complete' as for claw trees the problem is still open }
 description of the ideal of phylogenetic invariants for $G$-equivariant models was provided by Draisma and Kuttler \cite{draisma2008}. There are two ways to obtain the whole ideal of phylogenetic invariants for a given tree, assuming the ideals for star trees are known.
The first description relies on the ideals of star trees associated to inner vertices of a given tree - so-called flattenings. The second description is inductive, where we regard a big tree as a join of two smaller trees.

Our approach is based on the second method. %Instead of considering flattenings associated to vertices, we further flatten the tensor to a matrix corresponding to a decomposition of leaves given by an edge. A priori, the obtained matrix has rank equal to $\dim V$. However, thanks to the equivariant decomposition of $V$, the matrix has block diagonal form and we can explicitly bound the rank of each block.
We start by inducing phylogenetic invariants from smaller trees.
%There are several ways to achieve this and in our method we choose one leaf in each of the two smaller trees.
The induced phylogenetic invariants are of course not enough to provide a description for the larger tree. We complement them by so-called thin flattenings \cite{casfer2010}. They are very explicit, however still numerous. It turns out that the choice of leaves in smaller trees distinguishes specific thin flattenings. Combining those with induced invariants yields our
main result: under a minor assumption on claw trees (see \ref{assumption}) which is satisfied by the tripod on the most popular equivariant models (Jukes-Cantor, Kimura 3ST, strand symmetric) and also by the general Markov model, we provide an explicit local description of the variety associated to a model and a tree (see Theorem~\ref{invariantsCI}). Moreover, both in the starting point and in the induction process, the choices we make are almost canonical so that the complete intersection we produce is a natural one and could be reasonably used in practice.

%\red{mm: I would add somewhere that both in the starting point and in the induction process the choices we make are almost canonical. I.e. theoretical mathematicians may object to our paper that when we have a lot of equations we could just take some random combination of those and this should work. But even they should appreciate how naturally the numbers add up to a complete intersection.} %We show inductively that such a variety, in a Zariski neighborhood of a no evolution point, is a complete intersection of explicitly defined edge invariants coming from thin flattenings.

The methods used in this paper rely on basic algebraic geometry and group representation theory. It is important to note that the results of the paper hold for any $G$-equivariant model, $G \subset \mathfrak{S}_{\kappa}$, for any $\kappa$, and therefore representation theory has been  the necessary tool to deal with all these models at the same time. On the other hand, our results also hold for trees with any  number of leaves and any degrees at the interior nodes.

The approach adopted to prove the main result \ref{invariantsCI} also produces a computationally effective list of elements in the ideal of the phylogenetic variety. Indeed, the list of equations provided in Theorem \ref{invariantsCI} for a tree $T$ with $n$ leaves is constructed from equations describing locally the phylogenetic varieties of claw trees of the interior nodes of $T$ and from certain minors of the thin flattenings mentioned above. The number of equations from the thin flattenings grows exponentially with $n$, but the number of equations corresponding to claw trees does not (it grows exponentially with the maximum degree of an interior node of $T$, see Remarks \ref{rem_Nab} and \ref{rem_notexpon}). However, evaluating the minors of the thin flattenings is not the optimal way of evaluating the rank of a matrix and these equations could therefore be substituted by a numerical method such as the singular value decomposition (see \cite{eriksson2005}). The remaining equations form a set that can be useful in practice, for example for the estimation of the parameters that maximize the likelihood via Lagrange multipliers (see the tools used in \cite{chorhendysnir2006} and \cite{chor2000}).

The structure of the paper is as follows. In the following section 2, we recall the background on linear representation theory of finite groups that is needed in the sequel. In section 3, we recall the definition of equivariant models and of phylogenetic varieties. In this section we prove as well two key results that shall allow us to provide a complete intersection as desired: first we compute the dimension of the phylogenetic varieties for any equivariant model $\mathcal{M}_G$, $G\subset \mathfrak{S}_{\kappa}$, and any tree $T$, and then we prove that these varieties are smooth at generic points of no evolution. Then in section 4 we describe the set of equations that shall be used to prove our main result. The setup for this description is conceived towards the induction steps that are needed in the proof of the main result. In section 5 we describe the induction step and the ``claw tree hypotheses'' needed to prove our main result, Theorem \ref{invariantsCI} in the largest generality. The proof of this theorem is constructive and provides an algorithm for obtaining the desired complete intersection assuming the claw tree hypotheses is satisfied. In section 6 we prove that this claw tree hypothesis holds for trivalent trees on the general Markov model, the strand symmetric model, and the Jukes-Cantor model (the Kimura 3-parameter case was already considered in \cite{casfer2008}). For these models we also specify complete intersections (following the algorithm provided in section 5) that describe the variety for quartet trees around generic points of no evolution.

\section*{Acknowledgements}
M.Casanellas and J.Fern\'andez-S\'anchez were partially supported by Spanish government MTM2012-38122-C03-01/FEDER
 and Generalitat de Catalunya 2014 SGR-634. 
M.Micha{\l}ek was supported by National Science Center grant SONATA UMO-2012/05/D/ST1 /01063. 
Part of the work was conducted while Micha{\l}ek was visiting Freie Universit\"at, Berlin and UC Berkeley (PRIME DAAD program 2015-2016) and CRM Barcelona (EPDI program 2013).

\section{Background on representation theory}

%\begin{itemize}
% \item rooted, vertices (interior nodes and leaves) and edges. Denote by $T_n$ a tree with $n$ leaves.
%  \item Compute dimensions of equivariant models (Chang)
% \item Complete intersections
% \end{itemize}

%\paragraph{Representation theory}

In this section we recall the basic concepts of representation theory that will be needed in the sequel. The reader is referred to \cite{serre} or \cite{fulton1991} for details and proofs. Throughout the paper we work over the field of complex numbers $\CC$.

Let $G$ be a finite group. %
A \textit{representation} of $G$ is a group homomorphism $\rho:G\rightarrow GL(V)$, where $V$ is a $\CC$-vector space of finite dimension. We will refer to $V$ as the representation itself (or also as a $G$-module) if the map $\rho$ can be understood from the context, and for $g\in G$ and $u\in V$ we shall denote by $g u$ the vector $\rho(g)(u)$.
%; in this case the homomorphism shall be denoted as $\rho_V$.
A $G$-\textit{equivariant map} is a linear map $f:V\lra V'$ between two representations of $G$ that satisfies $f \circ \rho(g)= \rho_{V'}(g)\circ f$ for all $g\in G$. The set of all $G$-equivariant maps between $V$ and $V'$ is denoted as $\Hom_G(V,V')$. %and coincides with the elements of $Hom(V,V')$ that are invariant under the action of $G$.
Two $G$-modules $V$ and $V'$  are said to be \textit{isomorphic} (denoted as $V\cong_G V'$) if there is a $G$-equivariant isomorphism of vector spaces $f:V\rightarrow V'$.
A representation $V$ is \emph{irreducible} if it does not contain any proper $G$-invariant subspace. Otherwise, $V$ is said to be \emph{reducible}.
We will denote by $V^G$ the subspace of vectors of $V$ that are $G$-\emph{invariant} under the action of $G$, that is, $g u=u$ for all $g\in G$.

\begin{lem}[Schur]
 Let $V,V'$ be two irreducible representations of $G$. If $f:V\rightarrow V'$ is $G$-equivariant, then either $f=0$ or $f$ is an isomorphism, in which case ${\Hom_G(V,V')\cong \CC}$.
\end{lem}

\begin{notation}\label{nota}
 Let $\rho_k:G\rightarrow GL(\N{k})$, $k=1,\ldots,t$ be the irreducible representations of $G$ (up to isomorphism).
%, this meaning the equivalence classes of irreducible representations.
We write
%$d_k$ for the dimension of $N_k$, and
$\chi_k$ for the  character corresponding to $\N{k}$: $\chi_k(\sigma)=\operatorname{trace}(\rho_k(\sigma))$. We adopt the convention that $(\rho_1,\N{1})$ refers to the \emph{identity} (or \emph{trivial}) representation.
%\private{
%
%We Removed from above: ``Moreover, we will assume an \emph{orthonormal} basis of $N_k$ is fixed throughout the paper $B_k=\{w ^k_1,\ldots, w^k_{d_k}\}$, $k=1,\ldots, t$.''
%}
%
\end{notation}

% Schur's lemma implies that if $f:N_k\rightarrow V$ is an isomorphism of $G$-modules, the image of the basis $B_k$ by $f$ is determined up to multiplicative constant, with all the vectors in $B_k$ having the same constant. In particular, we can write $\Hom_G(N_k,V)=[ f ]_{\CC}$.

\begin{thm}[Maschke]\label{Maschke}
 If $\rho:G\rightarrow GL(V)$ is a representation of $G$, then there exists a unique decomposition
 $ V= \oplus_{k=1}^t V[\chi_k]$,   where each $V[\chi_k]$ is isomorphic to $\oplus^{m_k(V)}\N{k}$ for some \emph{multiplicity} $m_k(V)\geq 0$. We call $V[\chi_k]$  the \emph{isotypic component} of $V$ associated to $\N{k}$.
\end{thm}
Notice that $V^G$ is equal to the isotypic component of $V$ associated to the trivial representation of $G$: $V^G=V[\chi_1]$.

\begin{rem}\label{rem_Vk}
By virtue of these fundamental results, for any representation $V$ of $G$, the dimension of $\Hom_G(\N{k}, V)$ equals the multiplicity of $\N{k}$ in $V$. Moreover, the collection of the images of a chosen vector $v_k \in \N{k}$ under maps in $\Hom_G(\N{k}, V)$  form a subspace $\FF{k}{V}$ in $V$ of dimension equal to the multiplicity of the isotypic component, $\FF{k}{V}\cong \CC^{m_k(V)}$.
Analogously to highest weight spaces, the spaces $\FF{k}{V}$ will represent the whole isotypic components. In particular, for any two representations $V, V'$ we can identify $\Hom_G(V, V')$ with $\bigoplus_k\Hom_\CC(\FF{k}{V},\FF{k}{V'})$.
\end{rem}

%Now, if there is an inner product $\langle \cdot, \cdot \rangle$ in $V$, we can identify $V$ with its dual $V^*$.
%\private{One has to be very careful here! This is true on the level of vector spaces (and some representations). However, in general the dual representation is different from the representation itself (it depends on the fact if the character is real). (over C orthonormal basis does not give duality!). However, it works e.g. for $W\simeq W^*$ if our group is a subgroup of symmetric group and $W$ is the permutation representation. }
%\tp{We agree. Indeed, if we decompose $W$ into irreducible, this will involve non-real  characters.}
%

%Denote by $\m=(m_1,m_2,\ldots,m_t)$ the multiplicities of the irreducible representations of $G$ in the Maschke's decomposition of the space $W$.
%\begin{lem}\label{lem:dual}

%\end{lem}

%\private{\tp{I've reordered the text here.}}

\paragraph{Permutation representation.}

From now on we focus on the following setting. Given a finite set $\Sigma$
%=\{\xx_1,\ldots,\xx_{\kappa}\}$}
of cardinality $\kappa$, we define $W=\langle \Sigma \rangle _{\CC}$ as the $\CC$-vector space generated by the elements of $\Sigma$.  In this way, the elements of $\Sigma$ play the role of the standard basis of $W$, so that an element $\xx \in \Sigma$ and the corresponding vector of the standard basis shall be denoted in the same way. Motivated by biology, in our examples we consider $\Sigma=\{\aa,\cc,\gg,\tt\}$ but our work holds for any finite set.
We denote $\mathbf{1}:=\sum_{\xx\in\Sigma} \xx$. %, which will be also taken as the vector with all its coordinates equal to one.
By abuse of notation, $\mathbf{1}$ will be sometimes taken as the column-vector with all its $\kappa$ coordinates equal to one.
Henceforth, $G$ shall be a permutation group of $\Sigma$, that is, $G$ is a subgroup of $\mathfrak{S}_{\kappa}$. The restriction to $G$ of the \textit{permutation representation} $W$, given by the permutation of the elements in $\Sigma$, induces a representation $\rho(s)$ of $G$ on any tensor power $\op{s}$ by extending linearly the action
$ \sigma (\xx_{i_1} \o \ldots \o \xx_{i_s}):=\sigma \xx_{i_1} \o \ldots \o \sigma \xx_{i_s}$ for $\sigma \in G, \xx_{i_j}\in \Sigma$.
In this paper, we will only deal with such representations $\rho(s):G\lra GL(\op{s})$ together with the irreducible representations $\N{1}, \dots,\N{t}$ of $G$.
%
%

% \private{
% \red{We'll try to avoid using the following notation:} We write $\LL_s=\o^s W$, and $\LL_s^G$ for the space of tensors that are invariant by the action of $G$.%
% }
%By the isomorphism (\ref{equiv:tensors_maps}) above, these tensors translate into $G$-equivariant homomorphisms.
%, that will be denoted by  $\Hom_G(\cdot,\cdot)$.

According to Masckhe's theorem, any tensor power $\op{s}$ will decompose into a direct sum of modules $(\op{s})[{\chi_k}]$ (the isotypic components) each of them being a number of copies of one of the irreducible modules $\N{k}$. This number is the multiplicity of $\N{k}$ in $\o^{s}W$ and will be denoted by $m_k(s)$. In the particular case $G=\mathfrak{S}_k$, explicit formulas for $m_k(s)$ can be provided in terms of Kronecker coefficients.
We write $\m(s)=(m_1(s),\ldots,m_t(s))$ for the vector of multiplicities of $\o^s W$. As the case $s=1$ will play a special role, we simplify notation and write $\m=(m_1,\ldots, m_t)$ for the vector of multiplicities of $W$.

From now on, we fix subspaces $\FF{k}{W} \subset W$ for $k=1,\dots t$ according to Remark \ref{rem_Vk} by fixing a vector $v_k\in \N{k}$ and taking its images by maps in $\Hom_G(\N{k}, W)$. This vector also defines subspaces $\FF{k}{\o^{s}{W}}   \subset \o^{s}W$, which shall be considered fixed from now on.

We consider the Hermitian inner product in $W$ that makes $\Sigma$ into an orthonormal basis, and denote it by $v\cdot w$ for any $v,w\in W$. %
This inner product will be used to identify $W$ with $W^*$ by sending a vector $v$ to the linear form $v^* \in W^*$ that maps $u$ to $v\cdot u$.

The inner product in $W$ induces an inner product in $\op{s}$  defined as
\begin{eqnarray*}
 \xx_{1} \o \ldots \o \xx_{n} \cdot \mathtt{Y}_{1}\o \ldots \o \mathtt{Y}_{n} :=\prod_{k=1}^s  \xx_{k} \cdot\mathtt{Y}_{k},
\end{eqnarray*}
for $\xx_{i},\mathtt{Y}_{j}\in \Sigma$ and extending it sesquilinearly.

\paragraph{Dual representation.}
If $\rho: G\lra GL(V)$ is a representation of $G$ with character $\chi$, then its dual $V^*$ is also a representation via the homomorphism
$\rho^* : G \lra GL(V^*)$ that maps $g$ to $^t \rho(g^{-1})$, and $V^{*}$ has  character $\chi^*$ (the conjugate of $\chi$).
%
% Moreover, there is an isomorphism of representations $V\cong_G V^*$ induced by mapping each element $X\in \Sigma$ to the corresponding form $\omega_X$ given by $\omega_X(u)=\langle u, X \rangle$.
% %
% For any representation $V$ we have natural isomorphisms of representations $V^*[\chi_k]\simeq_G (V[\chi_{k^*}])^*$ and of linear spaces $(V^*)_k\simeq (V_{k^*})^*$, where $k^*$ is the index of the (irreducible) representation dual to $N_k$.

%Note that the identification of $V$ and $V^*$ is on the level of linear spaces.
At the level of vector spaces, the inner product above provides an isomorphism  $V \cong V^*$. Nevertheless, if $V$ is a representation of a group $G$, then it may happen that it is not $G$-isomorphic to $V^*$. This will force us to distinguish between the space and its dual in the sequel. However, the permutation representations $V=\op{s}$ that we consider in this paper satisfy $V\cong_G V^*$ because they have real characters.
In particular, the $G$-isomorphism $V\cong_G V^*$  induces $G$-isomorphisms  $V^*[\chi_k]\cong_G  V[\chi_k]$ and $\FF{k}{V}\cong \FF{k}{{V}^{*}}$ for all $k$.

Thus, the reader may freely ignore all the dual signs in our article. We decided to keep them, as most of the arguments we provide hold without the assumption $V\cong_G V^*$ on a representation theoretic level. There is also a natural $G$-isomorphism $V \o V' \cong_G  \Hom(V^*,V')$ which at the level of $G$-invariant vectors translates to $\left( V \o V'\right)^G \cong_G \Hom_G(V^*,V')$.
%\tp{mm: I would vote for removing last sentence. MARTA: Why? mm: I mean at least the isomorphisms $\cong_G \Hom(V,V')$ and $\cong_G \Hom_G (V,V')$ - this I think does not look good.}

Representation theory allows us to decompose the ambient space $(\o^n W)^G$ in terms of the irreducible representations of $G$ as follows. This decomposition will be fundamental for us and will play a key role in the paper.

%\tp{mm: just a small notation problem - maybe we do not have to even care about it. For any representation $V$ we have the space $V_k$. On the other hand $\N{k}$ is the representation, i.e. lower index $k$ plays two roles, once as a functor, once as an index}

\begin{prop}\label{prop:W^nG}
%Let $A|B$ be a bipartition of $L(T)$. Then,
For any %pair $(a,b)$ with
$a+b=n$, there is a natural isomorphism of vector spaces
 \begin{eqnarray*}
 \left (\otimes^{a+b} W \right )^G
 \cong \bigoplus_{k=1}^t \;  \FF{k^*}{\otimes^a W}  \o  \FF{k}{\otimes^b W} .
 \end{eqnarray*}
 In particular, the dimension of $(\o^{a+b} W)^G$ is $m_1(a+b)=\sum_{k=1}^t m_{k^*}(a)m_{k}(b)$, where $k^*$ is the index of the irreducible representation dual to $\N{k}$, that is, $\N{k^*}=(\N{k})^*$.
  Using the language of category theory, the functors $(\otimes^n \cdot)^G$ and $\bigotimes_{k=1}^t\FF{k^*}{\otimes^a\cdot}\otimes\FF{k}{\otimes^b \cdot}$ from the category of $G$ representations to the category of vector spaces are isomorphic.
\end{prop}

\begin{proof}
% A consequence of Maschke's theorem and Schur's lemma.
%We keep the notation introduced in Theorem \ref{Maschke}.
%For each irreducible representation $\N{k}$ of $G$, we have isomorphisms
%\begin{eqnarray*}
%(\o^a W)[\chi_k]\cong \oplus^{m_k(a)}\N{k} \qquad \mbox{ and }\qquad
%(\o^b W)[\chi_k] \cong \oplus^{m_k(b)}\N{k}.
%\end{eqnarray*}
%\red{As the representations $W$ and $W^*$ are isomorphic, }
Applying Maschke's theorem and Schur's lemma, we infer
\begin{eqnarray*}%\label{tensors_maps}
\opG{n} & \cong & \left ( (\o^a W) \o (\o^b W)\right )^G \cong \Hom_G((\o^a W)^*,\o^b W) \\
& \cong & \oplus_{i,j} \Hom_G((\o^a W)^*[\chi_i],(\o^b W)[\chi_j]) \cong  \oplus_{k=1}^t \Hom_{\CC}( \FF{k}{\o^a W^*}, \FF{k}{\o^b W}) \\ & \cong & \oplus_{k=1}^t \Hom_{\CC}\left( (\FF{k^*}{\o^a W})^*, \FF{k}{\o^b W}\right) \cong \bigoplus_{k=1}^t \; \FF{k^*}{\otimes^a W}  \o \FF{k}{\otimes^b W}.
\end{eqnarray*}
%the last isomorphism because of Schur's lemma.
\end{proof}

%
%
% \begin{rem}
% Recall that by the Schur-Weyl duality
% $$W^{\otimes n}=\sum_{\lambda\vdash n} S_\lambda(W)\otimes [\lambda]$$
% as a $GL(W)\times \mathfrak{S}_n$ representation. In particular, no evolution points are symmetric tensors invariant under the action of $G\subset GL(W)$.
%
%  Relate this to Schur-Weyl duality. Non-evolutionary points are points invariant under the action of $G$ and also the action of $\mathfrak{S}_n$, that permutes the leaves of the tree.
% \end{rem}

%\private{ Is $W_A^*[\chi_i]\cong_G  W_A[\chi_i]$ for any $i$? }
%In particular, we always have $V\o V'\simeq \Hom_G(V^*,V')$.

% Let $A|B$ be a bipartition of $L(T)$. Then,
% \begin{enumerate}
%  \item (Maschke)
%  $W_A=\oplus_{k=1}^t W_A[\chi_k]$,
%   and $W_A[\chi_k]\cong \bigoplus^{m_k(a)} N_k$. Similarly, for $B$. We may identify $\CC^{m_k(a)}$ with $(W_A)_k$.
% \item $\dim_{\CC}(W_A[\chi_k])=d_k m_k(a)$. Similarly, for $B$.
% %\item Schur's lemma: $\Hom_G(W_A^{\chi_k},W_B^{\chi_k})\cong \Hom_{\CC}(\CC^{m_k(a)},\CC^{m_k(b)})$
% \end{enumerate}

%\private{\tp{proposition moved to next section}}

\section{Equivariant evolutionary models and phylogenetic varieties}

A tree is a connected finite graph without cycles, consisting of vertices and edges. Given a tree $T$, we write $V(T)$ and $E(T)$ for the set of vertices and edges of $T$. The\textit{ degree} of a vertex is the number of edges incident to it. The set $V(T)$ splits into the set of leaves $L(T)$ (vertices of degree one) and the set of interior vertices $Int(T)$ : $V(T) =L(T) \cup Int(T)$. One says that a tree is \textit{trivalent} if each vertex in $Int(T)$ has degree 3. A tree topology is the topological class of a tree where every leaf has been labeled. Given a subset $A$ of $L(T)$, the subtree induced by $A$ is just the smallest tree composed of the edges and vertices of $T$ in any path connecting two leaves in $A$.
A tree $T$ is \textit{rooted} if a specific node $r$ is labeled as the root.

In order to model the substitution of the states in $\Sigma$ according to a Markov process on a rooted tree $T$, one has to  specify a distribution $\pi$ at the root of the tree and a collection of substitution matrices $\mathbf{A}=(A^e)_{e\in E(T)}$ \cite{chang1996,casferked}.
The set of possible root distributions and substitution matrices for a tree $T$ is called the \textit{set of parameters}. In the applications to biology, one has to restrict the set of parameters to stochastic vectors and matrices, but this restriction is unnecessary for the core of this paper. Below we describe \textit{equivariant} models of evolution, which include some of the most well-known models.

As above, let $\Sigma$ be a finite set of cardinal $\kappa$, $W$ be the $\CC$-vector space $\langle \Sigma \rangle_{\CC}$, and $G\leq \mathfrak{S}_{\kappa}$  be a permutation group of $\Sigma$.
%For biological reasons, $\Sigma$ is usually taken as the set $\{\aa,\cc,\gg,\tt\}.$ %\tp{This has already been said in page 2: is it necessary to insist? }
In this section we use the distinguished basis $\Sigma$ of $W$ to identify $\kappa\times\kappa$ matrices with complex entries with $\Hom(W,W)$.
%\tp{We should decide once and for all if matrices are elements of $\Hom(W^*,W)$ or $\Hom(W,W)$. See next definition.}.

\begin{defn}{(cf. \cite{draisma2008})}
A rooted tree $T$ evolves under the \textit{equivariant model} $\MM_G$
if a $G$-invariant vector $\pi$ is associated to the root of $T$ and substitution matrices $A^e$ in $\Hom_G(W,W)$ are associated to each edge $e$ of $T$.
%if the root distribution $\pi$ is a $G$-invariant vector in $W$ and the substitution matrices $A^e$ belong to $\Hom_G(W,W)$.
For the equivariant model $\MM_G$, the set of \textit{parameters} is
\[\Par_G(T)=W^G \times \prod_{e\in E(T)}\Hom_G(W,W).\]
If one wants to talk about \textit{stochastic parameters}, $s\Par_G(T)$ one has to restrict the root distribution to $s W^G:= W^G \cap \{\pi\in W \mid \pi \cdot \mathbf{1} =1\}$, and the substitution matrices to $s\Hom_G(W,W):= \Hom_G(W,W) \cap \{A \mid A \cdot \mathbf{1}=\mathbf{1} \}$ (and then require that all entries are real, nonnegative, but this is not relevant for our purposes). %
As a special case, if the group $G$ is trivial, $G=\{1\}$, we will denote by $\Par(T)$ and $s\Par(T)$ the corresponding spaces of parameters.
The \textit{parametrization map} that assigns a distribution at the leaves of $T$ to each set of parameters is
\begin{eqnarray}\label{param}
\Psi_T: \, \Par(T) \lra \o^n W
\end{eqnarray}
% \begin{eqnarray}\label{param}
% \Psi_T: \, W \times \prod_{e\in E(T)} \Hom_{\CC}(W,W) \lra \LL_n
% \end{eqnarray}
defined by
\begin{eqnarray*}
\Psi_T\left(\pi, \mathbf{A}\right)=\sum_{\xx_i \in \Sigma}p_{\xx_1...\xx_n}\xx_1 \otimes \dots \otimes \xx_n,
\end{eqnarray*}
where
\begin{eqnarray}\label{form_param}
p_{\xx_1\dots \xx_n}=\sum_{\xx_v \in \Sigma, v \in Int(T)}
\pi_{\xx_r}\prod_{e\in E(T)}
A^{e}_{\xx_{pa(e)}, \xx_{ch(e)}} \, ,
\end{eqnarray}
$\xx_v$ denotes the state at the vertex $v$,
$pa(e)$ (respectively $ch(e)$) is the parent (respectively, child) node of $e$, and $(\pi_{\xx})_{\xx \in \Sigma}$ are the coordinates of the root distribution $\pi$.
When we restrict this map to the set of parameters $\Par_G(T)$, we denote it as $\Psi_T^{G}.$ In this case the image lies in $(\o^n W)^G$.

When the parametrization  (\ref{param}) is restricted to the set of stochastic parameters, we obtain
\[\phi_T:  s\Par(T) \longrightarrow H\cap \o^n W,\]
where $H \subset \o^n W$  is the hyperplane defined as
\[H=\left\{ p \in \o^n W \mid  \sum_{\xx_i\in \Sigma} p_{\xx_1 \dots \xx_n}
=1 \right\}.\]
The analogous restrictions to $s\Par_G(T)$ are denoted as $\phi_T^{G}$.
The word ``stochastic'' here has a broader meaning than usually, because for our aim  we only need entries summing to one and not necessarily nonnegative entries.

The \emph{phylogenetic variety associated to a tree $T$ evolving under $\MM_G$} is the (affine) algebraic variety
 \begin{eqnarray*}
  CV_G(T):=\overline{\{\Psi_T^G(\pi,\mathbf{A}) \mid (\pi,\mathbf{A})\in \Par_G(T)\}} \subset (\o^n W)^G.
 \end{eqnarray*}
 where $\overline{S}$ represents the Zariski closure of a set $S$.
%An \emph{equivariant model} $\MM_G$ of evolution is a pair $(G,W)$ where $W=[ \Sigma ]_{\CC}$, $G\leq \mathfrak{S}_n$. Trees evolving under this equivariant model are phylogenetic trees on $(G,W)$ together with a space of $G$-evolutionary presentations. \tp{missing def. of presentations.... we'll add it.}
Similarly, the \emph{stochastic phylogenetic variety
associated to a tree }$T$  is the smallest algebraic variety $V_{G}(T)$ containing the set $$\Im \, \phi_T^{G}=\left\{ \phi_T^{G} (\pi,\mathbf{A}) : (\pi,\mathbf{A})\in  s\Par_{G}(T) \right\}.$$ One has
$V_{G}(T)= CV_G(T)\cap H$ (see, for example, \cite{casferked}). In particular, the equations defining $V_{G}(T)$ are the same equations defining $CV_G(T)$ plus the equation defining $H$.
\end{defn}

Notice that in the definition of the phylogenetic variety we have not specified the root of the tree. It is well known that different root placements give rise to the same phylogenetic variety. Indeed, it is clear that a matrix $M$ belongs to $\Hom_G(W,W)$ if and only if so does its transpose $M^t$. Now, it can be seen that if we move the root from one node to a neighboring node and we replace the matrices $A^e$ of the edges with inverted orientation with their transpose, the image of the new parameters will remain the same. Moreover, we may assume that $\pi=\mathbf{1}$ when parameterizing $CV_G(T)$ since choosing an edge $e_0$ attached to the root and changing $A^{e_0}$ by $diag(\pi)A^{e_0}$ gives rise to the same image point. Hence $CV_G(T)$ does not depend on the root. As $H$ also does not depend on the root, neither does $V_G(T)$.
%\private{Need to talk about root independence.}

%\private{\red{Marta changed non-evolutionary by no evolution everywhere. It is more used in Google!} \tp{I would suggest %"points of no evolution" or "points of no mutation", instead.}}

In case we take all matrices $A^e$ equal to the identity, the image by $\Psi_T^G$ represents no evolution at all.

\begin{defn}
Given an equivariant model $\MM_G$, a
point $\pi_n$ in $\o^n W$ is a \emph{point of no evolution} if
$\pi_n=\sum_{\xx \in \Sigma} \pi_{\xx} \xx \o \ldots \o  \xx$ and $\pi_n$ is %where $\pi:=\pi_A \aa+ \pi_C \cc+\pi_G \gg +\pi_T \tt$ is
invariant under $G$. %\red{I do not get it. Is there a misprint? Is $\pi$ the $X$? What is $\pi_X$?}
%and $(p_{b_1 \ldots b_1},\ldots,p_{b_n \ldots b_n})$ is $G$-invariant.
\end{defn}

If $\pi_n=\sum_{\xx \in \Sigma} \pi_{\xx} \xx \o \ldots \o  \xx$ is a point of no evolution, then it belongs to $CV_G(T)$ for any tree $T$ on $n$ leaves because $\pi_n=\Psi_T^G(\pi,\mathbf{I})$ where $\pi=\sum_{\xx \in \Sigma}\pi_{\xx} \xx \in W^G$ and $\mathbf{I}$ corresponds to the identity matrix at each edge: $\mathbf{I}=(Id)_{e\in E(T)}$.  %Moreover, we observe that if $\pi_n \in \op{n}$ is a point of no evolution then $\pi_n$ lies in $\op{n}^G$ and therefore in $\op{n}[\chi_1]$.

%\tp{After the lines in red in the previous paragraph, I don't see the point of the sentence in black in the same paragraph. I mean, if we know that $\pi_n$ lies in the variety associated to any tree, of course it lies in the space of invariant tensors! Hence, I would remove it. }
\begin{remark0}\rm
For biological applications, we are interested in real/stochastic points in $\o^n W$ that are close  to points of no evolution (in the complex euclidean distance). Indeed, as $\Psi_T^G$ is a continuous map, if $p$ is close to a point of no evolution $\pi_n$, then there is a preimage of $p$ close to $(\pi,\mathbf{I})$. The parameters close to $(\pi,\mathbf{I})$ are precisely those that are interesting biologically speaking because they account for probabilities of no mutation greater than probabilities of mutation (that is, diagonal entries greater than off-diagonal entries in transition matrices). To this end, the main goal of this paper is to provide a local description of the phylogenetic varieties around the points of no evolution.
\end{remark0}

\begin{example}\label{main_equivariantmodels}
 The definition of equivariant model includes important evolutionary models used in phylogenetics for $\kappa=4$ like
 \begin{enumerate}
  \item Jukes-Cantor \cite{JC69}, for $G=\mathfrak{S}_4$ or $G=\mathfrak{A}_4$ (the alternating group);
  \item Kimura 2-parameter \cite{Kimura1980}, for $G=\langle (\aa\cc\gg\tt),(\aa\gg)\rangle$;
  \item Kimura 3-parameter \cite{Kimura1981}, for $G=\langle (\aa\cc)(\gg\tt),(\aa\gg)(\cc\tt)\rangle $;
  \item strand-symmetric \cite{cassull}, for $G=\langle (\aa\tt)(\cc\gg)\rangle $;
  \item general Markov (briefly GMM)  \cite{barry1987}, for $G=\{1\}$.
   \end{enumerate}
 We say that $\MM_G$ is a \textit{submodel} of $\MM_H$ if $H \leq G$. With this terminology, all the models above are submodels of the general Markov model and we have inclusions from top to bottom on the sets of corresponding parameters (and algebraic varieties).
%Notice that under our definition of equivariant models, we require the root distribution to be an invariant vector under $G$. In particular, models like F81 or HKY85 do not fit in our definition (cf. \cite{draisma2008}).
\end{example}

\subsection{Dimension of phylogenetic varieties for equivariant models}

In this subsection we compute the dimension of the phylogenetic variety associated to any $G$-equivariant model on any tree $T$, $G\leq \mathfrak{S}_{\kappa}$. This dimension was already known in the particular cases of the Jukes-Cantor, Kimura 2 and 3 parameters and general Markov model. The result yields the codimension of these varieties and hence it is the first step towards providing a complete intersection containing them.

\begin{thm}\label{dimension} For any group $G\leq \mathfrak{S}_{\kappa}$ and any tree $T$ without nodes of degree 2,
the dimension of $V_G(T)$ is $|E(T)|(m_1(2)-m_1)+m_1-1$ and the dimension of $CV_G(T)$ is $|E(T)|(m_1(2)-m_1)+m_1$.
\end{thm}

%\private{\tp{Due to Chang's result \cite{chang1996}, the other bound can also be proven for any equivariant model on $4$ states. Chang's result implies finite fibers for generic parameters for GMM. These genericity (non-zero determinants of transition matrices and non-zero coordinates at the root distribution) is also satisfied for any equivariant model, and therefore generic fibers are also finite for any equivariant model. We'll write this carefully.}}

\begin{proof}
The dimension of $V_G(T)$ is upper bounded by the dimension of $s\Par_G(T)$, which we compute in the following.\\
%Because the tree $T$ is trivalent, the number of edges is $2n-3$.
Each transition matrix $M$ is an element of $$\Hom_G(W,W)\cong (W^*\otimes W)^G,$$ so, as in our case $W\cong W^*$, the number of parameters is $m_1(2)$. However, because of the stochastic assumption, the sum of the rows of each matrix $M$ is fixed to one. Notice that $(W^*\otimes W)^G$ surjects onto $W^G$ by the map $M\mapsto M \mathbf{1}$. Hence, there are $m_1$ independent restrictions on the parameters of $M$. This makes $|E(T)|(m_1(2)-m_1)$ free parameters for the choice of the transition matrices.
\\
On the other hand, the distribution of the root is given by a vector $\pi\in W^G$. The stochastic condition implies that the sum of the coordinates is equal to one. This makes $m_1-1$ free parameters for the choice of the root distribution.
\\
Summing up, we have that $\dim V_G(T)$ is less or equal than $|E(T)|(m_1(2)-m_1)+(m_1-1)$.

In order to prove the other inequality we  use Chang's result (\cite{chang1996}) and its generalization (\cite{allman2003}) on the ``generic identifiability of parameters`` of the general Markov model $\MM_1$ on trees without nodes of degree 2. This result says that the fiber of $\phi_T(\mathcal{P})$ is finite for parameters $\mathcal{P}=\left(\pi,(A^e)_e\right)$ that satisfy:
(1) no entry of $\pi$ is zero;  (2) all $A^e$ are non-singular;  and
 (3) $\det A^e \neq \pm 1$ for all $e$.
% \begin{enumerate}
%  \item[(1)] no entry of $\pi$ is zero,
%  \item[(2)] all $A^e$ are non-singular, and
%  \item[(3)] $\det A^e \neq \pm 1$ $\forall e,$
% \end{enumerate}
 These generic conditions $(1)$-$(3)$ are also generic for the parameters of any equivariant model $\MM_G$. That is, if $\Sigma=\{\xx_1,\dots,\xx_{\kappa}\}$, for any group $G \leq  \mathfrak{S}_{\kappa}$ we have
\begin{enumerate}
 \item[(i)] $sW^G$ is not included in $\{\pi\in W \,\mid  \pi_{\xx_1}\cdot \ldots \cdot \pi_{\xx_{\kappa}}=0 \}$ (indeed, $\frac{1}{\kappa}\mathbf{1} \in sW^G$ for example),
 \item[(ii)] $s\Hom_G(W,W)$ is not contained in the set of singular matrices (indeed, $Id \in s\Hom_G(W,W)$), and
 \item[(iii)] $s\Hom_G(W,W)$ does not only contain matrices with determinant $1$ or $-1$ (for example, the matrix with all entries equal to $1/{\kappa}$ belongs to $s\Hom_G(W,W)$).
\end{enumerate}
This means that for generic parameters $\mathcal{P} \in s\Par_G(T)$, if we set $p=\phi_T^{G}(\mathcal{P})$,
%(which coincides with $\phi_T^{id}(\mathcal{P})$)
the preimage $(\phi_T)^{-1}(p)$ is finite. As the preimage $(\phi_T)^{-1}(p)$ contains $(\phi_T^{G})^{-1}(p)$, this implies that the generic fiber of $\phi_T^{G}$ is finite.
Therefore, the dimension of $V_G(T)$ is upper bounded by the dimension of the domain of $\phi_T^{G}$, which has been computed above.

The claim for $CV_G(T)$ follows because it is the closure of the cone over $V_G(T)$.
\end{proof}

This result implies the generic identifiability of the stochastic parameters for trees evolving under equivariant models (see \cite[Definition 1]{Allman2008a} for example). %This was already known for Jukes-Cantor, Kimura 2  and 3 parameters, and the general Markov model.
\begin{cor}
 The stochastic parameters of a tree $T$ evolving under an equivariant model $\mathcal{M}_G$, $G\leq \mathfrak{S}_{\kappa}$, are generically identifiable if $T$ has no nodes of degree 2.
\end{cor}

\begin{rem}
It can be checked easily that for the evolutionary models listed in Example \ref{main_equivariantmodels}, the
dimension for a trivalent tree on $n$ leaves (and hence with $2n-3$ edges)
%upper bound of the previous proposition becomes an equality for trivalent trees. Namely, if
$T$ is
\begin{enumerate}
 \item $\dim_{\CC} V_G(T)=2n-3$ for the Jukes-Cantor model;
 \item $\dim_{\CC} V_G(T)=4n-6$ for the Kimura 2-parameter model;
 \item $\dim_{\CC} V_G(T)=6n-9$ for the Kimura 3-parameter model;
 \item $\dim_{\CC} V_G(T)=12n-17$ for the strand symmetric model;
 \item $\dim_{\CC} V_G(T)=24n-33$ for the general Markov model.
\end{enumerate}

\end{rem}

If $T$ has $n$ leaves, we will write $\codim(T)$ for the codimension of $CV_G(T)$ in  $(\o^n W)^G$ (equal to the codimension of $V_G(T)$ in $H$), that is,
$\codim(T):=\dim (\o^n W)^G-\dim CV_G(T)$. The dimension of $(\o^n W)^G$ is $m_1(n)$, which has been computed in \cite[Prop. 20]{casferked} for the models listed in Example \ref{main_equivariantmodels}.

%Notice that for trivalent trees this codimension does not really depend on $T$ but on the number of leaves, so we will often denote it by $\codim(n)$. As trivalent trees have $2n-3$ edges, we have
%\begin{equation}\label{codim}
%\codim(n)=m_1(n)- (2n-3)(m_1(2)-m_1)-m_1.
%\end{equation}

% \begin{example}
%  JC, GM, SM? Non-evolutionary points in ``Fourier'' coordinates.
% \end{example}

\subsection{Smoothness at points of no evolution}\label{subsec:smooth}

Let $T_n$ be the \emph{claw $n$-tree}, that is, the tree with one inner vertex and $n$ leaves. In what follows we prove that the variety corresponding to $T_n$ is smooth at generic points of no evolution. In particular, it can be locally defined by a complete intersection.

Given a permutation subgroup $G$ of $\mathfrak{S}_{\kappa}$, we denote by $GL(\kappa)^G$ the group of $G$-equivariant $\kappa \times \kappa$ invertible matrices. %,  $GL(\kappa)^G = \{A \in GL(\kappa) | P_g A =A P_g \, , \forall g\in G\}$.  
Clearly, $GL(\kappa)^G$ defines an action on $\Hom_G(W,W)$ by $(A,M) \rightarrow AM$.

\begin{thm}\label{thm:orbit}
The variety $CV_G(T_n)$ is the Zariski closure of the orbit of $\Psi_T^G(\mathbf{1},\mathbf{I})$ under the group action of $(GL(\kappa)^G)^n$.
\end{thm}
\begin{proof}
The conditions
\begin{itemize}
\item[1.] $\det A^e\neq 0$ for all $e\in E(T)$, and
\item[2.] all coordinates $\pi$ are different from $0$
\end{itemize}
define open sets in $\Par_G(T)$. As $\mathbf{I}\in GL(\kappa)^G$ and $\mathbf{1}\in W^G$, the intersection of these open sets is non-empty and a generic point $(\pi, A^e)\in\Par_G(T)$ satisfies both conditions.
%A generic point of $CV_G(T_n)$ is $\Psi_T^G(\pi,A^e)$.
Let us fix one edge $e_0$. Notice that $diag(\pi)A^{e_0}$ is invertible (as the coordinates of $\pi$ are nonzero)
and is $G$-invariant (as $\pi\in V^G$). This means $diag(\pi)A^{e_0}\in GL(\kappa)^G$. Notice that $\Psi_T^G(\pi,A^e)=\Psi_T^G(\mathbf{1},\tilde A^e)$, where
$\tilde A^e=A^e$ for $e\neq e_0$ and $\tilde A^{e_0}=diag(\pi)A^{e_0}$. However, $\Psi_T^G(\mathbf{1},\tilde A^e)=(\tilde A^e)\cdot\Psi_T^G(\mathbf{1},\mathbf{I})$.
\end{proof}
% \begin{cor}
% The varieties $CV_G(T)$ and $V_G(T)$ do not depend on the choice of the root.
% \end{cor}
% \begin{proof}
% As in the proof of Theorem \ref{thm:orbit} we see that in order to parameterize $CV_G(T)$ we may assume $\pi=\mathbf{1}$ putting the root distribution into
% one transition matrix. Hence $CV_G(T)$ does not depend on the root. As $H$ also does not depend on the root, neither does $V_G(T)$.
% \end{proof}
\begin{cor}
If $\sum_{\xx \in \Sigma} \pi_{\xx} \xx \o \ldots \o  \xx$ satisfies $\pi_{\xx} \neq  0$ for all $\xx$, then it is nonsingular. In particular, a generic point of no evolution of $CV_G(T_n)$ and $V_G(T_n)$ is nonsingular.
\end{cor}
\begin{proof}
For $CV_G(T_n)$ the statement follows directly from Theorem \ref{thm:orbit}. Let us also notice that $CV_G(T_n)$ is a cone, hence a point of $V_G(T_n)$ is smooth
if and only if it is smooth as a point of $CV_G(T_n)$.
\end{proof}

%\tp{I've removed a remark here; the claim has been included in the statement of the previous corollary}

% \begin{rem}\label{rem_generic}
% Above a point of no evolution $\sum_{\xx \in \Sigma} \pi_{\xx} \xx \o \ldots \o  \xx$ is \textit{generic} if
% $\pi_{\xx} \neq  0$ for all $\xx$.
% \end{rem}

\begin{rem}
When $G$ acts on the basis of $V$ transitively and freely then $GL(V)^G$ is a torus. This is the case of so-called group-based models and, as follows from Theorem \ref{thm:orbit} the variety $V_G(T_n)$ is toric \cite{sturmfels2004}.
\end{rem}

\section{Equations for the complete intersection}

A \emph{bipartition} $A|B$ of the set of leaves is just a decomposition $L(T)=A\cup B$, where $A$ and $B$ are disjoint sets. Throughout the paper, we write $a=|A|$, $b=|B|$ and, to avoid trivialities, we will assume that $a,b \geq 2$.  A bipartition is an \emph{edge split of $T$} if it arises by removing one of the edges of $T$.
If $A\subset L(T)$, we denote $\o_{i\in A} W_i$ by  $W_A$.

%\private{We have changed old notation : $T_1$,$T_2$ by $T_A$, $T_B$ and $L_1$, $L_2$ by $\alpha$, $\beta$. }

\begin{figure}
\begin{center}
\includegraphics[scale=0.25]{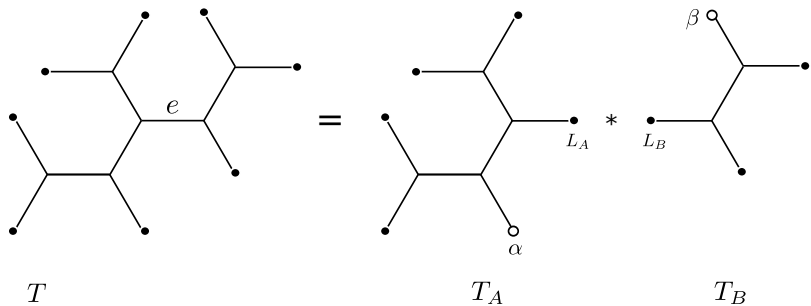}
\end{center}
\caption{\label{figure:decomposition2} Decomposition of $T$ into two subtrees $T_A$, $T_B$: $T=T_A \ast T_B$.}
\end{figure}

Given a tree $T$, in this section we proceed to construct equations that will define a complete intersection for $CV_G(T)$. We choose an internal edge $e$ of $T$, which induces an edge split of the set of leaves $L(T)=A\cup B$. This allows us to view the tree $T$ as the gluing of two trees $T=T_A*T_B$ where $L(T_A)=A\cup L_A$, $L(T_B)=B\cup L_B$, and $L_A$,$L_B$ are the two vertices of the edge $e$, see Figure \ref{figure:decomposition2} and \cite{allman2004}.
We assume that the leaves of $T$ are ordered so that those in $A$ appear in the first place, and those in $B$ appear afterwards. We call $\alpha\in L(T_A)$ the last leaf of $A$ and $\beta\in L(T_B)$ the first leaf of $B$. %

% Throughout this section the vector $\mathbf{1}$ plays the role of a $G$-invariant vector $\sum_{X\in\Sigma} X$ and the results of this section (and the following) still hold if we use any other $G$-invariant vector.

A complete intersection for the variety $CV_G(T)$ will be obtained by joining  equations of a complete intersection for $CV_G(T_A)$, of a complete intersection for $CV_G(T_B)$ and specific edge invariants. %
All results of this section (and the following) still hold if we replace the vector $\mathbf{1}$ by any other $G$-invariant vector of $W$.

\subsection{A basis linked to an edge split}\label{sec:basis}

In order to provide specific equations for the varieties associated to phylogenetic trees, we proceed to construct a basis $\mathcal{B}_{A|B}$ of $(\o^n W)^G$ related to a given edge split $A|B$ as above.
This basis shall be used to specify coordinates and provide the equations as polynomials in these specific indeterminates.
The important point is that this basis must be consistent with  the  decomposition
\begin{eqnarray}\label{dec}
 (\o^n W)^G %& \cong & \Hom_G(W_A^*,W_B) \cong \oplus_{k=1}^t \Hom_G(W_A^*[\chi_k], W_B[\chi_k]) \cong \\ & \cong &
 \cong \oplus_{k=1}^t \Hom_{\CC}(\FF{k}{W_A^*},\FF{k}{W_B}).
\end{eqnarray}
given by Proposition \ref{prop:W^nG}.

We construct the desired basis of $(\o^n W)^G$ compatible with \eqref{dec} as follows.

\vspace*{2mm}

\noindent
\textbf{Algorithm to construct a basis of $(\o^n W)^G$ linked to  an edge split.}
\begin{enumerate}
  \item Choose bases {$\{u^k_i\}_{i=1\div m_k}$} of each $\FF{k}{W}$, $k=1,\ldots,t$.
 %\item As $W_{\alpha}$, $W_{\beta}$ are copies of $W$, we consider the natural basis in $(W_{\alpha}\otimes W_{\beta})^G$ induced by the bases above via the isomorphism $(W_{\alpha}\otimes W_{\beta})^G\cong \bigoplus_k (W_{\alpha})_{k^*}\otimes (W_{\beta})_k$.
 \item For each $k=1,\dots, t$, the vectors $u^k_{B,i}:=u^k_i \otimes \mathbf{1}^{b-1}$ of $\FF{k}{W_B}$, $i=1\div m_k$, are linearly independent.
 Indeed, the monomorphism $W_{\beta} \stackrel{\cdot \o \mathbf{1}^{b-1}}
 {\longrightarrow} W_B$ obtained by tensoring with a power of $\mathbf{1}$ induces a monomorphism
 $\FF{k}{W} \stackrel{\iota}{\hookrightarrow} \FF{k}{W_B}$.  We extend them to a basis $\{u^k_{B,i}\}_{i=1\div m_k(b)}$ of $\FF{k}{W_B}$.
% by vectors in the orthogonal to $\Im(\iota)$ %(that is, vectors in the kernel of the contraction map $f_{\mathbf{1}}:(W_B)_k\rightarrow W_k$).\tp{mm: was $f_{\mathbf{1}}$ already defined?}
 \item We repeat step 2 for $A$ to obtain a basis $\{u^{k^*}_{A,i}\}_{i=1\div m_k(a)}$ of  $\FF{k^*}{W_A}$ for each $k$ (but now tensoring at the left, $\mathbf{1}^{a-1}\otimes u^{k^*}_i $).
 %\item We apply the symmetrization operator
%\begin{eqnarray}\label{S_operator}
%S=\frac{1}{|G|} \sum_{g\in G} \rho_a(g)\o \rho_b(g)
%\end{eqnarray}
%to the basis $\{u^k_{A,i}\o u^k_{B,j}\}_{i,j,k}$, to produce the desired basis $\Ba_{A|B}=\{S(u^k_{A,i}\o u^k_{B,j})\}_{i,j,k}$ of $(\o^n W)^G$.

%  \item Let $S$ be the natural isomorphism of Proposition \ref{prop:W^nG} from \[\FF{k^*}{W_A} \otimes \FF{k}{W_B}\cong \Hom_\mathbb{C}(\FF{k^*}{W_A}^*,\FF{k}{W_B})\] to \[\Hom_G(W_A^*[\chi_k],W_B[\chi_k])\cong (W_A[\chi_{k^*}] \otimes W_B[\chi_k])^G.\]%

\item Write $S$ for the inverse of the isomorphism of Proposition \ref{prop:W^nG}. Its restrictions induce natural isomorphisms from \[\FF{k^*}{W_A} \otimes \FF{k}{W_B}\cong \Hom_\mathbb{C}(\FF{k^*}{W_A}^*,\FF{k}{W_B})\] to \[\Hom_G(W_A^*[\chi_k],W_B[\chi_k])\cong (W_A[\chi_{k^*}] \otimes W_B[\chi_k])^G.\]%

 We call $\Ba_{A|B}$ the desired basis $\{S(u^{k^*}_{A,i}\o u^k_{B,j})\}_{i,j,k}$ of $(\o^n W)^G$.
 \end{enumerate}
From now on, we will denote by $q^k_{ij}$ the coordinate corresponding to the basis vector $S(u^{k^*}_{A,i}\o u^k_{B,j})$.

\begin{rem}\label{S_operator}
%\red{If the chosen basis $\{u^k_{A,i}\}$ satisfies $\rho_a(g)u^k_{A,i}\neq u^k_{A,j}$ for $i\neq j$ and for any $g\in G$, then we can explicitly describe the isomorphism $S$ above. Indeed, in this case we can guarantee that the vectors

We describe here the isomorphism $S$ mentioned above.  Let $f$ be the morphism in  $\Hom_\mathbb{C}((\FF{k^*}{W_A}^*,\FF{k}{W_B})$ corresponding to $u^{k^*}_{A,i}\o u^k_{B,j}$ (this is, $f(\omega)=\omega(u^{k^*}_{A,i})u^k_{B,j}$). To present $f $ as an element $S(f)\in \Hom_G(W_A^*[\chi_k],W_B[\chi_k])$ we proceed as follows:
\begin{enumerate}
% \item We choose a subset $H\subset G$ such that, for any $i=1\div m_k(a)$, $\{h(u^{k}_{A,i})^*\}_{h\in H}$ is a basis of a $G$-representation (isomorphic to $\N{k}$). Here $(u^{k}_{A,i})^*\in \FF{k}{W_A^*}$ is simply a dual vector to $(u^{k}_{A,i})$ with respect to the chosen basis.
\item Denote by $\{(u^{k^*}_{A,i})^*\}_i\subset  \FF{k}{W_A^*}$ the dual basis for $\{u^{k^*}_{A,i}\}_i$. Choose a subset $H\subset G$ such that, for any $i=1\div m_k(a)$, $\{h(u^{k^*}_{A,i})^*\}_{h\in H}$ is a basis of a subrepresentation in $W_A^*[\chi_k]$ (necessarily isomorphic to $\N{k}$).
% \item By changing $i$ above we obtain a basis of $W_A^*[\chi_k]$.
\item Then, $\{h(u^{k^*}_{A,i})^*\}_{h\in H,i=1\div m_k(a)}$ is a basis of $W_A^*[\chi_k]$.
%\red{If we add this we should remove the last sentence of the following point?}
\item We define $S(f)(h(u^{k^*}_{A,i})^*)=hu^k_{B,j}$ and $S(f)(h(u^{k^*}_{A,l})^*)=0$ for $l\neq i$. This is the natural $G$-equivariant morphism associated to $f$.
%Notice that $S(f)$ is well defined as the vectors $h(u^{k^*}_{A,l})$, $h\in H$, $1\leq l\leq m_k(a)$ are a basis of $W_A^*[\chi_k]$. %\red{I do not know how detailed we should be here in the explenation?}

% \item We identify $S(f)$ with its image through the canonical isomorphism $\Hom_G(W_A^*[\chi_k],W_B[\chi_k]) \cong (W_A[\chi_{k^*}]\otimes W_B[\chi_k])^G$.
\end{enumerate}
\end{rem}
The following statement claims that under stronger assumptions, the image of $S(u^{k^*}_{A,i}\o u^k_{B,j})$ under the canonical isomorphism $\Hom_G(W_A^*[\chi_k],W_B[\chi_k]) \cong (W_A[\chi_{k^*}]\otimes W_B[\chi_k])^G$ has the particularly nice form of the averaging operator.
%This is, in this case we have
Namely,
\begin{lem}
If we can choose $H$ to be a subgroup of $G$, then
\begin{eqnarray*}
S(u^{k^*}_{A,i}\o u^k_{B,j})=\frac{n_k}{|G|} \sum_{g \in G} (g u^{k^*}_{A,i}) \otimes (g u^k_{B,j}),
\end{eqnarray*}
where $n_k$ is the cardinality of $H$ (equal to the dimension of $N_k$).
\end{lem}

\begin{proof}
First we notice that the kernel of the averaging operator always contains the kernel of $S(u^{k^*}_{A,i}\o u^k_{B,j})$. Moreover, the result of the averaging operator is a $G$-equivariant homomorphism. It remains to evaluate it on $(u^{k^*}_{A,i})^*$.
Notice that for each $h\in H$, we have $G=\{g^{-1}h: g\in G\}$. Hence we have the following equality
\begin{eqnarray*}
\sum_{g \in G} \left  (u^k_{A,i} \right)^*\hspace{-1mm} \left (g u^{k^*}_{A,i} \right ) \; g u^k_{B,j}=
\sum_{g \in G} \left  (u^k_{A,i} \right)^*\hspace{-1mm} \left (g^{-1}h u^{k^*}_{A,i} \right ) \; g^{-1}h u^k_{B,j}.
\end{eqnarray*}
Therefore, the right hand part of the equality in the lemma when evaluated at $(u^{k^*}_{A,i})^*$ is:
\begin{eqnarray}\label{eq:one}
\frac{n_k}{|G|} \sum_{g \in G} \left ( u^{k^*}_{A,i} \right )^*\hspace{-1mm} \left ( g u^{k^*}_{A,i} \right )  \; g u^k_{B,j} & = & %\frac{1}{|G|} \sum_{g \in G,h\in H} <g^{-1}h u^{k*}_{A,i},(u^{k*}_{A,i})^*> g^{-1}h u^k_{B,j} =
\frac{1}{|G|} \sum_{g \in G,h\in H}g^{-1}  \left ( \; g \left (u^{k^*}_{A,i} \right )^* \hspace{-1mm} \left (h u^{k^*}_{A,i} \right ) \; h u^k_{B,j} \; \right )= \nonumber \\ & = & \frac{1}{|G|} \sum_{g \in G}g^{-1} \left ( \; g \left (u^{k^*}_{A,i} \right )^* \hspace{-1mm} \left (h u^{k^*}_{A,i} \right ) \; h u^k_{B,j} \; \right ).
\end{eqnarray}
On the other hand we have
\begin{equation}\label{eq:two}
S \left (u^{k^*}_{A,i}\o u^k_{B,j} \right ) \left ((u^{k^*}_{A,i})^* \right )=\frac{1}{|G|} \sum_{g \in G}g^{-1} \left (\; S \left (u^{k^*}_{A,i}\o u^k_{B,j} \right ) (g(u^{k^*}_{A,i})^*) \; \right ),
\end{equation}
as $S(u^{k^*}_{A,i}\o u^k_{B,j})$ is equivariant.
% This equals:
%$$\frac{1}{|G|} \sum_{g \in G,h\in H} (g^{-1}h u^k_{A,i})((u^k_{A,i})^*) (g^{-1}h u^k_{B,j})=\frac{1}{|G|} \sum_{g \in G,h\in H} (h u^k_{A,i})(g(u^k_{A,i})^*) (g^{-1}h u^k_{B,j})$$
So far we did not use the fact that $H$ is a subgroup. However, in such case $Hu^{k^*}_{A,i}$ is the dual basis to $H((u^{k^*}_{A,i})^*)$, i.e.~$h(u^{k^*}_{A,i})^*=(hu^{k^*}_{A,i})^*$. In particular,
\begin{eqnarray*}
g(u^{k^*}_{A,i})^*=\sum_{h\in H} g\left (u^{k^*}_{A,i} \right )^* \hspace{-1mm} \left (hu^{k^*}_{A,i}\right ) \; (hu^{k^*}_{A,i})^*=\sum_{h\in H}  g\left (u^{k^*}_{A,i} \right )^* \hspace{-1mm} \left (hu^{k^*}_{A,i}\right ) \; h(u^{k^*}_{A,i})^*.
\end{eqnarray*}
Substituting this in (\ref{eq:two}) clearly yields the expression in (\ref{eq:one}), by point 3 in Remark \ref{S_operator}.
%are linearly independent.}
%\red{ We can explicitly describe the isomorphism $S$ above. We have that $W_A[\chi_k] $ is spanned by vectors of type $\rho_a(g) u^k_{A,i}$, $g \in G$. As $\{u^k_{A,i}\}$ are linearly independent, by Steinitz theorem there is a basis of $W_A[\chi_k] $ formed by $\{\rho_a(g) u^k_{A,i} | g \in G^k_{A,i}\}_i$ for certain subsets $G^k_{A,i}\subset G$ (these subsets can be explicitly found in each example).
%Let $e_{i,j}$ be the morphism in $Hom_\mathbb{C}(((W_A)_{k^*})^*,(W_B)_k)$ corresponding to $u^k_{A,i}\o u^k_{B,j}$ (that is, $e_{i,j}({u^k_{A,s}}^*)=u^k_{B,j}$ if $s=i$ and is 0 otherwise).
%Then the natural isomorphism from $Hom_\mathbb{C}(((W_A)_{k^*})^*,(W_B)_k)$ to  $Hom_G(W_A[\chi_k]^*,W_B[\chi_k])$ sends $e_{i,j}$ to $\tilde{e_{i,j}}$, which is the map that sends $\rho^*(g) {u^k_{A,s}}^*$ to $\rho(g) e_{i,j}({u^k_{A,s}}^*)$. Thus, $\tilde{e_{i,j}}(\rho(g) {u^k_{A,s}}^*)=\rho(g) {u^k_{B,j}}^*$ if $s=i$ and is equal to zero otherwise. Now, using the natural isomorphism $Hom_G(W_A[\chi_k]^*,W_B[\chi_k])\cong (W_A[\chi_k] \otimes W_B[\chi_k])^G$, $\tilde{e_{i,j}}$ is mapped to $$\sum_{s, g \in G^k_{A,i}} \rho_a(g) u^k_{A,s} \otimes \tilde{e_{i,j}}(\rho^*(g) {u^k_{A,s}}^*)=\sum_{g \in G^k_{A,i}} \rho_a(g) u^k_{A,i} \otimes \rho_b(g) u^k_{B,j}.$$
%Therefore,
%$$S(u^k_{A,i}\o u^k_{B,j})=\sum_{g \in G^k_{A,i}} \rho_a(g) u^k_{A,i} \otimes \rho_a(g) u^k_{B,j}.$$
%}
\end{proof}

% \begin{rem}
% %Each bipartition $A|B$ induces a Schur basis.
% If all the irreducible representations of $G$ are 1-dimensional (this is the case if $G$ is abelian), the basis $\Ba_{A| B}$ does not depend, on the particular bipartition $A|B$ chosen., \tp{mm: I am not sure, I think I wrote previous sentence, but now I would remove this remark. I think that there is some dependance of choices made. I mean in abelian case there is an (up to scalar) canonical choice, but as we present the construction, it seems there is freedom.}
% This is the case for Kimura 3-parameter, strand symmetric and general Markov models for which we can consider the Fourier basis (see section 7). %can be taken to be of the form $\{ \bxx_1 \ldots \bxx_n \mid \xx_i\in \Sigma\}$.
% \end{rem}

\subsubsection{Some examples}

For the models of Example \ref{main_equivariantmodels}, we consider the Fourier basis of the space $W$, defined as $\underline{\Sigma}:=\{\baa,\bcc,\bgg,\btt\}$ of $W$ where
\begin{eqnarray*}
 \baa =  \aa+\cc+\gg+\tt; \quad  \bcc =  \aa+\cc-\gg-\tt; \quad \bgg = \aa-\cc+\gg-\tt; \quad   \btt  =   \aa-\cc-\gg+\tt.
\end{eqnarray*}
Notice that $\baa$ equals the vector $\mathbf{1}$ introduced above and is invariant under the action of any permutation of  $\mathfrak{S}_4$. %
Notice also that the permutation groups associated to these models have only real characters; so for every irreducible representation, it holds $k^*=k$.
Throughout this section, we adopt the following notation: given $\underline{\mathtt{X}}_i\in \underline{\Sigma}$, we write $\underline{\mathtt{X}}_1\ldots \underline{\mathtt{X}}_m$ for the tensor $\underline{\mathtt{X}}_1\o \ldots \o \underline{\mathtt{X}}_m\in \o^m W$.
%Moreover, the basis $\bar{\Sigma}=\{\baa,\bcc,\bgg,\btt\}$ is orthonormal.

\begin{example}\label{examples2}
To illustrate the first step of the algorithm of the previous section, we proceed to obtain basis of each space $\FF{k}{W}$ when the group $G$ is chosen according to some of the models of Example  \ref{main_equivariantmodels}. All these models satisfy the following property:
\begin{itemize}
\item[(*)] the isotypic components of $W$ can be spanned by some elements of the Fourier basis above.
\end{itemize}
%$W$ can be decomposed into $G$-irreducible modules, each one of them being spanned
Namely,
\begin{enumerate}
 \item $G=\{1\}$ (GMM). In this case, the only representation is the \emph{identity} representation. We can take  $u^1_1=\baa$, $u^1_2=\bcc$, $u^1_3=\bgg$, $u^1_4=\btt$, which form a basis of $\FF{1}{W}=W$.

\item $G=\langle (\aa \tt)(\cc \gg)\rangle \cong \mathbb{Z}_2$ (strand symmetric model). There are two irreducible representations: the \emph{identity} and the \emph{sign} representations.
%Then, $W=W[\chi_1]\oplus W[\chi_2]$.
By taking $u^1_1=\baa$, $u^1_2=\btt$,
$u^2_1=\bcc$, $u^2_2=\bgg$, we have that $\{u^1_1,u^1_2\}$ and $\{u^2_1,u^2_2\}$ are basis of $\FF{1}{W}=W[\chi_1]$ and $\FF{2}{W}=W[\chi_2]$, respectively.

\item $G=\langle  (\aa \cc)(\gg \tt), (\aa\gg)(\cc\tt)\rangle \cong \mathbb{Z}_2 \times \mathbb{Z}_2$ (Kimura 3-parameter). There are four irreducible representations, each with dimension one (since $G$ is abelian). Then, we can take
$u^1_1=\baa$,  $u^2_1=\bcc$,  $u^3_1=\bgg$ and  $u^4_1=\btt$,
so that each $\FF{k}{W}$ is spanned by the corresponding $u^k_1$.
%
% $W=W[\chi_1]\oplus W[\chi_2]\oplus W[\chi_3]\oplus W[\chi_4]$, where $W[\chi_1]=\langle \baa \rangle$,
% $W[\chi_2]=\langle \bcc \rangle$, $W[\chi_3]=\langle \bgg \rangle$ and $W[\chi_4]=\langle \btt \rangle$.

%\tp{Note that each vector of the Fourier basis spans some isotypic components of $W$. This corresponds to the identification of $\Sigma$ with the group $(\mathbb{Z}_2 \times \mathbb{Z}_2,+)$ defined by
%\begin{eqnarray}
% \label{id_group}
% \aa \mapsto (0,0); \quad  \cc \mapsto (1,0); \quad
% \gg\mapsto (0,1); \quad  \tt \mapsto (1,1).
%\end{eqnarray}
%The characters of the irreducible representations of $G$ will play an important role in Section 7, and we will  denote them as $\chi_{\aa}=\chi_1$, $\chi_{\cc}=\chi_2$, $\chi_{\gg}=\chi_3$ and $\chi_{\tt}=\chi_4$.}

\item $G=\langle  (\aa \cc \gg \tt), (\aa\gg)\rangle$  (Kimura 2-parameter). There are two irreducible representations for $G$ with dimension 1. Taking $u^1_1=\baa$, $u^2_1=\bgg$, we obtain  $\FF{k}{W}=\langle u^k_1 \rangle$, for $k=1,2$. There is still a 2-dimensional irreducible representation; we can take $u^3_1=\bcc$ to get a basis of the corresponding space $\FF{3}{W}$ (a different possibility would be to take $u^3_1=\btt$). %
%
% $=with Then, $W=W[\chi_{1}]\oplus W[\chi_{2}]\oplus W[\chi_3]$, where $W[\chi_{1}]=\langle \baa \rangle$,  $W[\chi_{2}]=\langle \bgg \rangle$, and $W[\chi_3]=\langle \bcc, \btt \rangle$ ($\chi_3$ corresponds to the 2-dimensional irreducible representation of $G$).

\item $G=\mathfrak{S}_4$  (Jukes-Cantor). There are five irreducible representation, but only two of them appear in the Maschke decomposition of $W$: the identity representation and one 3-dimensional representation with character $\chi_4$ in Table \ref{table:char_tables}. By taking $u^1_1=\baa$ and $u^4_1=\bcc$, we obtain bases for the spaces $\FF{1}{W}$ and $\FF{4}{W}$.
% only  with characters $\chi_i$, $i=1,\ldots,5$. Then, $W=W[\chi_1]\oplus W[\chi_{4}]$, where $W[\chi_1]=\langle \baa \rangle$, and  $W[\chi_4]=\langle \bcc, \bgg, \btt \rangle$.

\end{enumerate}
%The above assumption does not hold for any equivariant model. For example,
%\begin{rem}
 \end{example}

\begin{rem}
There exist equivariant models that do not satisfy the property (*) above. For example, if $G=\langle  (\aa\cc) \rangle  \cong \mathbb{Z}_2$, there are two irreducible representations $\chi_1,\chi_2$ and the Maschke decomposition of $W$ becomes
$  W=W[\chi_1]\oplus W[\chi_2]$,   where $W[\chi_1]=\langle \baa \rangle \oplus \langle \bcc \rangle \oplus \langle \bgg+\btt \rangle$ and $W[\chi_2]=\langle \bgg-\btt \rangle$.
%Therefore, in this case, the Fourier elements do not span $G$-irreducible modules.
\end{rem}
%\end{rem}

\begin{example}\label{example_SS}
\textbf{A basis linked to a bipartition for the strand symmetric model.}
Take $G=\langle (\aa \tt)(\cc \gg )\rangle \cong \mathbb{Z}_2$, so we deal with the strand symmetric model.
%
% There are two irreducible representations of the group $G$: $\N{1}, \N{2}$. Each representation $\N{i}$ has the character $\chi_i$ shown in the table \ref{table:char_tables}.
The character table of $G$ is shown in Table \ref{table:char_tables}.

\begin{table}
\begin{center}
\begin{tabular}{ l || c |  c |  c |  c | c  } $\Omega_{\mathfrak{S}_4}$ & $\mathrm{id}$ & $(\aa\cc)$ & $(\aa\cc\gg)$ & $(\aa\cc\gg\tt)$ & $(\aa\cc)(\gg\tt)$ \\ \hline
$\chi_1$ & 1 & 1 & 1 & 1 & 1 \\
$\chi_2$ & 1 & -1 & 1 & -1 & 1 \\
$\chi_3$ & 2 & 0 & -1 & 0 & 2 \\
$\chi_4$ & 3 & 1 & 0 & -1 & -1 \\
$\chi_5$ & 3 & -1 & 0 & 1 &- 1 \\
\hline $\chi$     & 4 & 2 & 1 & 0 & 0
  \end{tabular} \quad
  \begin{tabular}{ l || c |  c } $\Omega_{\mathbb{Z}_2}$ & $\mathrm{id}$ & $(\aa\tt)(\cc \gg)$ \\ \hline
$\chi_1$ & 1 & 1  \\
$\chi_2$ & 1 & -1 \\
\hline $\chi$  & 4 & 0
  \end{tabular}
  \end{center}
\caption{\label{table:char_tables} Character tables of the groups $\mathfrak{S}_4$ and $G=\langle (\aa\tt)(\cc\gg)\rangle$. The character $\chi$ corresponds to the permutation representation of the group on the space $W$.}
\end{table}

\begin{center}
\end{center}
The permutation representation of $G$ decomposes as $\chi=2\chi_1+2\chi_2$, and $W=W[\chi_1] \oplus W[\chi_2]$,
with $W[\chi_1]=\langle \baa, \btt \rangle$ and $W[\chi_2]=\langle \bcc,\bgg\rangle$.

On the tree $12|34$, we consider the edge split $A=\{1,2\}$, $B=\{3,4\}$, $\alpha=2$, $\beta=3$.
%As in the Jukes-Cantor case, we construct first a convenient basis for the space $(\o^3 W)^G$.
%
The vectors $u^1_1=\baa$, $u^1_2=\btt$, $u^2_1=\bcc$, $u^2_2=\bgg$ regarded as vectors in $W_{\alpha}$
%$(W_{\alpha})_1=\langle \baa, \btt \rangle$ and $(W_{\alpha})_2=\langle\bcc,
induce tensors in $\FF{1}{W_A}$ and $\FF{2}{W_A}$, just by tensoring with $\mathbf{1}=\baa$ on the left: $u^1_{A,1}=\baa \baa$,  $u^1_{A,2}=\baa \btt$, $u^2_{A,1}=\baa \bcc$ and $u^2_{A,2}=\baa \bgg$.  We extend these tensors to a basis of $\FF{1}{W_A}$ and $\FF{2}{W_A}$ with
\begin{eqnarray*}
 u^1_{A,3}=\btt\baa, \; u^1_{A,4}=\btt\btt, \;u^1_{A,5}=\bcc\bcc, \; u^1_{A,6}=\bcc\bgg, \; u^1_{A,7}=\bgg\bcc, \;u^1_{A,8}=\bgg\bgg, \\
u^2_{A,3}=\btt\bcc, \;u^2_{A,4}= \btt\bgg, \;u^2_{A,5}=\bcc\baa, \; u^2_{A,6}=\bcc\btt, \;u^2_{A,7}=\bgg\baa, \;u^2_{A,8}=\bgg\btt.
\end{eqnarray*}
% %
% This gives two distinguished copies of $\N{1}$ and $\N{2}$ in $W_A$.
% %
% % mapping $\baa=b^1_1\o u^1_{\alpha,1}$ to $\baa\baa=b^1_1\o u^1_{A,1}$, where $u^1_{A,1}$ is the composition of $u^1_{\alpha,1}$ with the map $\cdot \mapsto \baa \o \cdot$. %
% There are still six other copies of $\N{1}$ and $\N{2}$ in $W_A$. Following the algorithm of Section 5.1, we extend Actually,
% \begin{eqnarray*}
%  (W_A)_1=\langle \baa\baa,\baa\btt, \btt\baa, \btt\btt,\bcc\bcc, \bcc\bgg, \bgg\bcc,\bgg\bgg \rangle; \\
%  (W_A)_2=\langle \baa\bcc,\baa\bgg, \btt\bcc, \btt\bgg, \bcc\baa,\bcc\btt, \bgg\baa, \bgg\btt \rangle.
% \end{eqnarray*}
% %
% and we take $u^1_{A,3}=\btt\baa$, $u^1_{A,4}=\btt\btt$, $u^1_{A,5}=\bcc\bcc$, $u^1_{A,6}=\bcc\bgg$, $u^1_{A,7}=\bgg\bcc$,$u^1_{A,8}=\bgg\bgg$,
% $u^2_{A,3}=\btt\bcc$, $u^2_{A,4}= \btt\bgg$, $u^2_{A,5}=\bcc\baa$, $u^2_{A,6}=\bcc\btt$, $u^2_{A,7}=\bgg\baa$, $u^2_{A,8}=\bgg\btt.$
We proceed similarly for $B$, and then we construct the basis $\{S(u^k_{A,i}\otimes u^k_{B,j})\}_{k,i,j}$ of $(\o^4 W)^G$. As the two irreducible representations of $G$ are $1$-dimensional, the $S$ operator has no effect and $\{u^k_{A,i}\otimes u^k_{B,j}\}_{k,i,j}$ is already a basis linked to $A|B$.
\end{example}

\subsection{Explicit Edge invariants}\label{sub:EdgeInv}

Once an edge split $A|B$ of the tree topology $T$ is given, edge invariants associated to it arise as restrictions on the rank of some matrices $M_k$, $k=1,\ldots, t$.
Our goal here is to explain how these matrices arise, and investigate how these rank restrictions look like.

% \tp{Somewhere, clarify (the second isomorphism by Schur)
% \begin{eqnarray*}
% \Hom_G(W_A^*[\chi_k],W_B[\chi_k]) \cong  \Hom_G(\N{k}, N_k) \o \Hom_{\CC}((W_A^*)_k,(W_B)_k) \cong \Hom_{\CC}((W_A)_k,(W_B)_k).
% \end{eqnarray*}
% \red{Moved to Proposition}
% }
The decomposition (\ref{dec}) allows us to understand any tensor $p\in (\o^n W)^G$ as a collection $(g^1_p,g^2_p,\ldots g^t_p)$, where each $g^k_p: \FF{k}{W_A^*}\rightarrow \FF{k}{W_B}$ is a linear map.

\begin{defn}[Thin flattening]
The collection of linear maps constructed above is referred to as the \emph{thin flattening of $p$ relative to the bipartition $A|B$}:
$ \Tf_{A|B}(p) =(g_p^1,g_p^2,\ldots g_p^t)$.
\end{defn}

The main result of \cite{casfer2010} claims that if $p$ is a (general) point in $CV_G(T)$, then the bipartition $A|B$ is an edge split in $T$ if and only if
\begin{eqnarray*}
 \rk \, g_p^k \leq m_k, \qquad \mbox{for every } k=1,\ldots, t.
\end{eqnarray*}
The $(m_k+1) \times (m_k+1)$ minors of matrices representing $g_k$ are usually known as \textit{edge invariants}.

%
%Each invariant tensor $\LL_n^G$ may be identified with an element of $$\Hom_G(W_A^*,W_B)=\sum_k\Hom_G(W_A^*[\chi_k],W_B[\chi_k])=\sum_k\Hom_{\CC}((W_A^*)_{k},(W_B)_k).$$
%\begin{rem}\label{rem:distsub}
We consider the basis $\mathcal{B}_{A|B}$ of $(\op{n})^G$ linked to the edge split $A|B$  constructed in section \ref{sec:basis}.
%$(W_A)^G$ and $(W_B)^G$.
%{\red{Do you mean of
%$(W_A)_{k^*}$ and $(W_B)_k$. %?? Or I misundersootd something ?}}
%
%Moreover, $(W_A)_{k^*}$ (respectively $(W_B)_k$)  contains a distinguished $m_k-$dimensional subspace spanned by the vectors of the basis $\baa^{a-1}\o u^k_{i}$, $i=1\div m_k$, (resp. $u^k_{j}\o \baa^{b-1}$, $j=1\div m_k$).
%\end{rem}
%
%By considering the dual basis, we have a natural basis of $(W_A^*)_{k}$.
As we have fixed bases $\{u^{k^*}_{A,i}\}_i$ of $\FF{k^*}{W_A}$ and $\{u^k_{B,j}\}_j$ of $\FF{k}{W_B}$, each tensor in $(\o^n W)^G$ naturally induces  matrices $M_k$ representing the morphisms $g^k_p\in\Hom_{\CC}(\FF{k^*}{W_A},\FF{k}{W_B})$ of the thin flattening.
Each rank restriction for $M_k$, is an equation on the coordinates $q^k_{i,j}$ introduced in section \ref{sec:basis}.
%$$\left\{q^k_{ij}\mid 1\leq i \leq m_{k^*}(a), 1 \leq j \leq m_k(b) \right\}.$$

In order to obtain a complete intersection, we shall now choose specific  minors of order $m_k+1$ in the matrices $M_k$.
%By Remark \ref{rem:distsub}\red{This remark was removed (I do not know when, but I think not by me). Should we add it back? Or just refer to the first basis vectors?},
The basis we constructed for $\FF{k}{W_B}$ (respectively $\FF{k^*}{W_A}$) has a distinguished set of $m_k$ (respectively $m_{k^*}$) elements, namely the first $m_k$ (resp. $m_{k^*}$) elements. We call $M_k ^0$  the submatrix of $M_k$ corresponding to these elements. %$G$-modules isomorphic to $N_k$ (resp. $N_{k^*}$) .

%In particular,
We choose only the $(m_{k^*}+1)\times(m_k+1)$ minors of $M_k$ that contain the distinguished $m_{k^*}\times m_k$-submatrix $M_k^0$.
As in our setting we have $k=k^*$, we observe that $M_k^0$ is a square matrix.

For the purpose of the next section, we need to write these minors in terms  of the determinant of $M_k^0$, $\Delta_k(p)=\operatorname{det} M_k^0$.
%To simplify the notation let us say
Note that $M_k^0$ is the upper left $m_k \times m_k$ submatrix of $M_k$ so that $\operatorname{det} M_k^0$ is a polynomial in indeterminates $q^k_{ij}$ for $1\leq i\leq m_{k^*}$, $1\leq j\leq m_k$.
%\begin{eqnarray}\label{Mk0}
%M_k^0= \left (
% \begin{array}{ccc}
%  & \vdots & \\
%  \ldots & q^k_{xy} & \ldots \\
%  & \vdots &
% \end{array}
% \right )_{1\leq x,y\leq m_k }
%\end{eqnarray}
%with rows and columns indexed by the coordinates $\{q^k_{i,j}\}_{1\leq i,j \leq m_k}$.
%
%
%\private{\red{Wouldn't it be better to reorder the basis above? So that the coordinates already appear in this order. DONE!!!}}

We note by $E^k_{i,j}$ the minor containing $M^0_k$, the row indexed by $u^{k^*}_{A,i}$ and the column indexed by $u^k_{B,j}$.
Then the minors $E^k_{i,j}$
%\begin{eqnarray*}
%E^k_{ij}=0, \qquad k=1 \div  t, \; i=m_{k^*}+1 \div m_{k^*}(a),  \;
%j=m_k+1 \div m_k(b)
%\end{eqnarray*}
 containing $M_k^0$, $k=1 \div  t$, $i=m_{k^*}+1 \div m_{k^*}(a)$,  $j=m_k+1 \div m_k(b)$, can be written as
\begin{eqnarray}\label{edge_invariants}
E^k_{ij}= q^k_{ij} \, \Delta_k(p)+\sum_{s=1}^{m_k} (-1)^{j+s}\,  q^k_{sj}\, \Delta^{(s)}_k(p)
\end{eqnarray}
where $\Delta_k^{(s)}(p)$ is the determinant of the matrix obtained by removing the $s$-th row of ${M}_k^0$ and adding the first $m_k$ entries of the $i$-th row of $M_k$.
The set of equations $E^k_{i,j}=0$ (which are a particular subset of the edge invariants for $A|B$) will be denoted by $\eq_{A|B}$. There are
$$N_{A|B}:= \sum_k (m_{k^*}(a)-m_{k^*})(m_k(b)-m_k)$$ such minors.

\begin{rem}\label{rem_Nab}
Notice that the cardinality $N_{A|B}$ of this set depends only on $a=|A|$ and $b=|B|$, and not of the particular choice of leaves in $A$ or $B$. Moreover, by Proposition \ref{prop:W^nG}, we have $N_{A|B}=m_1(n)-m_1(a+1)-m_1(b+1)+m_1(2)$. For the models of Example \ref{main_equivariantmodels}, explicit formulas for $m_1(s)$ are given in \cite[Prop. 5]{casferked}. From these formulas, it is easy to see that $N_{A|B}$ grows exponentially with $n$. Therefore, the list of local phylogenetic invariants we give in the following section has exponential cardinality in $n$, which would make it useless for practical applications when $n$ is big. %
However, it is well known that rank conditions do not need to be checked directly by evaluating minors; they can be checked by using the singular values of the matrix instead (this is the approach followed in \cite{eriksson2005,casfer2015}). Thus, for the reader interested in applying the local phylogenetic invariants provided in this paper, we suggest using singular values instead of $\eq_{A|B}$.
\end{rem}

The following result, which is needed in the next section, easily follows from \eqref{edge_invariants}.
\begin{lem}\label{derivatives}
%  For $k=1,\ldots,t$, we have
%  \begin{eqnarray*}
%   \frac{\partial \, E^k_{ij}}{\partial \, q^k_{ij}}=\Delta_k(p), \qquad \mbox{and} \qquad
%  \frac{\partial \, E^k_{ij}}{\partial \, q^k_{i'j'}}=0,
%  \end{eqnarray*}
%  where $m_{k^*}+1 \leq i,i'\leq m_{k^*}(a)$, $m_k+1\leq j,j'\leq m_k(b)$ and $(i,j)\neq (i',j')$.
 %
   For any $k$ and $m_{k^*}<i',j',i,j\leq m_{k^*}(a)$, we have
  \begin{eqnarray*}
  \frac{\partial \, E^k_{ij}}{\partial \, q^k_{i'j'}}= \left \{
   \begin{array}{cl}
     \Delta_k(p), & \mbox{ if }(i,j)=(i',j');  \\
     0, & \mbox{ otherwise.}
       \end{array}
   \right .
  \end{eqnarray*}

 \end{lem}

\subsection{Equations from $CV_G(T_A)$ and $CV_G(T_B)$}
%
%The reader may note that any $G$-invariant tensor in $W$ would equally work.
The following result is essentially well-known (see Lemma 1 of \cite{FuLi}, or \cite{allman2004b}). However, we prove it in our setting. 

\begin{figure}
\begin{center}
\includegraphics[scale=0.35]{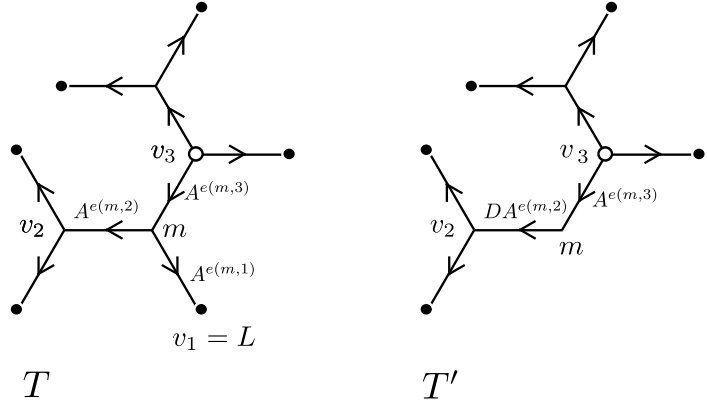}
\end{center}
\caption{\label{figure:remove_leaf}. Illustration of the proof of Lemma \ref{lem:cutleaf}. The tree $T'$ is obtained by taking the leaf $v_1=L$ off the tree $T$.}
\end{figure}

\begin{lem}\label{lem:cutleaf}
Let $T$ be a tree and let  $L$ be one of its leaves. Let $T'$ be the subtree of $T$ that has the same vertices, apart from $L$ (Fig. \ref{figure:remove_leaf}).
%We may regard $\mathbf{1}^*$ as a contraction map:
The following contraction map
\[\begin{array}{rcl}
f_L:\o^n W=\bigotimes_{l\in L(T)}W_l &\rightarrow &\bigotimes_{l\in L(T')}W_l\\
\otimes_{l \in L(T)} v_l &\mapsto & \mathbf{1}\cdot v_L  \left(\otimes_{l \in L(T')} v_l\right)
\end{array}\]
satisfies $f_L(\Im \Psi_T^G)=\Im \Psi_{T'}^G$ and, as a consequence, %where $\mathcal{P}_{|T'}$ is the restriction of the parameters on $T$ to $T'$.
$\overline{f_L(CV_G(T))}= CV_G(T')$.
\end{lem}

In stochastic terms, this map is called the \emph{marginalization over the random variable at $L$}.

\vspace{2mm}

%\private{I write the proof not only for trivalent. And I assume that $T'$ has one vertex of degree 2 because it says that $T'$ has the same vertices as $T$ a part from $L$}
\begin{proof}
The map $f_L$ is induced by the multilinear map
\[\begin{array}{rcl}
\prod_{l\in L(T)}W_l &\rightarrow &\bigotimes_{l\in L(T')}W_l\\
\left( v_l \right)_{l \in L(T)} &\mapsto & \mathbf{1} \cdot v_L  \left(\otimes_{l \in L(T')} v_l\right)
\end{array}\]
and therefore $f_L$ is well defined (by the universal property of tensor products).
%The contraction map $f_{\mathbf{1}}$ is well defined because it obviously comes from ... (this has to be written still)

Without loss of generality, we can assume that the interior node $m$ adjacent to $L$ in $T$ has degree $\geq 3$ (indeed, if it had degree two, then $CV_G(T)$ would be isomorphic to the variety associated to the tree with this vertex removed and two adjacent edges joined into a single edge).

%If the vertex $m\in L(T)$ adjacent to the leaf $L$ is trivalent, then the same vertex $m$ in $T'$ has degree $2$. In particular $CV_G(T')$ is isomorphic to the variety associated to the tree with this vertex removed and two adjacent edges joint into one edge. Therefore, without loss of generality, we can assume that the interior node $m$ adjacent to $L$ in $T$ has degree $\geq 3$.
%\tp{mm: I do not understand something, shoudn't the beginning of the paragraph be about the case, when the vertex has degree 2 and becomes a leaf?} \newline
%\tp{jf: I agree with Mateusz. Actually, I think we could skip the assumption ''degree of $m\geq 3$`` by defining $B^{e(m,3)}=A^{e(m,3)}D$ (that is, getting  $A^1 \mathbf{1}$ into $A^3$ instead of $A^2$).}

We call $v_1,\dots,v_t$ ($t\geq 3$) the vertices adjacent to $m$ and we set $v_1=L$. We root the tree $T$ at $v_t$ and call $e(m,i)$ the edges from $m$ to $v_i$, $i=1\dots,t-1$ (see figure \ref{figure:remove_leaf}).

Let $\mathcal{P}=\left(\pi,\left(A^e\right)_{e \in E(T)}\right)$ be a point in $\Par_G(T)$ (rooted at $v_t$). For the edge $e(m,2)$ from $m$ to $v_2$, we consider a new matrix $B^{e(m,2)}:= D A^{e(m,2)}$ where $D$ is the diagonal matrix $diag(A^{e(m,L)}\mathbf{1})$ formed by the entries of $A^{e(m,L)}\mathbf{1}$. Since $A^{e(m,2)}$ is $G$-equivariant, the vector $A^{e(m,2)} \mathbf{1}$ is $G$-invariant, and $D$ is $G$-equivariant again. It follows from this that the new matrix $B^{e(m,2)}$ is $G$-equivariant. For all other edges of $T'$, take  $B^e=A^e$. It is not difficult to check that
\begin{eqnarray*}
f_L \left( \Psi_{T}^G(\mathcal{P})\right)=\Psi_{T'}^G\left(\pi,\left(B^e\right)_{e \in E(T')}\right).
\end{eqnarray*}
Therefore $f_L(\Im \Psi_T^G)\subset \Im \Psi_{T'}^G$.
The other inclusion $f_L(\Im \Psi_T^G) \supseteq \Im \Psi_{T'}^G$ follows easily by adding the identity matrix at the edge $e(m,1).$

The equality of the parametrized part of the varieties implies equality in the closures, hence $\overline{f_L(CV_G(T))}= CV_G(T')$.
\end{proof}

By successive applications of Lemma \ref{lem:cutleaf} to all leaves in $A \setminus \{\alpha\} \subset L(T)$ (that is, marginalizing over all leaves in $A\setminus \{\alpha\}$) we obtain a map
$$f_{A\setminus \alpha}:\bigotimes_{l\in L(T)}W_l \lra \bigotimes_{l\in L(T_B)}W_l$$
that sends the variety $CV_G(T)$ to $CV_G(T_B)$.

In order to induce equations from $T_A$, $T_B$ to $T$, below we translate this map in terms of the corresponding affine coordinate rings. We do not explicitly write indeterminates nor coordinates because the above map $f_{A\setminus \alpha}$ is basis independent. This fact will play an important role in the proof of the main result in the next section.

%In order to induce equations from $T_i$ to $T$,
%we fix leaves $\alpha\in L(T_A)$, $\beta\in L(T_B)$, where $\alpha\neq L_A$ and $\beta\neq L_B$.

%After identifying $\op{n}$  and $\bigotimes_{l\in L(T_B)}W_l$ with its duals, we see that
The map $f_{A\setminus \{\alpha\}}$  above is dual to the map
$$%\left(\o \mathbf{1} \right)_{A \setminus \alpha}:
\bigotimes_{l\in L(T_B)}W_l^*  \rightarrow \bigotimes_{l\in L(T)} W_l^* $$
that maps $t$ to $\mathbf{1}^{a-1} \o t$ if the leaf  $L_B$ of $T_B$ is identified with the leaf $\alpha$ of $T$, so $f_{A\setminus \{\alpha\}}^*$ %$\left(\o \mathbf{1} \right)_{A \setminus \alpha}$
is the map corresponding to $f_{A\setminus \{\alpha\}}$ in terms of coordinates. Moreover, both maps restrict to $G$-invariant vectors. Summing up we have:

\begin{cor}\label{equations_subtree}
Any equation vanishing on $CV_G(T_B)$ extends to an equation vanishing on $CV_G(T)$ via the map:
% \begin{eqnarray*}
% f^*_{A \setminus \alpha}: \left(\bigotimes_{l\in L(T_B)}W_l^* \right)^G & \lra & \left (\bigotimes_{l\in L(T)} W_l^* \right )^G\\
% t&\mapsto& \mathbf{1}^{a-1} \o t
% \end{eqnarray*}
$$\begin{array}{ccc}
f^*_{A \setminus \alpha}: \left(\bigotimes_{l\in L(T_B)}W_l^* \right)^G & \lra & \left (\bigotimes_{l\in L(T)} W_l^* \right )^G\\
t&\mapsto& \mathbf{1}^{a-1} \o t
\end{array}
$$
where the leaf $L_B$ of $T_B$ is identified with the leaf $\alpha$ of $T$.
%and a tensor $t$ is mapped to %$t\otimes (\mathbf{1})^{|L(T_A)|-1}$.
%$$.
\end{cor}

Similarly, for any subsets $R \subsetneq S\subset L(T)$ of leaves of $T$ we have $f_S^*$ is the map
$$\begin{array}{rcl}
   f_S^*: \bigotimes_{l\in R\setminus S}W_l^*  & \lra & \bigotimes_{l\in R} W_l^* \\
   \otimes_{l\in R\setminus S} v_l &\mapsto & \otimes_{l\in R} w_l \, ,
  \end{array}
$$
where $w_l=v_l$ if $l\in R\setminus S$ and $w_l=\mathbf{1}$ if $l\in S$.
%\private{ We should be more precise maybe in this statement or in the proof???}

\begin{rem}
 It is convenient to write the equations of $T_B$ in the basis related to a bipartition $L_B |B$. In this way, the extension of a coordinate as defined in Corollary \ref{equations_subtree} gives rise to a coordinate that is already in the basis
 $\mathcal{B}_{A|B}$ of $(\o W^n)^G$ (indeed, as $\mathbf{1}$ is $G$-invariant, the operator \eqref{S_operator} does not affect it and it is easy to check that the extended basis are elements of $\mathcal{B}_{A|B}$).

%\red{mm: Here I just removed the word averaging. Indeed when we are dealing with invariants, then $\FF{1}{W_A}=W_A^G$ etc for other spaces. So $S$ does not do anything on invariants - it just maps a tensor to itself (in other words every vector in the whole representation is already in $\FF{1}$ so there is no need of extenstion).}

% WARNING: NOW THE AVERAGING OPERATOR CAN ONLY BE APPLIED IF WE HAVE THAT SUBGROUP H SUCH THAT ....(REMARK 5.1). SO THIS PART OF THIS REMARK CAN ONLY BE STATED IN THAT CASE. DO YOU AGREE?

% mm: so indeed we cannot refere to averaging operator, but I guess that the remark that a coordinate (basis vector) gives a good basis vector is correct - I will think about it and try to rephrase it.

% This change of basis will not be necessary in the case of abelian groups (for example, for GMM, strand symmetric model and Kimura 3-parameter)}
\end{rem}

\section{The main result}

Given a phylogenetic tree under an equivariant model $\MM_G$, the goal of this section is to construct a complete intersection for $CV_G(T)$ on a neighborhood of a point of no evolution. This will be done by using induction on the number of leaves of the tree.

Let $T$ be a tree with at least one interior edge and leaves $L(T)=\{l_1,\ldots,l_n\}$. Reordering the set of leaves (if needed) we can assume that there exists a node with children $\{l_{n-l},\dots, l_{n}\}$,  $l \geq 1$.
We take the edge split $A|B$ given by $A=\{l_1,\ldots,l_{n-l-1}\}$, $B=\{l_{n-l},\dots,l_n\}$. Write $e$ for the interior edge of $T$ associated with this split.
%and write $e_{A|B}$ for the interior edge corresponding to it.
%
% Let $z$ denote the middle point in $e_{A|B}$ and consider the decomposition
% (see Figure \ref{figure:decomposition})
% $T=\tilde{T}_A \ast_z \tilde{T}_B$.
% Equivalently, $\tilde{T}_A$ is the tree with leaves $A \cup \{z\}$, and $\tilde{T}_{B}$ is the tree with leaves $z\cup B$.
Keeping the notation already used throughout the paper, $T_A$ has leaves $A \cup \{L_A\}$, $T_B$ has leaves $B \cup \{L_B\}$ (where $L_A$, $L_B$ are defined as in figure 1). The variety associated to $T_A$ is the closure of the image of the polynomial map
\begin{eqnarray*}
\Psi_{T_A}: \mathrm{Par}_G(T_A)\rightarrow \op{n-l}^G.
\end{eqnarray*}
and the variety $CV_G(T_B)\subset \op{l+2}^G$ is the closure of the image of
\begin{eqnarray*}
 \Psi_{T_B}: \mathrm{Par}_G(T_B)\rightarrow \op{l+2}^G.
\end{eqnarray*}

%{\red{ We first need to prove that marginalizations of a point no evolution give points of no evolution of the corresponding subtrees:

\begin{lem}\label{marginal_noevol}
Given a leaf $L$ in $T$, the image of a point of no evolution $\pi_n \in \op{n}$ under the map $f_L$ is a point of no evolution in $\op{n-1}$. %I particularn $CV_G(T_A)$ (analogously in $CV_G(T_B)$).
\end{lem}
\begin{proof}
% , the image of a point of no evolution $\pi_n=\sum_{\xx \in \Sigma} \pi_{\xx} \xx  \ldots  \xx$ under the map $f_{\mathbf{1}}$ is a point in $CV_G(T')$. %(analogously in $CV_G(T_B)$).
%  We prove now that $f_{\mathbf{1}}(\pi_n)$ is a point of no evolution. Indeed,
If $\pi_n=\sum_{\xx \in \Sigma} \pi_{\xx} \xx  \o \ldots  \o \xx$, then its image under the map $f_L$ is
 $$f_L(\pi_n)=\sum_{\xx \in \Sigma}\pi_{\xx} (\mathbf{1}\cdot \xx)^a  \left(\xx \o\stackrel{n-a}{\ldots}\o \xx \right)=\sum_{\xx \in \Sigma} \pi_{\xx} \xx\o \stackrel{n-a}{\ldots} \o \xx .$$
As $\pi_n$ was invariant by the action of $G$, so is $f_L(\pi_n)$ and therefore it is a point of no evolution in $\op{n-a}$.
\end{proof}

% By successively applying lemma \ref{lem:cutleaf}, the image of a point of no evolution $\pi_n$ under the map $f_{\mathbf{1}}$ is a point in $CV_G(T_A)$ (analogously in $CV_G(T_B)$). Moreover its images are still points of no evolution. Indeed, if $\pi_n$ is a point of no evolution in $CV_G(T)$, then $\pi_n \in \op{n}^G$ and $\pi_n=\sum_{\xx \in \Sigma} \pi_{\xx} \xx  \ldots  \xx$. Therefore $$f_{\mathbf{1}}(\pi_n)=\sum_{\xx \in \Sigma}\pi_{\xx} \mathbf{1}\cdot \xx  \left(\xx \o\stackrel{n-1}{\ldots}\o \xx \right)=\sum_{\xx \in \Sigma} \pi_{\xx} \xx\o \stackrel{n-1}{\ldots} \o \xx $$
% and as $\pi_n$ was invariant by the action of $G$, so is $f_{\mathbf{1}}(\pi_n)$. Successively applying this we obtain the claim.}}

By successively applying the above lemma and lemma \ref{lem:cutleaf} we obtain points of no evolution in $CV_G(T_A)$ and in $CV_G(T_B)$ from a point of no evolution $\pi_n$ in $CV_G(T)$. This shall allow us to apply an induction argument.

%Let $T_{a+2}$ be the $(a+2)$-pod tree evolving under $\MM_G$ and $L(T_{a+2})=\{x_0,\dots,x_{a+1}\}$, and let $\eq_{T_{a+2}}:=\{h_1, h_2, \ldots, h_{\codim(a+2)}\}$ be a set of equations of a complete intersection
%that defines $CV_G(T_{a+2})\subset (\o^{a+2} W)^G$ on an open subset containing general points of no evolution. As we already proved, general points of no evolution are smooth by Theorem \ref{thm:orbit}. In particular, the variety is locally a complete intersection, which guaranties the existence of $\eq_{T_{a+2}}$.
Let $T_{d}$ be the claw tree with $d$ leaves (claw $d$-tree) evolving under $\MM_G$ and $L(T_{d})=\{x_0,\dots,x_{d-1}\}$, and let $\eq_{T_{d}}:=\{h_1, h_2, \ldots, h_{\codim(CV_G(T_{d}))}\}$ be a set of equations of a complete intersection
that defines $CV_G(T_{d})\subset (\o^{d} W)^G$ on an open subset containing general points of no evolution. As we already proved, general points of no evolution are smooth by Theorem \ref{thm:orbit}. In particular, the variety is locally a complete intersection, which guaranties the existence of $\eq_{T_{d}}$.
Before proceeding with induction, we need the following assumption about the equivariant model $\MM_G$ on the claw tree with $d$ leaves, which shall be checked for every particular equivariant model and every $d$ equal to a degree of one of the interior nodes of $T$. For the GMM, the strand symmetric model, and the Jukes-Cantor model,
we prove in section 6 that this assumption holds for the tripod (and hence our result is valid for trivalent trees evolving on these models). A local complete intersection for the Kimura 3-parameter model for trivalent trees was already given in \cite{casfer2008}.

% \begin{figure}
% \begin{center}
% \includegraphics[scale=0.3]{decomposition2.eps}
% \end{center}
% \caption{\label{figure:decomposition}.}
% \end{figure}

%\red{I renamed the Assupmtion to Tripod hipothesis --Marta}
\begin{assumption}\label{assumption}
We write equations $\eq_{T_{d}}$ on a basis of type $\mathcal{B}_{x_0\mid \{x_1,\dots,x_{d-1}\}}$ following subsection \ref{sub:EdgeInv}. The jacobian  of these new equations $\eq_B$ , which we denote as $\J{x_0\mid x_1,\cdot,x_{d-1}}(T_{d})$,
has rank equal to $\codim(CV_G(T_d))$ at any general point of no evolution. We denote by $\J{x_0\mid x_1,\cdot,x_{d-1}}^*(T_{d})$  the matrix obtained from $\J{x_0\mid x_1,\cdot,x_{d-1}}(T_{d})$ by removing the columns corresponding to $S(u^k_{x_0,i}\o u^k_{x_1,\cdot,x_{d-1},j})$ for $k=1\div t$ and $i,j=1\div m_k$.
We say that the equivariant model $\mathcal{M}_G$ satisfies \emph{the claw $d$-tree hypothesis} if
%\subset \CC[q^k_{i,j}]$, so that
\begin{eqnarray}\label{assumption_eqn}
 \rk \, \J{x_0\mid x_1,\cdot ,x_{l-d} }^*(T_{d})=\codim(CV_G(T_d)),
\end{eqnarray}
whenever this matrix is evaluated at a generic point of no evolution. %In particular, $\eq^{T_3}$ is a local minimal system of generators for the ideal of $CV_G(T_3)$ in $\LL_3^G$.
\end{assumption}

%We shall require that the claw $d$-tree hypothesis is satisfied for any $d$ in the set $D$ of degrees of interior nodes of $T$.
%
\paragraph{Induction hypothesis.}
We will use the following induction hypothesis:

%\begin{enumerate}
 %\item[(a)]
$(\ast)$ There is a set of equations $\eq_{T_A}=\{g_1,g_2,\ldots,g_{\codim(CV_G(T_A))}\}$ that defines  the variety $CV_G(T_A)\subset (\o^{n-l} W)^G$ scheme theoretically on an open subset containing general points of no evolution. %

\vspace*{3mm}
By Corollary \ref{equations_subtree}, the map $\tau \mapsto \tau \o \mathbf{1}^l$ induces new equations for $CV_G(T)$ from $\eq_{T_A}$. These equations shall be written in the coordinates $q^k_{i,j}$ corresponding to the basis $\mathcal{B}_{A|B}$ linked to the edge split $A|B$ and shall be called $\eq_{A}=\{f^{A}_1,\ldots,f^{A}_{\codim(CV_G(T_A))} \}\subset \mathbb{C}[q^k_{i,j}]$.

% \item[(b)] There is a set of equations $\eq_B=\{h_1,\ldots,h_{\codim(3)}\}\subset \CC[q^k_{i,j}]$ that defines the variety $CV_G(T_B)\subset \LL_{3}^G$ {\red{on a neighborhood containing generic points of no evolution and that satisfies the assumption \ref{assumption}}}.
% %
%
% \end{enumerate}

As above, by Corollary \ref{equations_subtree}, the map $\tau \mapsto \mathbf{1}^{n-l-2} \o \tau$
induces new equations
\[\eq_B= \{f^B_1,\ldots,f^B_{\codim(CV_G(T_B)}\}\subset \mathbb{C}[q^k_{i,j}]\]
  for $CV_G(T)$ from the set of equations  $\eq_{T_B}$ of the underlying model assumption.
Besides, we still need to consider the set of polynomials coming from the edge split.

\paragraph*{Edge invariants.}
As in subsection 4.2, for each $k=1,\ldots,t$, write
$M_k$ for the ${m_k(n-l-1)}\times {m_k(l+1)}$-matrix
with rows indexed by the $u^{k^*}_{A,i}$, $i=1,\ldots,m_k(n-l-1)$, columns indexed by $u^k_{B,j}$, $j=1,\ldots, m_k(l+1)$, and whose
 $(i,j)$-entry is the coordinate $q^k_{i,j}$.
For each of these matrices, take
the set of all the %$(m_k(n-2)-m_k)\times (m_k(2)-m_k)$
$(m_k+1)\times (m_k+1)$-minors containing the sub matrix $M_k^0$ defined in Section \ref{sub:EdgeInv}, with rows and columns indexed by $\{u^{k^*}_{A,i}\}_{i=1,\ldots,m_k}$ and $\{u^k_{B,j}\}_{j=1,\ldots,m_k}$, respectively. We obtain $N_{A|B}$
polynomials in $\CC[q^k_{i,j}]$ of the
form (\ref{edge_invariants}).
%Write $\eq_{A|B}$ for the whole set of these equations:
%$\eq_{A|B}=\bigcup_{k=1}^t \{E^k_{ij}\}_{i,j}$.

 \begin{lem}\label{number}
 We have
% \begin{enumerate}
%\item $N_{A|B}=  m_1(n)-m_1(n-a)-m_1(a+2)+m_1(2)$.
 % \item $N_{A|B}=  m_1(n)-m_1(n-l)-m_1(l+2)+m_1(2)$.
%  \item $\codim(n-a)+\codim(a+2)+N_{A|B}=\codim(n)$.
 $\codim(CV_G(T_A))+\codim(CV_G(T_B)+N_{A|B}=\codim(CV_G(T))$.
 %\end{enumerate}
\end{lem}

\begin{proof}
%1. It follows by direct computation  from Proposition \ref{prop:W^nG}.

%\noindent
We assume that $T$ has no vertices of degree $2$, as such nodes can be removed. By Theorem \ref{dimension} and Remark \ref{rem_Nab}, we have
\begin{eqnarray*}
 \codim(CV_G(T_A))+\codim(CV_G(T_B)+N_{A|B}= \\
 \left (m_1(l+2)-(|E(T)|-(n-l-1))(m_1(2)-m_1)-m_1\right ) +\\ \left ( m_1(n)-m_1(n-l)-m_1(l+2)+m_1(2)\right ) + \\ + \left ( m_1(n-l)-(n-l)(m_1(2)-m_1)-m_1\right ).
\end{eqnarray*}
The sum above equals
$$m_1(n)-|E(T)|(m_1(2)-m_1)-m_1=\codim(CV_G(T)).$$
%It follows easily from \eqref{codim}.
\end{proof}

%{\red {Define open subset $$\Omega=\{p \mid \Delta_k(p)\neq 0, k=1,\dots, t\}$$.
%Lemma \ref{nonzero_derivatives} proves that it contains no evolution points. And it contains also all other  biologically interesting points???}

\begin{thm}\label{invariantsCI} Let $T$ be a phylogenetic tree on $n$ leaves, $n\geq 3$; let $D$ be the set of degrees of its interior nodes, and assume that $d\geq 3$ for any $d \in D$.
Let $\mathcal{M}_G$ be an equivariant model that satisfies the claw $d$-tree hypothesis for any $d \in D$.
The %ideal generated by the set
set of equations
$\eq_{T}:= \eq_{A} \cup \eq_{B} \cup \eq_{A|B}$
defines the variety $CV_G(T)$ scheme theoretically on an open subset that contains general points of no evolution.

%{\red{But we cannot use induction hypothesis then!!!}}
%is a local complete intersection of the variety $CV_G(T)$.
\end{thm}

% \begin{rem}
% Notice that the polynomials in $K$ have degree 2 and those in $J(3)$ have degree 4. Since the polynomials in  $J(n-1)$ are computed from a subtree of $T$ with $n-1$ leaves, it follows that there are only quadrics and quartics in the above set of generators.
% \end{rem}

\begin{proof}
%First of all, direct computation shows
%that the number of polynomials being
%considered equals the  codimension of
%the variety $W$, i.e.
We proceed by induction on the number of leaves. The first step is $n=3$, which is covered by the claw $d$-tree assumption for $d=3$.
We assume thus $n>3$ and that $T$ has at least one interior edge which splits the leaves $L(T)=\{l_1,\ldots,l_n\}$ into two sets $A=\{l_1,\ldots,l_{n-l-1}\}$ and $B=\{l_{n-l},\dots,l_n\}$ (reordering leaves if necessary).
Consider the trees $T_A$ and $T_B$ as defined in the beginning of this section.
Note that  we are able to use the induction hypothesis stated above because the set of degrees for the interior nodes of the tree $T_A$ is included in $D$.
By Lemma \ref{number}, we know that
\begin{eqnarray*}
|\eq_A|+|\eq_B| + |\eq_{A|B}|  = \codim(CV_G(T)),
\end{eqnarray*}
that is, the number of equations equals the codimension of the variety.
% ........................................................................................................
We already know that
%the  ideal $\wp$ generated by
%$\eq^n_{A|B}$
$\eq_{T}$ are equations satisfied by all points in %is contained in $I(T)$, and
$CV_G(T)$. %is contained in the variety $V'$ defined by $\wp$.
%We claim that the tangent space to $V'$ of a generic point of no-evolution $\pi_n$  has codimension equal to the codimension of $CV_G(T)$.

Let $V'$ be the variety defined by $\eq_{T}$. %$\eq_{A|B}$.
Now, consider the jacobian matrix $\J{A|B}(T)$ obtained by taking the partial derivatives of the polynomials in $\eq_{A|B}$
with respect to the coordinates $q^k_{i,j}$ of $(\o^n W)^G$. We claim that the rank of this matrix at a generic  point of no evolution $\pi_n$ is maximal.
From this, we will deduce that $V'$ is non-singular in a neighborhood $U$ of $\pi_n$. Since $V'$ and $CV_G(T)$ have the same dimension, it follows that both varieties are equal in $U$.
%
%Now we provide the description of the tangent space to $V'$ at $\pi_n$.
%
% Write
% \begin{eqnarray*}
%  Q_0:= \left \{\frac{1}{\sqrt{d_k}} \sum_{l=1}^{d_k}  \baa^{n-3}\o \xx^k_{i,l} * \xx^k_{j,l} \o \baa : 1\leq i,j\leq m_k \right \}_{k=1,\ldots, t}
% \end{eqnarray*}

By reordering the
columns of the jacobian matrix if necessary, we may
assume that columns are indexed  as follows:
\begin{itemize}
\item[--] the first $m_1(n-l)-m_1(2)$ columns are indexed by $q^k_{i,j}$ with
${i=m_k+1\div m_k(n-l-1)}$, ${j=1\div m_k}$, ${k=1\div t}$;
\item[--] then,  $m_1(2)$ columns indexed by $q^k_{i,j}$ with
$i=1\div m_k$, $j=1\div m_k$, $k=1\div t$;
\item[--] then, $m_1(l+2)-m_1(2)$ columns indexed by $q^k_{i,j}$ with
$i=1\div m_k$, $j=m_k+1\div m_k(l+1)$, $k=1\div t$;
\item[--] the remaining columns correspond to
$i=m_k+1\div m_k(n-l-1)$, $j=m_k+1\div m_k(l+1)$, $k=1\div t$.
\end{itemize}
%\private{\tp{I've corrected the numbers of columns in the previous list. I realized yesterday they were all wrong. Still I would like someone to verify them. }}

Notice that the equations in $\eq_{A}$ only have coordinates in the first $m_1(n-l)$ columns and $\eq_{B}$ only in the middle $m_1(l+2)$ columns.
With this ordering, the jacobian matrix has the form
\begin{eqnarray*}
\newcommand*{\temp}{\multicolumn{1}{c|}{0}}
\renewcommand{\arraystretch}{2}
\J{A|B}(T)=\left (\; \begin{array}{ccccccc}
\cmidrule[2pt]{1-4}
\multicolumn{4}{|c|}{\J{A| L_A}(T_A)  } & 0 & \multicolumn{2}{c}{\quad 0}  \\
\cmidrule[2pt]{1-5}
\multicolumn{3}{c}{\quad * \qquad \,}  &\multicolumn{2}{|c|}{\J{L_B | B}(T_B)}& \multicolumn{2}{c}{\quad 0}  \\
\cmidrule[2pt]{1-7}
%0 &0 &\temp & 1 &\ast &\ast \\
\multicolumn{7}{|c|}{  \frac{\partial }{\partial \, q^k_{i,j} }E^k_{i,j}}
\\ \cmidrule[2pt]{1-7}
\end{array}\; \right )
\end{eqnarray*}%\cline{1-7}
where:
\begin{enumerate}
 \item The first block $\J{A|L_A}(T_A) $ has $\codim(CV_{T_A})$ rows, and $m_1(n-l)$ columns indexed by the coordinates in $(\o^{n-l} W)^G$ extended to $(\o^{n} W)^G$.
% : rows are indexed by the equations $\eq_A^n$, and columns by the variables $\{q^k_{i,j}\}$, with $i=1\div m_k(n-2)$, $j=1\div m_k$, and ordered so that the coordinates $q^k_{i,j}$ with $i=1\div m_k$ and $j=1\div m_k$ are left to the end.
 %
 \item The second block $\J{L_B | B}(T_B)$ has  $\codim(CV_{T_B})$ rows, and $m_1(l+2)$ columns indexed by the coordinates in $(\o^{l+2} W)^G$ extended to $(\o^{n} W)^G$.
%  Rows are indexed by the equations $\eq_B^n$ and columns by the coordinates $\{q^k_{i,j}\}$ with $i=1\div m_k$, $j=1\div m_k(2)$ so that the coordinates $q^k_{i,j}$ with $i=1$.
\end{enumerate}
Notice that the first two blocks share the columns indexed by the coordinates $\{q^k_{i,j}\}_{1\leq i,j \leq m_k}$.

\begin{enumerate}
\setcounter{enumi}{2}
\item The third block $\big ( \frac{\partial }{\partial \, q^k_{i,j} }E^k_{i,j} \big )$ has
 $N_{n-l-1| l+1}$ rows, indexed by the equations of $\eq^n_{A| B}$, and $m_1(n)$ columns, indexed by all the coordinates above:  $\{q^k_{i,j}\}_{1\leq i\leq m_k(l+1), 1\leq j \leq m_k(n-l-1), 1\leq k \leq t}$.
\end{enumerate}

From now on, these coordinates refer to a generic point of no evolution $\pi_n.$

We proceed by induction on the number of leaves and the induction hypothesis applied is the one explained above the statement of the theorem.
By the induction hypothesis we know that the rank of $\J{A\mid L_A}(T_A)$ is equal to $\codim(CV_{T_A})$. By the claw $d$-tree hypothesis (\ref{assumption_eqn}), we know that
\begin{eqnarray*}
 \J{L_B| B}(T_B)=\left[
 \begin{array}{cc|r}
 * & * & \\
 * & * & \J{L_B|B}^*(T_B) \\
 * & *  &
\end{array}
\right]
\end{eqnarray*}
and the rank of $\J{L_B|B}^*(T_B)$ is equal to $\codim(CV_{T_B})$. This assures that the first $\codim(CV_{T_B})+\codim(CV_{T_A})$ rows in the matrix are linearly independent. %
For the third block and by virtue of Lemma \ref{derivatives}, we have
\begin{eqnarray*}
 \Big [ \frac{\partial }{\partial \, q^k_{i,j} }E^k_{i,j} \Big ] =
 \left[
 \begin{array}{cccc}
 * & \ldots &\ldots & * \\
 * & & & *  \\
 * & \ldots &\ldots & *
\end{array}
\right | \left .
 \begin{array}{c}
\; \\
 \operatorname{Diag}(\Delta_k(\pi_n)) \\
\;
 \end{array}
 \right]
\end{eqnarray*}
% \begin{eqnarray*}
%  \Big ( \frac{\partial }{\partial \, q^k_{i,j} }E^k_{i,j} \Big ) =
%  \left[
%  \begin{array}{cccc}
%  * & \ldots &\ldots & * \\
%  * & & & *  \\
%  * & \ldots &\ldots & *
% \end{array}
% \right | \left .
%  \begin{array}{cccc}
% \Delta_k(\pi_n) & 0 & \ldots & 0  \\
% 0 & \Delta_k(\pi_n) & \ldots & 0   \\
% \ldots &  & & \ldots  \\
% 0 & 0 & \ldots & \Delta_k(\pi_n)
%  \end{array}
%  \right]
% \end{eqnarray*}
where $\operatorname{Diag}(\Delta_k(\pi_n))$ is the diagonal matrix with entries $\{\Delta_k(\pi_n)\}_k$ and columns indexed by the coordinates $q^k_{i,j}$ for $i,j\geq m_k+1$. By Lemma \ref{nonzero_derivatives} below, these entries are nonzero. We conclude that the rank of $\J{A|B}(T)$ is maximal and equal to $\codim(CV_{T})$.
\end{proof}

\begin{lem}\label{nonzero_derivatives}
Let $\pi_n$ be a generic point of no evolution.
Then, $\Delta_k(\pi_n)\neq 0$ for every $k=1,\ldots, t$.
\end{lem}

%\paragraph{Coordinates of no evolution points}
\begin{proof}
The matrix $M^0_k$ for any tensor in $p\in (\o^n W)^G$ represents a tensor in $\FF{k^*}{W_{\alpha}} \o \FF{k}{W_{\beta}}$ obtained as follows:
\begin{enumerate}
\item first contract $p$ with $f_{L(T)\setminus \{\alpha,\beta\}}$ to obtain a tensor $p'$ in $(W_{\alpha}\o W_{\beta})^G$,
\item project $p'$ according to the decomposition $(W_{\alpha}\o W_{\beta})^G\cong \bigoplus_k  \FF{k^*}{W_{\alpha}}\o \FF{k}{W_{\beta}}$.
\end{enumerate}
The determinant of $M_k^0$ is nonzero if and only if the associated map $\FF{k}{W_{\alpha}^*} \rightarrow \FF{k}{W_{\beta}}$ has maximal rank, i.e. is an isomorphism. However, this is the case for all $k$ if and only if $p'$ defines a $G$-isomorphism $W_{\alpha}^*\rightarrow W_{\beta}$.

By lemma \ref{marginal_noevol}, the marginalization of  $\pi_n=\sum_{\xx \in \Sigma} \pi_{\xx} \xx  \o \ldots  \o \xx$ over all leaves different from $\alpha$ and $\beta$ provides a tensor $p'=\sum_{{\xx}\in \Sigma} \pi_{\xx} {\xx}\o {\xx}$. This tensor corresponds to the map from $W_{\alpha}^*$ to $W_{\beta}$ whose matrix in basis $\{X_i\}$ is diagonal with entries $\pi_{X_i}$. Therefore, if the coordinates $\pi_X$ of $\pi_n$ are all non-zero, this is an isomorphism and therefore the matrices $M^0_k$ have non-zero determinant for all $k$.
%and hence it clearly defines an isomorphism, mapping a basis to a dual basis, between .
%
 This proves the claim.
\end{proof}

%\begin{rem}
% \red{From this lemma, can we infer that a point of no evolution $\pi_n=\sum_{\xx \in \Sigma} \pi_{\xx} \xx \o  \ldots \o  \xx$ is  \textit{generic} in our sense if it has coordinates $\pi_{\xx}$ different from 0. But in order to obtain a complete intersection for the tripod this condition is not enough (see Remark \ref{rem_generic}). (MM: My guess is that if will be often enough). Marta: I agree. Let's remove this remark!}
%\end{rem}

\begin{rem}\label{rem_notexpon}
The list of equations for a phylogenetic tree as in Theorem \ref{invariantsCI} is obtained from edge invariants and from local equations for $d$-claw trees, $d \in D$.  As pointed out in \ref{rem_Nab}, the number of edge invariants  is exponential in $n$ (at least for the usual models) and can be substituted by a direct evaluation of the rank of the thin flattening. It is worth noticing that the other subset of equations, the ones coming from $d$-claw trees, are at most exponential in $d$ and therefore this subset can reasonably be used in practice.
\end{rem}

\section{Explicit equations for usual models}

The aim of this section is to provide explicit examples of complete intersections of the particular models listed in Example \ref{main_equivariantmodels}. For the Kimura 3-parameter, this was already done in \cite{casfer2008}; here we deal with GMM, strand symmetric and Jukes-Cantor models (the only remaining case would be Kimura 2-parameter, for which the tripod assumption could be checked using computational algebra software).

%
%\begin{notation}
%\begin{itemize}
%%  \item Write $B_i$ for a basis of $W^{\chi_i}$.
%% \item $B_i \o B_i:=\{u\o v: u,v\in B_i\}$
%%
%% \item Write $B^{(2)}:=\cup_i B_i \o B_i$. This is a basis for $\opG{2}$. \\
%% $B^{(n)}$ is a basis of $\opG{n}$.
%
%\item $\CC^{q_{B^{(n)}}}=\CC^{\{q_u\}_{u\in B^{(n)}}}$
%
%\item $\J{\BB{n}}{E_{T}}$ is the jacobian matrix for equations $E_{T}$ in basis $\BB{n}$.
%\end{itemize}
%\end{notation}

As mentioned in Example \ref{examples2}, for these models we can use the Fourier basis to span the isotypic components of $W$. In this cases, we can identify $\underline{\Sigma}$ with the group $(\mathbb{Z}_2 \times \mathbb{Z}_2,+)$ via
\begin{eqnarray}
 \label{id_group}
 \baa \mapsto (0,0); \quad  \bcc \mapsto (1,0); \quad
 \bgg\mapsto (0,1); \quad  \btt \mapsto (1,1).
\end{eqnarray}
%The characters of the irreducible representations of $G$ will play an important role in Section 7, and we will  denote them as $\chi_{\aa}=\chi_1$, $\chi_{\cc}=\chi_2$, $\chi_{\gg}=\chi_3$ and $\chi_{\tt}=\chi_4$.
We denote by $\chi_{\aa}$, $\chi_{\cc}$, $\chi_{\gg}$, $\chi_{\tt}$ the characters associated to this group $(\mathbb{Z}_2 \times \mathbb{Z}_2,+)$ (see table \ref{tableFourier}). %
These characters are useful to describe the coordinates of a point of no evolution.

\begin{table}[hb]
\begin{center}
\begin{tabular}{ l || c |  c |  c |  c }  & $\baa$ & $\bcc$ & $\bgg$ & $\btt$ \\ \hline
$\chi_{\aa}$ & 1 & 1 & 1 & 1  \\
$\chi_{\cc}$ & 1 & 1 & -1 & -1  \\
$\chi_{\gg}$ & 1 & -1 & 1 & -1 \\
$\chi_{\tt}$ & 1 & -1 & -1 & 1 \\
  \end{tabular} \quad
  \end{center}
\caption{\label{tableFourier} Description of characters $\chi_{\aa}$, $\chi_{\cc}$, $\chi_{\gg}$, $\chi_{\tt}$ for the group $(\mathbb{Z}_2 \times \mathbb{Z}_2,+)$.}
\end{table}

\begin{lem}\label{Fourier_coord}
 If $\pi_n$ is a point of no evolution, then
 $ \pi_n=\sum_{\byy_1,\ldots,\byy_n} q_{\byy_1,\ldots, \byy_n} \byy_1 \ldots  \byy_n$,
 where
 \begin{eqnarray*}
  q_{\byy_1,\ldots, \byy_n}=\frac{1}{4^n} \left (\pi_{\aa} +
  \chi_{\cc}(\byy)\pi_{\cc} +
  \chi_{\gg}(\byy)\pi_{\gg}+
  \chi_{\tt}(\byy)\pi_{\tt} \right ),
 \end{eqnarray*}
 $\byy:=\byy_1+\ldots+\byy_n$ with the operation given by the identification (\ref{id_group}).
%  \red{But we have underlined X's in the q's !!!} \\
%  \tp{Characters are applied to elements of the group $Z_2\times Z_2\cong \Sigma$; the relevant equality is $\xx=\frac{1}{4}\sum_{\byy\in \underline{\Sigma}} \chi_{\xx}(\yy) \byy$. Anyway, there is a well-stablished bijection $\aa \leftrightarrow \baa$, $\cc \leftrightarrow \bcc$, etc.  }\\
%  \tp{Another possibility would be to identify $Z_2\times Z_2$ with $\underline{\Sigma}$ and define the sum directly there. Then $\xx=\frac{1}{4}\sum_{\byy\in \underline{\Sigma}} \chi_{\xx}(\byy) \byy$ }
%
%
\end{lem}
\begin{proof}
From the definition of the Fourier basis in 4.1.1 and using Table \ref{tableFourier}, we have $\xx=\frac{1}{4}\sum_{\byy\in \underline{\Sigma}} \chi_{\xx}(\byy)\, \byy$, for any $\xx \in \Sigma$.
%where $\yy$ is identified with the element of $(\mathbb{Z}_2 \times \mathbb{Z}_2,+)$ according to (\ref{id_group}). %
Now, using that characters of 1-dimensional representations are multiplicative, if $\pi_n$ is a point of no evolution we have
\begin{eqnarray*}
 \pi_n & = &  \sum_{\xx} \pi_{\xx} \xx \o \ldots \o \xx = \\
 & = & \sum_{\xx} \; \frac{\pi_{\xx}}{4^n} \; \Big ( \sum_{\byy_1} \chi_{\xx}(\byy_1) \,\byy_1 \Big ) \o \ldots \o \Big ( \sum_{\byy_n} \chi_{\xx}(\byy_n)\, \byy_n \Big ) = \\
 & = & \sum_{\xx} \; \frac{\pi_{\xx}}{4^n} \; \sum_{\byy_1,\byy_2,\ldots, \byy_n} \chi_{\xx}(\byy_1)\chi_{\xx}(\byy_2)\ldots \chi_{\xx}(\byy_n)\;  \byy_1 \, \byy_2 \, \ldots \, \byy_n  = \\
 & = & \sum_{\xx} \; \frac{\pi_{\xx}}{4^n} \; \sum_{\byy_1,\byy_2,\ldots, \byy_n} \chi_{\xx}(\byy_1+\byy_2+\ldots + \byy_n)\;  \byy_1 \, \byy_2 \, \ldots \, \byy_n  = \\
  & = & \sum_{\byy_1,\byy_2,\ldots, \byy_n} \left ( \frac{1}{4^n} \sum_{\xx} \; \pi_{\xx} \;  \chi_{\xx}(\byy) \right ) \;  \byy_1 \, \byy_2 \, \ldots \, \byy_n.
\end{eqnarray*}
From this, the claim follows.
\end{proof}

%\textbf{Assumption}
%There is a basis $\BB{3}$ of $\opG{3}$ obtained by extending $\BB{2}$ and a system of $c(3)$ equations in $\CC^{q_{\BB{3}}}$ such that the rank of the matrix ... is maximal, equal to $c(3)$.
%
%% Consider the following basis for $\opG{2}$: $\{u_i\o v_i\}_i$, where $u_i, v_i$ are in a basis of $W^{\chi_i}$.
%
%\textbf{HI}
%There is a basis $\BB{n-1}$ of $\opG{n-1}$ and a system of equations in $\CC^{q_{\BB{n-1}}}$ such that the rank of $J_{\BB{n-1}}(E_{T_{n-1}})$ is maximal (equal to $c(n-1)$).
%
%\textbf{Claim}
%There is a basis $\BB{n}$ of $\opG{n}$ and a system of equations in $\CC^{q_{\BB{n}}}$ such that the rank of $J_{\BB{n}}(E_{T_n})$ is maximal (equal to $c(n)$).

% Given an inner product in $E$, we have an isomorphism
% \begin{eqnarray}\label{aux}
%  \Hom_G(E,F) & \cong & (E \o F)^G \\
%  f & \mapsto  & \sum_i u_i\o f(u_i) \nonumber
% \end{eqnarray}
% where $\{u_i\}$ is a given orthonormal basis of $E$. Moreover, two tensors $\tau_1=\sum_i u_i\o f(u_i)$ and $\tau_2=\sum_i v_i\o f(v_i)$ are equal if and only if there is an orthogonal (unitary?) map such that $f(u_i)=v_i$.

% Throughout this section, we will adopt the following assumption: we will assume $G$ is a permutation group $G\leq \mathfrak{S}_4$ such that
% $W=\oplus_{k=1}^t W^{\chi_k}$, where each $W^{\chi_k}$ decomposes as direct sum of $G$-modules generated by elements of the Fourier basis, i.e. for each $k$, we have
% \begin{eqnarray*}
% W^{\chi_k}=\oplus_{i=1}^{m_k} [ \bxx^k_{i,1},\ldots, \bxx^k_{i,d_k}
% ],
% \end{eqnarray*}
% and each $[ \bxx^k_{i,1},\ldots, \bxx^k_{i,d_k} ]$ being $G$-invariant.

\begin{figure}
\begin{center}
%EPS problem
\includegraphics[scale=0.5]{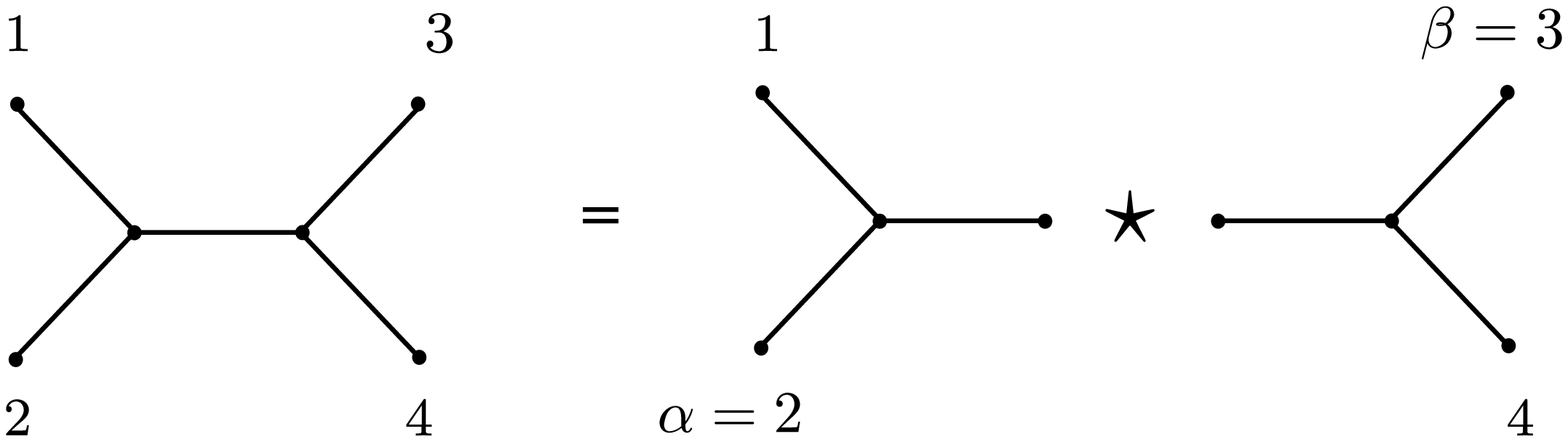}
\end{center}
\caption{\label{figure:quartet}.}
\end{figure}

%\private{
%\tp{The structure of these two examples is the same. First, we proceed to obtain the basis $\omega^k_{ij}$ of $(\o^3 W)^G$ to check the assumption \ref{assumption}. Then, we obtain the basis of $(\o^4 W)^G$ used to write the equations for the complete intersection of the quartet. }
%}

\subsection{Equations for trees evolving under the Jukes-Cantor model}

Take the Jukes-Cantor model, this is, $G=\mathfrak{S}_4$.
There are five irreducible representations of the group $G$: $\N{1}, \dots,\N{5}$ (see \cite{fulton1991} $\S$2.3). For each $i=1\div 5$, the representation $\N{i}$ has the character $\chi_i$ shown in Table \ref{table:char_tables}.
% \begin{center}  \begin{tabular}{ l || c |  c |  c |  c | c  } $\Omega_{\mathfrak{S}_4}$ & $\mathrm{id}$ & $(\aa\cc)$ & $(\aa\cc\gg)$ & $(\aa\cc\gg\tt)$ & $(\aa\cc)(\gg\tt)$ \\ \hline
% $\chi_1$ & 1 & 1 & 1 & 1 & 1 \\
% $\chi_2$ & 1 & -1 & 1 & -1 & 1 \\
% $\chi_3$ & 2 & 0 & -1 & 0 & 2 \\
% $\chi_4$ & 3 & 1 & 0 & -1 & -1 \\
% $\chi_5$ & 3 & -1 & 0 & 1 &- 1 \\
% \hline $\chi$     & 4 & 2 & 1 & 0 & 0
%   \end{tabular}
% \end{center}
As in Example \ref{examples2}.5, we have $u^1_1=\baa$, $u^4_1=\bcc$ and $\FF{1}{W}= \langle\baa\rangle$ and $\FF{4}{W}=\langle\bcc \rangle$. By letting the group $G$ act, it follows that $W=W[\chi_1] \oplus W[\chi_4]$, with $W[\chi_1]=\langle \baa\rangle$ and $W[\chi_4]=\langle \bcc,\bgg,\btt \rangle$. %

We take the tree $T$ with 4 leaves $12|34$ and take the bipartition $A=\{1,2\}$ and $B=\{3,4\}$.
We proceed to obtain a basis for the ambient space $(\o^4 W)^G$ of the variety $V_T$. Choose $\alpha=2$ and $\beta=3$ as in figure~\ref{figure:quartet}.
Following the algorithm described, and noting that $\chi^2=2\chi_1+\chi_3+3\chi_4+\chi_5$, we obtain the following basis
\begin{eqnarray*}
\begin{tabular}{lcl}
 $u^1_{A,1} = \baa\baa $ & \quad & $u^1_{B,1} = \baa\baa $ \\
 $u^1_{A,2} =   \bcc\bcc+\bgg\bgg+\btt\btt  $ & \quad &  $u^1_{B,2} =   \bcc\bcc+\bgg\bgg+\btt\btt $\\
 $u^3_{A,1} =  \bcc\bcc-\bgg\bgg  $ & \quad & $u^3_{B,1} =  \bcc\bcc-\bgg\bgg $\\
 $u^4_{A,1} =  \baa \bcc$ & \quad & $u^4_{B,1} =  \bcc\baa$ \\
 $u^4_{A,2} =  \bcc  \baa$ & \quad & $u^4_{B,2} =   \baa\bcc$\\
 $u^4_{A,3} =  \bgg\btt+\btt\bgg $& \quad &  $u^4_{B,3} =  \bgg\btt+\bgg\btt $\\
 $u^5_{A,1} =   \bgg\btt-\btt\bgg $& \quad &  $u^5_{B,1} =   \bgg\btt-\btt\bgg$
 \end{tabular}
\end{eqnarray*}

Now, for each of the irreducible representations of $\mathfrak{S}_4$, we can choose a subgroup $H_k$ so that $\{hu^k_{A,i} \mid h\in H_k\}$ is a basis of the representation spanned by the $\mathfrak{S}_4$-orbit of $u^k_{A,i}$. Namely,
\begin{eqnarray*}
 H_1 = \{e\} & \quad & n_1=1; \\
 H_3 = \{e,(\aa\cc)\} & \quad & n_3=2;\\
 H_4 = \{e,(\cc\gg\tt),(\cc\tt\gg) \} & \quad & n_4=3; \\
 H_5 = \{e,(\cc\gg\tt),(\cc\tt\gg) \} & \quad & n_5=3.
\end{eqnarray*}
A basis for $(W_A\o W_B)^G$
is inferred  from $\bigoplus_k \FF{k}{W_A}\o \FF{k}{W_B}$ by taking the image of the operator $S$ specified in Remark \ref{S_operator}  %$S=\frac{1}{|G|} \sum_{g\in G} \rho_A(G)\o \rho_B(g)$
applied to the tensors $u^k_{A,i}\o u^k_{B,j}$:
\begin{eqnarray*}
%S(u^k_{A,i}\o u^k_{B_j}):=
S(u^k_{A,i}\o u^k_{B,j}):=
\frac{n_k}{|G|}\sum_{g\in G} (g\cdot u^k_{A,i}) \o (g\cdot u^k_{B_j}),
\end{eqnarray*}
that is,
\begin{eqnarray*}
S(u^1_{A,1} \o u^1_{B,1}) & = & \baa \baa \baa \baa, \\
S(u^1_{A,1} \o u^1_{B,2})& = & \frac{1}{3} \baa \baa \o (\bcc\bcc+\bgg\bgg+\btt\btt), \\
S(u^1_{A,2} \o u^1_{B,1})& = & \frac{1}{3} (\bcc\bcc+\bgg\bgg+\btt\btt )\o \baa \baa,  \\
S(u^1_{A,2} \o u^1_{B,2})& = & \frac{1}{3} (\bcc\bcc+\bgg\bgg+\btt\btt)\o (\bcc\bcc+\bgg\bgg+\btt\btt), \\
S(u^3_{A,1} \o u^3_{B,1})& = &\frac{2}{3}\big (
(\bcc\bcc-\bgg\bgg)\o(\bcc\bcc-\bgg\bgg)+
(\bgg\bgg-\btt\btt)\o(\bgg\bgg-\btt\btt)+
(\btt\btt-\bcc\bcc)\o(\btt\btt-\bcc\bcc) \big ) \\
S(u^4_{A,1} \o u^4_{B,1})& = &\baa\bcc\bcc\baa+\baa\bgg\bgg\baa+\baa\btt\btt\baa, \\
S(u^4_{A,1} \o u^4_{B,2})& = &\baa\bcc\baa\bcc+\baa\bgg\baa\bgg+\baa\btt\baa\btt, \\
S(u^4_{A,1} \o u^4_{B,3})& = &\baa\bcc\o(\bgg\btt+\btt\bgg)+\baa\bgg\o(\bcc\btt+\btt\bcc)+\baa\btt\o(\bcc\bgg+\bgg\bcc), \\
S(u^4_{A,2} \o u^4_{B,1})& = &  \bcc\baa\bcc\baa+\bgg\baa\bgg\baa+\btt\baa\btt\baa , \\
S(u^4_{A,2} \o u^4_{B,2})& = & \bcc\baa\baa\bcc+\bgg\baa\baa\bgg+\btt\baa\baa\btt,\\
S(u^4_{A,2} \o u^4_{B,3})& = &  \bcc\baa\o(\bgg\btt+\btt\bgg)+\bgg\baa\o(\bcc\btt+\btt\bcc)+\btt\baa\o(\bcc\bgg+\bgg\bcc), \\
S(u^4_{A,4} \o u^4_{B,1})& = &  (\bgg\btt+\btt\bgg)\o\bcc\baa+(\bcc\btt+\btt\bcc)\o\bgg\baa+(\bcc\bgg+\bgg\bcc)\o\btt\baa, \\
S(u^4_{A,3} \o u^4_{B,2})& = &   (\bgg\btt+\btt\bgg)\o\baa\bcc+(\bcc\btt+\btt\bcc)\o\baa\bgg+(\bcc\bgg+\bgg\bcc)\o\baa\btt,\\
S(u^4_{A,3} \o u^4_{B,3})& = & (\bgg\btt+\btt\bgg)\o(\bgg\btt+\btt\bgg)+(\bcc\btt+\btt\bcc)\o(\bcc\btt+\btt\bcc)+
(\bcc\bgg+\bgg\bcc)\o(\bcc\bgg+\bgg\bcc), \\
S(u^5_{A,1} \o u^5_{B,1})& = & (\bgg\btt-\btt\bgg)\o(\bgg\btt-\btt\bgg)+(\bcc\btt-\btt\bcc)\o(\bcc\btt-\btt\bcc)+(\bcc\bgg-\bgg\bcc)\o(\bcc\bgg-\bgg\bcc).
%
% (W_A)_1\o (W_B)_1& = & [ u^1_{A,1} \o u^1_{B,1}, u^1_{A,1} \o u^1_{B,_1}, u^1_{A,1} \o u^1_{B,1}, u^1_{A,1} \o u^1_{B,1}] \\
%  (W_A)_3\o (W_B)_3& = &[u^3_{A,1}\o u^3_{B,1}] \\
%  (W_A)_4\o (W_B)_4& = &[ u^4_{A,1} \o u^4_{B,1}, u^4_{A,1} \o u^4_{B,2},u^4_{A,1} \o u^4_{B,3}, u^4_{A,2} \o u^4_{B,1}, u^4_{A,2} \o u^4_{B,2},\ldots \\ & &
%  \ldots, u^4_{A,2} \o u^4_{B,3},
%  u^4_{A,3} \o u^4_{B,1}, u^4_{A,3} \o u^4_{B,2},u^4_{A,3} \o u^4_{B,3}] \\
%  (W_A)_5\o (W_B)_5& = &[u^5_{A,1}\o u^5_{B,1}]
\end{eqnarray*}
As above, denote by $q^k_{i,j}$ the coordinates corresponding to this  basis $S(u^k_{A,i} \o u^k_{B,j})$.
We proceed to obtain a complete intersection for the tree $T$ with 4 leaves $12|34$. Take the bipartition $A=\{1,2\}$ and $B=\{3,4\}$.
%A convenient basis for the space $\o^4 W$ has been obtained in example \ref{}.

First of all, we proceed to obtain the edge invariants following the section 4.2.
If $\pi_4$ is a no evolution point, write $\pi_4=\sum_{k,i,j} q^k_{ij} \, S(u^k_{A,i} \o u^k_{B,j})$. Each irreducible representation $N_k$ of $\mathfrak{S}_4$ gives rise to a $m_k(2)\times m_k(2)$-matrix $M_k$:
\begin{eqnarray*}
 M_1=\left ( \begin{array}{cc}
q^1_{11} & q^1_{12} \\
q^1_{21} & q^1_{22}
\end{array}
\right ), \qquad
M_3=\left ( \begin{array}{c}
q^3_{11}
\end{array}
\right ), \qquad
M_4=\left ( \begin{array}{ccc}
q^4_{11} & q^4_{12} & q^4_{13}\\
q^4_{21} & q^4_{22} & q^4_{23}\\
q^4_{31} & q^4_{32} & q^4_{33}
\end{array}
\right ), \qquad
M_5=\left ( \begin{array}{c}
q^5_{11}
\end{array}
\right ) \quad
\end{eqnarray*}
where the rows of each $M_k$ are indexed by the $\{u^k_{A,i}\}$ and the columns are indexed by the $\{u^k_{B,j}\}$. Notice that there is no $M_2$ as there is no isotypic component corresponding to $\N{2}$ in $\o^4 W$. %
Moreover, from Lemma \ref{nonzero_derivatives}, we know that
$\Delta_1(\pi_4) = q^1_{11} \neq 0$ and $\Delta_4(\pi_4) =  q^4_{11}\neq 0$.
The resulting edge invariants arise as rank restrictions for these matrices:
\begin{center}
%\begin{table}
\begin{tabular}{lcl}
$\chi_1$  :  $q^1_{2,2}q^1_{1,1}-q^1_{1,2}q^1_{2,1}=0$ & $\quad$ &
$\chi_4$  :   $q^4_{2,2}q^4_{1,1}-q^4_{1,2}q^4_{2,1}=0$  \\
 $\chi_3$  :  $q^3_{1,1}=0$ & $\quad$  & $\qquad$ $q^4_{2,3}q^4_{1,1}-q^4_{1,3}q^4_{2,1}=0$ \\
$\chi_5$  :   $q^5_{1,1}=0$ & $\quad$ &  $\qquad$ $q^4_{3,2}q^4_{1,1}-q^4_{1,2}q^4_{3,1}=0$ \\
  & $\quad$ &  $\qquad$  $q^4_{3,3}q^4_{1,1}-q^4_{1,3}q^4_{3,1}=0$
\end{tabular}
  %\end{table}
 \end{center}
We also need to consider the equations obtained from the tripods associated to the bipartition: $T_A$ with leaves $\{1,\alpha,L_A\}$, and $T_B$ with leaves $\{L_B,\beta,4\}$ (see figure \ref{figure:quartet}). %
Now, it can be seen that a complete intersection for the tripod $T_3$ with leaves $x,y,z$ is given by
% \begin{eqnarray*}
%  q^4_{1,1}q^4_{2,1}q^1_{2,1}-q^1_{1,1}(q^4_{3,1})^2=0
% \end{eqnarray*}

\begin{eqnarray*}
  Q^4_{1,1}Q^4_{1,2}Q^1_{1,2}-Q^1_{1,1}(Q^4_{1,3})^2=0
 \end{eqnarray*}

where these $Q^k_{ij}$ are the coordinates corresponding to the basis linked to the edge split $x|yz$:
\begin{eqnarray*}
 %\omega^1_{1,1}
 S(u^1_{x,1} \o u^1_{yz,1}) & = & \baa \baa \baa, \\
 %\omega^1_{2,1}
 S(u^1_{x,1} \o u^1_{yz,2}) & = & \baa\bcc\bcc +\baa\bgg \bgg +\baa\btt \btt, \\
 %\omega^4_{1,1}
 S(u^4_{x,1} \o u^4_{yz,1}) & = & \bcc \baa \bcc +\bgg\baa \bgg + \btt \baa \btt, \\
 %\omega^4_{1,1}
 S(u^4_{x,1} \o u^4_{yz,2}) & = & \bcc \bcc \baa+\bgg \bgg \baa+ \btt \btt \baa, \\
 %\omega^4_{1,1}
S(u^4_{x,1} \o u^4_{yz,3})  & = & \bcc \bgg \btt +\bcc \btt \bgg +\bgg \bcc \btt + \bgg \btt \bcc + \btt  \bcc \bgg +  \btt \bgg \bcc.
\end{eqnarray*}

  From $V_{T_A}$, we take $x=L_A,y=1,z=2$ to obtain the extended equation in the original coordinates $q^k_{ij}$:
  \begin{eqnarray*}
 q^4_{1,1}\,q^4_{2,1}\,q^1_{1,2}-q^1_{1,1}\,(q^4_{3,1})^2=0.
   \end{eqnarray*}
Analogously, from $V_{T_B}$, we take $x=L_B,y=3,z=4$ to obtain the extended equation:
 \begin{eqnarray*}
  q^4_{1,1}\, q^4_{1,2}\, q^1_{1,2}-q^1_{1,1}\,(q^4_{1,3})^2=0.
 \end{eqnarray*}
These 9 equations
%together with the stochastic restriction $\sum_{k,i,j} q^k_{i,j}=1$
define a local complete intersection around generic points of no evolution for the variety $CV_G(T)$ in $(\o^4 W)^G$ .

\subsection{Equations for the tripod evolving under the general Markov model}\label{sec:gmm}
The aim of this subsection is to explicitly provide codimension many equations of the variety $X$ for the tripod for GMM that cut $X$ out in a neighborhood of no evolution points. As the salmon conjecture is well-studied and answered on the set theoretic level \cite{MR2097214,MR2836258,MR2891138, MR2500471} many of the equations of $X$ are known. However, it turns out that the simplest equations, going back to Strassen, are enough to obtain a description of the local complete intersection. As we will see not only are they enough - also their choice is astonishingly natural.

Recall that when $\Sigma$ has $\kappa$ elements, $X$ is the $\kappa$-th secant variety of the Segre product $\p(A)\times\p(B)\times\p(C)$, i.e. the closure of the locus of rank $\kappa$ tensors in the space $\p(A\otimes B\otimes C)$, where $\dim A=\dim B=\dim C=\kappa$.
%\red{To be consistent with the previous sections, we should change $e_i$ by $\xx_i$}.

For the sake of completeness let us recall Landsberg-Ottaviani's interpretation of Strassen equations \cite{MR3376667}.
Each tensor $\tau\in A\otimes B\otimes C$ is naturally identified with a map $\tau:A^*\rightarrow B\otimes C$. Let us tensor this map with the identity on $C$ obtaining a map $A^*\otimes C\rightarrow B\otimes C\otimes C$. Using the natural map $C\otimes C\rightarrow C\wedge C$ we obtain $f(\tau):A^*\otimes C\rightarrow B\otimes (C\wedge C)$.

When $\tau$ has rank one, i.e. $\tau=a\otimes b\otimes c$, then the rank of $f(\tau)$ (as a matrix) is at most $\kappa-1$ because the image of $f(\tau)$ is the subspace $b\otimes (c \wedge C)$. Hence, if $\tau$ is of rank $\kappa$, then $f(\tau)$ has rank at most $\kappa(\kappa-1)$. Using the matrix representation of $f(\tau)$ all $\kappa(\kappa-1)+1$ minors provide equations of the $\kappa$-th secant variety.

In coordinates, if $\tau=\sum a_{ijk} X_i\otimes X_j\otimes X_k$, in a certain basis $X_1,\dots,X_{\kappa}$ of $W$, then the entry of the column $X_i^*\otimes X_j$ and row $X_s\otimes (X_r\wedge X_t)$ equals (we assume $r<t$):
%\private{\red{Not very important, but I get the signs reversed below.}}
\begin{itemize}
\item $0$, if $r$ and $t$ are different from $j$;
\item $-a_{ist}$, if $r=j$;
\item $a_{isr}$, if $t=j$;
\end{itemize}
because  $f(\tau)(X_i^*\otimes X_j)=\sum_{p,q} a_{ipq}X_p\otimes(X_q\wedge X_j)$.
A display of this matrix is shown in the Table~\ref{matrix}.

%
%
%We obtain the following lemma.

\begin{defn}
In the matrix representation of $f(\tau)$ exactly $\kappa(\kappa-1)$ columns contain an entry $a_{iii}$ or $-a_{iii}$ for some $i$.
Namely, these are  the columns indexed by $X_i^*\otimes X_j$, where $i\neq j$. Each such column contains exactly one such entry.
Moreover, these entries are contained in $\kappa(\kappa-1)$ different rows: those indexed by $X_i\otimes (X_j\wedge X_i)$ (for $j<i$) or $X_i\otimes (X_i\wedge X_j)$ (for $j>i$).
These precise rows and columns will be called \emph{distinguished}.
In table \ref{matrix} distinguished rows and columns are marked with $*$ and depicted in gray.
\end{defn}

Consider any minor $M$ of $f(\tau)$ of order $\kappa(\kappa-1)+1$. If it does not contain all the distinguished rows and columns, then all its derivatives vanish on any point of no evolution. Indeed, each monomial in $M$ will contain at least a degree two factor in variables different from $a_{iii}$, hence any derivative of such monomial (if nonzero) will contain such a variable and will vanish on the no evolution point.

Thus from now on we will be interested only in those minors $M$ that contain all the distinguished rows and columns. Such minors are of course specified by choosing a non-distinguished row $r$ and column $c$.
By the same argument as above, only the derivative of $M$ with respect to the $(r,c)$ entry  can be nonzero at a point of no evolution:
this derivative equals the determinant of the submatrix given by distinguished rows and columns, i.e. it equals
\begin{eqnarray}\label{defn_c}
c:=\pm \left (\prod_i a_{iii} \right )^{\kappa-1}.
\end{eqnarray}
% where $\varepsilon=\pm 1$ and depends only on the variable considered to differentiate.

Let us consider a nondistinguished row indexed by $X_s\otimes (X_r\wedge X_t)$, $r<t$, $r\neq s$, $t\neq s$. It contains exactly two nonzero variables in the nondistinguished columns: -$a_{rst}$ and $a_{tsr}$.
In particular, there are $2(\kappa{\kappa\choose 2}-\kappa(\kappa-1))=\kappa(\kappa-1)(\kappa-2)$ variables in the submatrix indexed by nondistinguished rows and columns, and they are all different.
Hence, the corresponding minors have independent differentials at a generic point of no evolution (which does not have any coordinate equal to zero).
Notice however that $\kappa(\kappa-1)(\kappa-2)=\kappa^3-1-(3\kappa(\kappa-1)+\kappa-1)$ is the codimension of the variety by Terracini's lemma \cite{MR2677616, MR3011324}.
Thus we can conclude that the given minors provide locally a description of the variety as a complete intersection at the generic points of no evolution.

For $r,s,t$ in $\{1,\dots,\kappa\}$, with $s\neq r$, $t\neq r,s$, we call $\eeq_{X_r,X_s,X_t}$ the equation given by the minor formed by the distinguished rows and columns, plus row $X_s\otimes (X_r\wedge X_t)$, and column $X_r^*\otimes X_r$ if $r<t$ or $X_t^*\otimes X_t$ if $t<r$.
We have proven the following:

\begin{lemma0}
The equations $\eeq_{X_r,X_s,X_t}$ for $s\neq r$, $t\neq r,s$ describe the variety $CV_{GMM}(T_3)$ as a complete intersection locally at a generic point of no evolution.
\end{lemma0}

Next, we proceed to prove that Assumption \ref{assumption} is also satisfied. Consider $\kappa=4$ and let $\{\baa,\bcc,\bgg,\btt\}$ be the Fourier basis. %
We have $24$ equations $\eeq_{\xx,\yy,\zz}$ which are indexed by 3-element subsets of $\{\aa,\cc,\gg,\tt\}$.
By the previous discussion, each equation $eq_{\xx,\yy,\zz}$ has only one nonzero directional derivative (after evaluated at a generic point of no evolution) in the standard basis, namely with respect to $\xx\o \yy\o \zz$.
Notice that if we removed all the columns of the Jacobian matrix indexed by variables of type $\aa \otimes \yy \otimes \zz$ the matrix would drop rank: all rows indexed by equations $\eeq_{\aa,\yy,\zz}$ would be zero.

Let us consider the basis $\bxx \otimes \yy \otimes \zz$, where $\bxx \in \{\baa,\bcc,\bgg,\btt\}$, but $\yy,\zz\in \{\aa,\cc,\gg,\tt\}$. % (this is precisely the basis choice in the assumption (except roles of leaves 1 and 3 interchanged), the second and third index do not matter}.
Let us call $J$ the Jacobian matrix with $24$ rows indexed by the above equations written in this new basis
% \red{(Do you mean after substituting variables $a_{\xx \yy \zz}$ by the corresponding combination of the new variables $a_{\bxx '\otimes \yy'\otimes \zz'}$?)
% MM: So I was thinking about equations abstractly - independent of basis - as functions on a vector space.
% }
and $64$ columns indexed by basis elements $\bxx \otimes \yy\otimes \zz$. Let $\tilde J$ be the submatrix of $J$ obtained by removing all the columns indexed by a variable of type $\baa\otimes \yy\otimes \zz$ for any $\yy,\zz$.

\begin{lemma0}\label{lem:maxrank}
The rank of the matrix $\tilde J$ is maximal, i.e. equal to $24$.
\end{lemma0}
\begin{proof}
Let us fix distinct $\yy,\zz\in \{\aa,\cc,\gg,\tt\}$. There are precisely two equations $\eeq_{\xx_1,\yy,\zz}, \eeq_{\xx_2,\yy,\zz}$ as described above.
% \red{what are $X_1$ and $X_2$? MM: I guess I have to be more precise, as in the general discussion the indexing of basis were numbers like $i,j,k$ and here these are letters (in both cases these are indices of basis elements). Here explicitly $\xx_1,\xx_2$ are the two elements of $\Sigma$ different from $\yy,\zz$}.
The evaluation at a generic point of no evolution of the directional derivatives of these equations with respect to  $\underline{\ss}\vv \ww$ equals zero unless $\vv=\yy$ and $\ww=\zz$. These give us $3$ variables that can give nonzero derivatives (as we assume $\underline{\ss}  \neq\baa$).
Actually, the directional derivative of $\eeq_{\xx \yy \zz}$ with respect to $\underline{\ss} \o \vv \o \ww$ is equal to
\begin{eqnarray*}
 \frac{d\, \eeq_{\xx \yy \zz}}{d\, \underline{\ss} \o \vv \o \ww}=
 \left \{ \begin{array}{cc}
\varepsilon  c & \mbox{ if } \vv=\yy, \ww=\zz \\
0 & \mbox{otherwise.}
\end{array}
 \right.
\end{eqnarray*}
where $c$ is the amount defined in (\ref{defn_c}) and $\varepsilon$ represents the sign of $\xx$ in $\underline{\ss}$: $\varepsilon=1$ if $\xx=\aa$ or $\xx=\ss$, and $\varepsilon=-1$ otherwise.

On the other hand, at a generic point of no evolution, the evaluation  of the derivatives with respect to any variable $\underline{\ss}\otimes \yy\otimes \zz$ gives a nonzero value only when applied to the equations $eq_{\xx_1,\yy,\zz}$ or $eq_{\xx_2,\yy,\zz}$. Hence the nonzero entries of the $24\times 48$ matrix $\tilde J$ are contained in $12$ rectangles of shape $2\times 3$, not sharing rows or columns. To finish the proof it remains to show a $2\times 2$ submatrix with nonzero determinant in each rectangle.

If $\yy=\aa$ or $\zz=\aa$, then $\xx_1,\xx_2\neq \aa$. Without loss of generality, we may assume that $\yy=\aa$. Then, we choose the derivatives with respect to $\bxx_1 \yy  \zz$, and $\underline{\zz}  \yy  \zz$. The submatrix obtained has the form:
\[\left( \begin{array}{cc}
            c & -  c \\
           -  c & - c
          \end{array} \right).\]
% where $c\neq 0$ is the derivative of the $\eeq_{\xx \yy \zz}$ with respect to the variable $\xx \otimes \yy \otimes \zz$ evaluated at a generic point of no evolution.
%
If $\yy,\zz\neq \aa$, then $\aa \in \{\xx_1,\xx_2\}$ and we can take $\xx_1=\aa$. We choose the derivatives with respect to $\bxx_2  \yy  \zz$ and $\underline{\yy}  \yy  \zz$. The obtained submatrix is of the form:
\[
\left( \begin{array}{cc}
           c & c \\
           c & - c
          \end{array} \right).\]
In any case, the determinant is $-2c^2\neq 0$ and we are done.
\end{proof}

%If we change the basis arbitrarily on the second and third tensor factor, for example to the Fourier basis, then the matrix $\tilde M$, after removing columns indexed by variables of type $\baa\otimes \dots$, has still maximal rank, by Lemma \ref{lem:maxrank}. Notice, that a priori, it may be hard to point which subminor is nonzero, but this is not needed to verify Assumption \ref{assumption}.

\subsection{Equations for trees evolving under the strand symmetric model}\label{sec:ssm}
Take $G=\langle (\aa \tt)(\cc \gg )\rangle \cong \mathbb{Z}_2$, corresponding to the strand symmetric model as in Example \ref{main_equivariantmodels}.
%
%There are two irreducible representations of the group $G$: $N_1, N_2$. Each representation $N_i$ has the character $\omega_i$ shown in the the following table:
%\begin{center}  \begin{tabular}{ l || c |  c } $\Omega_{\mathbb{Z}_2}$ & $\mathrm{id}$ & $(\aa\tt)(\cc \gg)$ \\ \hline
%$\omega_1$ & 1 & 1  \\
%$\omega_2$ & 1 & -1 \\
%\hline $\chi$  & 4 & 0
%  \end{tabular}
%\end{center}
%The permutation representation of $G$ decomposed as $\chi=2\omega_1+2\omega_2$, and $W=W[\omega_1] \oplus W[\omega_2]$,
%with $W[\omega_1]=\langle \baa, \btt \rangle$ and $W[\omega_2]=\langle \bcc,\bgg\rangle$.

In order to deduce equations for the tripod, we construct first a convenient basis for the space $(\o^3 W)^G$.
We take $\FF{1}{W_{\alpha}}=\langle \baa, \btt \rangle$ and $\FF{2}{W_{\alpha}}=\langle\bcc, \bgg \rangle$.
Take $A=\{1,2\}$ and keep the notation used in Example \ref{example_SS}.
A basis for $(W_A\o W)^G$ would be inferred  from $\bigoplus_k \FF{k}{W_A}\o \FF{k}{W}$ by taking the image of the operator $S$ specified in Remark \ref{S_operator}  %$S=\frac{1}{|G|} \sum_{g\in G} \rho_A(G)\o \rho_B(g)$
applied to the tensors $u^k_{A,i}\o u^k_j$ (as the irreducible representations have dimension 1, Remark \ref{S_operator} trivially applies):
\begin{eqnarray*}
\begin{array}{c}
 S(u^1_{A,1}\o u^1_1) = \baa\baa\baa \\
 S(u^1_{A,2}\o u^1_1) = \baa\btt\baa \\
 S(u^1_{A,3}\o u^1_1) = \btt\baa\baa \\
 S(u^1_{A,4}\o u^1_1) = \btt\btt\baa
 \end{array} \qquad
 \begin{array}{c}
 S(u^1_{A,5}\o u^1_1) = \bcc\bcc\baa  \\
 S(u^1_{A,6}\o u^1_1) = \bcc\bgg\baa \\
 S(u^1_{A,7}\o u^1_1) =\bgg\bcc\baa \\
 S(u^1_{A,8}\o u^1_1) = \bgg\bgg\baa
 \end{array}\qquad
 \begin{array}{c}
   S(u^1_{A,1}\o u^1_2) = \baa\baa\btt \\
 S(u^1_{A,2}\o u^1_2) = \baa\btt\btt \\
 S(u^1_{A,3}\o u^1_2) = \btt\baa\btt \\
 S(u^1_{A,4}\o u^1_2) = \btt\btt\btt
 \end{array}\qquad
 \begin{array}{c}
 S(u^1_{A,5}\o u^1_2) = \bcc\bcc\btt \\
 S(u^1_{A,6}\o u^1_2) = \bcc\bgg\btt \\
 S(u^1_{A,7}\o u^1_2) = \bgg\bcc\btt \\
 S(u^1_{A,8}\o u^1_2) = \bgg\bgg\btt
 \end{array} \\ \\
 \begin{array}{c}
 S(u^2_{A,1}\o u^2_1) = \baa\bcc\bcc \\
 S(u^2_{A,2}\o u^2_1) = \baa\bgg\bcc \\
 S(u^2_{A,3}\o u^2_1) = \btt\bcc\bcc \\
 S(u^2_{A,4}\o u^2_1) = \btt\bgg\bcc
 \end{array} \qquad
 \begin{array}{c}
 S(u^2_{A,5}\o u^2_1) = \bcc\baa\bcc  \\
 S(u^2_{A,6}\o u^2_1) = \bcc\btt\bcc \\
 S(u^2_{A,7}\o u^2_1) =\bgg\baa\bcc \\
 S(u^2_{A,8}\o u^2_1) = \bgg\btt\bcc
 \end{array}\qquad
 \begin{array}{c}
   S(u^2_{A,1}\o u^2_2) = \baa\bcc\bgg \\
 S(u^2_{A,2}\o u^2_2) = \baa\bgg\bgg \\
 S(u^2_{A,3}\o u^2_2) = \btt\bcc\bgg \\
 S(u^2_{A,4}\o u^2_2) = \btt\bgg\bgg
 \end{array}\qquad
 \begin{array}{c}
 S(u^2_{A,5}\o u^2_2) = \bcc\baa\bgg \\
 S(u^2_{A,6}\o u^2_2) = \bcc\btt\bgg \\
 S(u^2_{A,7}\o u^2_2) = \bgg\baa\bgg \\
 S(u^2_{A,8}\o u^2_2) = \bgg\btt\bgg
 \end{array} .
\end{eqnarray*}
In other words, the \textit{Fourier basis for SSM}  is the subbasis of the usual Fourier basis for $\otimes^3W$ formed by triplets that contain an even number of elements in $\{\bcc,\bgg\}$. Thus we denote its coordinates as the usual Fourier coordinates $q_{\underline{\xx}\underline{\yy}\underline{\zz}}$ (corresponding to the basis vector $\underline{\xx}\underline{\yy}\underline{\zz}$).
If $\pi_3=\sum_{\xx} \pi_{\xx} \xx\xx\xx $ is a no evolution point, then its Fourier coordinates are given by

\begin{eqnarray*}
q_{\underline{\xx}\underline{\yy}\underline{\zz}}=\left \{
 \begin{array}{cc}
 2 \pi^+, & \mbox{ if } \bxx+ \byy + \underline{\zz}=\baa \\
 2 \pi^-, & \mbox{ if }\bxx+ \byy + \underline{\zz}=\btt
 \end{array}
 \right .
\end{eqnarray*}
where $\pi^+:=\pi_{\aa}+\pi_{\cc}$ and $\pi^-:=\pi_{\aa}-\pi_{\cc}$ (see Lemma \ref{Fourier_coord}) and the sum of nucleotides is done according to (\ref{id_group}).
%. \tp{I still have to choose 12 equations doing the job. However, discussing with Marta about this, we don't think it's necessary to write them down here in the paper. We could prepare some website to include these equations, and more: the GMM case. What do you think Mateusz? }

A complete intersection for the variety corresponding to the tripod (which has codimension 12 in $(\otimes^3W)^G$ ) is defined by 12 equations and can be obtained from the 24 equations we described for the tripod evolving under GMM as follows.
As in Section \ref{sec:gmm} we consider the tensor $T$ and write the matrix $f(T)$ in the Fourier coordinates. As  $q_{\underline{\xx}\underline{\yy}\underline{\zz}}$ is 0 if $\underline{\xx}\underline{\yy}\underline{\zz}$ contains an odd
number of elements in $\{\bcc,\bgg\}$, the matrix $f(T)$ reduces to the matrix presented in Table \ref{matrixSSM}. The subindices $1,2,3,4$ refer now to $\baa,\bcc,\bgg,\btt$ respectively.
When we consider the same equations as in the previous section we observe that out of the 12 nondistinghuished rows, only 6 contain nonzero entries at the nondistinguished columns (and they contain exactly 2 nonzero entries). These rows are those labelled by
$\underline{\xx}_r\otimes (\underline{\xx}_s\wedge \underline{\xx}_t)$ such that $s<t$ and $\{\underline{\xx}_r,\underline{\xx}_s,\underline{\xx}_t\}$ contains an even number of $\bcc,\bgg 's.$ The same argument used for the general Markov model proves that these twelve $13 \times 13$ minors define a local complete intersection
at the generic points of evolution and that they satisfy assumption \ref{assumption}.
%The $13\times13$ minors considered for GMM by choosing the nondistinguished rows and columns, satisfy a certain symmetry in this new matrix.
%Indeed, considering the permutation $\sigma=(\baa,\btt)(\bcc,\bgg)$ to any  minor formed by the distinguished rows and columns, row $\xx_i\otimes (\xx_j\wedge \xx_k)$, and column $\xx_j^*\otimes \xx_j$ for $j<k$ coincides, up to sign, with the minor formed
%by the distinguished rows and columns, row $\sigma(\xx_i)\otimes (\sigma(\xx_j)\wedge \sigma(\xx_k))$ and column $\sigma(\xx_j)^*\otimes \sigma(\xx_j)$. Therefore, we only need to consider 12 of these minors,
%namely those those given by one representative of these permutation classes.

Now out of these equations we provide a local complete intersection for trees with four leaves. Take $T$ the tree with 4 leaves and choose $\alpha=2$ and $\beta=3$ as in the previous example (see figure~\ref{figure:quartet}). We proceed to obtain a basis for the ambient space $(\o^4 W)^G$ of the variety $V_T$.
Similar computations as above show that
\begin{eqnarray*}
\FF{1}{W_B} = \langle \baa\baa,\btt\baa,\baa\btt,  \btt\btt,\bcc\bcc, \bcc\bgg, \bgg\bcc,\bgg\bgg \rangle;\\
\FF{2}{W_B} = \langle \bcc\baa,\bgg\baa, \bcc\btt, \bgg\btt, \baa\bcc,\btt\bcc, \baa\bgg, \btt\bgg \rangle.
\end{eqnarray*}
So, a basis for $(W_A\o W_B)^G$ is given by all tensors of the form $\omega^k_{i,j}:=u^k_i \o u^k_j$ (as the irreducible representations have dimension 1, Remark \ref{S_operator} trivially applies).

Now, if $\pi_4 \in (\o^n W)^G$ is a no evolution point, and we write
%\begin{eqnarray*}
 $\pi_4=\sum_{k,i,j} q^k_{ij} \, \omega^k_{ij}$,
%\end{eqnarray*}
then each irreducible representation gives rise to a $m_k(2)\times m_k(2)$-matrix $M_k$, $k=1,2$:
\begin{eqnarray*}
 M_k=\left(
 \begin{array}{cccccccc}
  q^k_{11} & q^k_{12} & q^k_{13} & q^k_{14} & q^k_{15} & q^k_{16} & q^k_{17} & q^k_{18} \\
  q^k_{21} & q^k_{22} & q^k_{23} & q^k_{24} & q^k_{25} & q^k_{26} & q^k_{27} & q^k_{28} \\
  q^k_{31} & q^k_{32} & q^k_{33} & q^k_{34} & q^k_{35} & q^k_{36} & q^k_{37} & q^k_{38} \\
  q^k_{41} & q^k_{42} & q^k_{43} & q^k_{44} & q^k_{45} & q^k_{46} & q^k_{47} & q^k_{48} \\
  q^k_{51} & q^k_{52} & q^k_{53} & q^k_{54} & q^k_{55} & q^k_{56} & q^k_{57} & q^k_{58} \\
  q^k_{61} & q^k_{62} & q^k_{63} & q^k_{64} & q^k_{65} & q^k_{66} & q^k_{67} & q^k_{68} \\
  q^k_{71} & q^k_{72} & q^k_{73} & q^k_{74} & q^k_{75} & q^k_{76} & q^k_{77} & q^k_{78} \\
  q^k_{81} & q^k_{82} & q^k_{83} & q^k_{84} & q^k_{85} & q^k_{86} & q^k_{87} & q^k_{88}
 \end{array}
  \right)
%\quad
%  M_2=\left(
%  \begin{array}{cccccccc}
%   q^2_{11} & q^2_{12} & q^2_{13} & q^2_{14} & q^2_{15} & q^2_{16} & q^2_{17} & q^2_{18} \\
%   q^2_{21} & q^2_{22} & q^2_{23} & q^2_{24} & q^2_{25} & q^2_{26} & q^2_{27} & q^2_{28} \\
%   q^2_{31} & q^2_{32} & q^2_{33} & q^2_{34} & q^2_{35} & q^2_{36} & q^2_{37} & q^2_{38} \\
%   q^2_{41} & q^2_{42} & q^2_{43} & q^2_{44} & q^2_{45} & q^2_{46} & q^2_{47} & q^2_{48} \\
%   q^2_{51} & q^2_{52} & q^2_{53} & q^2_{54} & q^2_{55} & q^2_{56} & q^2_{57} & q^2_{58} \\
%   q^2_{61} & q^2_{62} & q^2_{63} & q^2_{64} & q^2_{65} & q^2_{66} & q^2_{67} & q^2_{68} \\
%   q^2_{71} & q^2_{72} & q^2_{73} & q^2_{74} & q^2_{75} & q^2_{76} & q^2_{77} & q^2_{78} \\
%   q^2_{81} & q^2_{82} & q^2_{83} & q^2_{84} & q^2_{85} & q^2_{86} & q^2_{87} & q^2_{88}
%  \end{array}
%  \right)
\end{eqnarray*}
where rows are indexed by the $\{u^k_{A,i}\}$ and columns are indexed by the $\{u^k_{B,j}\}$.
As above, from Lemma \ref{nonzero_derivatives}, we know that
\begin{eqnarray*}
\Delta_1(\pi_4) =  \left | \begin{array}{cc}
q^1_{11} & q^1_{12} \\
q^1_{21} & q^1_{22}
\end{array}
\right | \neq 0, \qquad \mbox{ and } \qquad
\Delta_2(\pi_4)  =  \left | \begin{array}{cc}
q^2_{11} & q^2_{12} \\
q^2_{21} & q^2_{22}
\end{array}
\right | \neq 0.
\end{eqnarray*}
The resulting edge invariants arise from each $M_k$ as rank restrictions for the 3-minors containing $\Delta_k(\pi_4)$, namely:
\begin{eqnarray*}
q^k_{ij} \Delta_k(\pi_4) -q^k_{2j}(q^k_{11}q^k_{i2}-q^k_{12}q^k_{i1})+
q^k_{1j}(q^k_{21}q^k_{i2}-q^k_{22}q^k_{i1})=0, \quad \mbox{ for }i,j\geq 3,\;  k=1,2.
\end{eqnarray*}
These invariants together with the 12 equations obtained from $T_A$ and the 12 equations obtained from $T_B$ define a complete intersection for the variety of $T$ locally near the points of no evolution.

%\newpage
%\begin{sidewaysfigure}
%\begin{figure}
\begin{landscape}
\begin{table}[t]
$$
\begin{array}{l|cccccccc|cccccccc}
& &(*)&(*)&(*)&(*)& &(*)&(*)&(*)&(*) & &(*)&(*)&(*)&(*)& \\
 & \aa^*\otimes \aa &\aa^*\otimes \cc & \aa^*\otimes \gg & \aa^*\otimes \tt & \cc^*\otimes \aa & \cc^*\otimes \cc & \cc^*\otimes \gg & \cc^*\otimes \tt & \gg^*\otimes \aa & \gg^*\otimes \cc & \gg^*\otimes \gg & \gg^*\otimes \tt &
 \tt^*\otimes \aa & \tt^*\otimes \cc & \tt^*\otimes \gg & \tt^*\otimes \tt \\
 \hline
\aa \otimes (\aa \wedge \cc) (*)&-a112&  \cellcolor{lightgray} a111&  \cellcolor{lightgray} 0 &  \cellcolor{lightgray}0 &  \cellcolor{lightgray} -a212& a211&  \cellcolor{lightgray} 0&  \cellcolor{lightgray}0&  \cellcolor{lightgray}-a312 &  \cellcolor{lightgray}a311&  0&  \cellcolor{lightgray}0&  \cellcolor{lightgray} -a412&  \cellcolor{lightgray}a411&  \cellcolor{lightgray}0 &0\\
\aa \otimes (\aa \wedge \gg) (*)&-a113 &  \cellcolor{lightgray}0 & \cellcolor{lightgray}a111&   \cellcolor{lightgray}0 & \cellcolor{lightgray}-a213& 0&  \cellcolor{lightgray}a211 & \cellcolor{lightgray}0&  \cellcolor{lightgray}-a313 & \cellcolor{lightgray}0&  a311&  \cellcolor{lightgray}0&  \cellcolor{lightgray}-a413&  \cellcolor{lightgray}0&  \cellcolor{lightgray}a411& 0\\
\aa \otimes (\aa \wedge \tt)(*)&-a114&  \cellcolor{lightgray}0&  \cellcolor{lightgray}0&  \cellcolor{lightgray}a111&  \cellcolor{lightgray}-a214 &0&  \cellcolor{lightgray}0&  \cellcolor{lightgray}a211&  \cellcolor{lightgray}-a314&  \cellcolor{lightgray}0&  0&  \cellcolor{lightgray}a311&  \cellcolor{lightgray}-a414&  \cellcolor{lightgray}0 & \cellcolor{lightgray}0 &a411\\
\aa \otimes (\cc \wedge \gg) & 0& -a113& a112& 0& 0& -a213 &a212& 0& 0& -a313& a312& 0& 0& -a413& a412& 0\\
\aa \otimes (\cc \wedge \tt) & 0 &-a114& 0& a112& 0 &-a214& 0& a212 &0 &-a314& 0& a312& 0& -a414& 0& a412\\
\aa \otimes (\gg \wedge \tt) & 0 &0 &-a114& a113& 0& 0& -a214 &a213& 0& 0& -a314 &a313& 0& 0& -a414& a413\\
\cc \otimes (\aa \wedge \cc) (*)&-a122& \cellcolor{lightgray}a121& \cellcolor{lightgray}0& \cellcolor{lightgray}0& \cellcolor{lightgray}-a222& a221& \cellcolor{lightgray}0& \cellcolor{lightgray}0& \cellcolor{lightgray}-a322& \cellcolor{lightgray}a321& 0& \cellcolor{lightgray}0& \cellcolor{lightgray}-a422& \cellcolor{lightgray}a421& \cellcolor{lightgray}0& 0\\
\cc \otimes (\aa \wedge \gg) &-a123& 0 &a121& 0 &-a223& 0& a221& 0& -a323& 0& a321& 0& -a423& 0 &a421& 0\\
\cc \otimes (\aa \wedge \tt) &-a124& 0& 0& a121& -a224& 0& 0& a221& -a324& 0 &0& a321& -a424& 0& 0& a421\\
\cc \otimes (\cc \wedge \gg) (*)& 0& \cellcolor{lightgray}-a123& \cellcolor{lightgray}a122& \cellcolor{lightgray}0& \cellcolor{lightgray}0& -a223& \cellcolor{lightgray}a222& \cellcolor{lightgray}0& \cellcolor{lightgray}0& \cellcolor{lightgray}-a323& a322& \cellcolor{lightgray}0& \cellcolor{lightgray}0& \cellcolor{lightgray}-a423& \cellcolor{lightgray}a422& 0\\
\cc \otimes (\cc \wedge \tt) (*)& 0 &\cellcolor{lightgray}-a124& \cellcolor{lightgray}0& \cellcolor{lightgray}a122& \cellcolor{lightgray}0& -a224& \cellcolor{lightgray}0& \cellcolor{lightgray}a222& \cellcolor{lightgray}0 &\cellcolor{lightgray} -a324& 0& \cellcolor{lightgray}a322& \cellcolor{lightgray}0& \cellcolor{lightgray}-a424& \cellcolor{lightgray}0& a422\\
\cc \otimes (\gg \wedge \tt) & 0& 0& -a124& a123& 0& 0& -a224& a223& 0& 0& -a324& a323& 0& 0& -a424& a423\\
\hline
\gg \otimes (\aa \wedge \cc)  &-a132& a131& 0& 0& -a232& a231& 0& 0& -a332& a331& 0& 0& -a432& a431& 0& 0\\
\gg \otimes (\aa \wedge \gg) (*) &-a133& \cellcolor{lightgray}0& \cellcolor{lightgray}a131& \cellcolor{lightgray}0& \cellcolor{lightgray}-a233& 0& \cellcolor{lightgray}a231& \cellcolor{lightgray}0& \cellcolor{lightgray}-a333& \cellcolor{lightgray}0& a331& \cellcolor{lightgray}0& \cellcolor{lightgray}-a433& \cellcolor{lightgray}0& \cellcolor{lightgray}a431& 0\\
\gg \otimes (\aa \wedge \tt) &-a134& 0& 0& a131& -a234& 0& 0& a231& -a334& 0& 0& a331& -a434& 0& 0& a431\\
\gg \otimes (\cc \wedge \gg) (*) & 0& \cellcolor{lightgray}-a133& \cellcolor{lightgray}a132& \cellcolor{lightgray}0& \cellcolor{lightgray}0& -a233& \cellcolor{lightgray}a232& \cellcolor{lightgray}0& \cellcolor{lightgray}0& \cellcolor{lightgray}-a333& a332 &\cellcolor{lightgray}0 &\cellcolor{lightgray}0 &\cellcolor{lightgray}-a433& \cellcolor{lightgray}a432& 0\\
\gg \otimes (\cc \wedge \tt) & 0 &-a134& 0& a132& 0& -a234& 0& a232& 0& -a334& 0& a332& 0& -a434& 0& a432\\
\gg \otimes (\gg \wedge \tt) (*) & 0& \cellcolor{lightgray}0& \cellcolor{lightgray}-a134& \cellcolor{lightgray}a133 &\cellcolor{lightgray}0& 0 &-\cellcolor{lightgray}a234& \cellcolor{lightgray}a233& \cellcolor{lightgray}0& \cellcolor{lightgray}0& -a334& \cellcolor{lightgray}a333& \cellcolor{lightgray}0& \cellcolor{lightgray}0& \cellcolor{lightgray}-a434& a433\\
\tt \otimes (\aa \wedge \cc) &-a142& a141& 0& 0& -a242& a241& 0& 0& -a342& a341& 0& 0& -a442& a441& 0& 0\\
\tt \otimes (\aa \wedge \gg) &-a143& 0& a141& 0& -a243& 0 &a241 &0& -a343& 0& a341& 0 &-a443 &0& a441& 0\\
\tt \otimes (\aa \wedge \tt) (*) &-a144 &\cellcolor{lightgray}0& \cellcolor{lightgray}0& \cellcolor{lightgray}a141& \cellcolor{lightgray}-a244& 0& \cellcolor{lightgray}0& \cellcolor{lightgray}a241& \cellcolor{lightgray}-a344& \cellcolor{lightgray}0& 0& \cellcolor{lightgray}a341& \cellcolor{lightgray}-a444& \cellcolor{lightgray}0& \cellcolor{lightgray}0& a441\\
\tt \otimes (\cc \wedge \gg) &0& -a143& a142& 0& 0& -a243& a242& 0& 0& -a343& a342& 0& 0& -a443 &a442& 0\\
\tt \otimes (\cc \wedge \tt) (*) & 0 &\cellcolor{lightgray}-a144& \cellcolor{lightgray}0& \cellcolor{lightgray}a142& \cellcolor{lightgray}0& -a244& \cellcolor{lightgray}0& \cellcolor{lightgray}a242& \cellcolor{lightgray}0& \cellcolor{lightgray}-a344& 0& \cellcolor{lightgray}a342& \cellcolor{lightgray}0& \cellcolor{lightgray}-a444& \cellcolor{lightgray}0& a442\\
\tt \otimes (\gg \wedge \tt) (*) & 0 & \cellcolor{lightgray}0& \cellcolor{lightgray}-a144 & \cellcolor{lightgray}a143 & \cellcolor{lightgray}0 &0 &\cellcolor{lightgray} -a244 & \cellcolor{lightgray}a243 &\cellcolor{lightgray} 0& \cellcolor{lightgray} 0& -a344 & \cellcolor{lightgray}a343 &\cellcolor{lightgray}0& \cellcolor{lightgray}0 & \cellcolor{lightgray} -a444 &a443
\end{array}
$$

%\caption{\label{fig_matrix} Matrix associated to $f(T)$ in rows labelled by $e_i \otimes (e_j \wedge e_k)$ and columns labelled by $e_i^* \otimes e_j$.}

%\end{figure}
%\end{sidewaysfigure}
\caption{\label{matrix} Matrix representation of $f(\tau)$ for the GMM used in section 6.2. Light gray cells represent the entries lying in distinguished rows and columns.}
\end{table}
\end{landscape}
%\newpage

%\newpage
%\begin{sidewaysfigure}
%\begin{figure}
\begin{landscape}
\begin{table}[t]
$$
\begin{array}{l|cccccccc|cccccccc}
& &(*)&(*)&(*)&(*)& &(*)&(*)&(*)&(*) & &(*) & (*)&(*)&(*)& \\
&\baa^*\otimes \baa &\baa^*\otimes \bcc & \baa^*\otimes \bgg & \baa^*\otimes \btt & \bcc^*\otimes \baa & \bcc^*\otimes \bcc & \bcc^*\otimes \bgg & \bcc^*\otimes \btt & \bgg^*\otimes \baa & \bgg^*\otimes \bcc & \bgg^*\otimes \bgg & \bgg^*\otimes \btt &
 \btt^*\otimes \baa & \btt^*\otimes \bcc & \btt^*\otimes \bgg & \btt^*\otimes \btt \\
 \hline
\baa \otimes (\baa \wedge \bcc) (*)& 0 & \cellcolor{lightgray}q111 & \cellcolor{lightgray} 0 & \cellcolor{lightgray} 0 & \cellcolor{lightgray} -q212& 0& \cellcolor{lightgray} 0&  \cellcolor{lightgray} 0&  \cellcolor{lightgray} -q312&  \cellcolor{lightgray} 0& 0& \cellcolor{lightgray}0& \cellcolor{lightgray}0& \cellcolor{lightgray}q411 & \cellcolor{lightgray}0& 0\\
\baa \otimes (\baa \wedge \bgg) (*) & 0 & \cellcolor{lightgray}0 & \cellcolor{lightgray}q111 & \cellcolor{lightgray}0 &\cellcolor{lightgray} -q213& 0& \cellcolor{lightgray}0& \cellcolor{lightgray}0& \cellcolor{lightgray}-q313& \cellcolor{lightgray}0& 0& \cellcolor{lightgray}0& \cellcolor{lightgray}0& \cellcolor{lightgray}0& \cellcolor{lightgray} q411& 0\\
\baa \otimes (\baa \wedge \btt) (*) &-q114& \cellcolor{lightgray}0& \cellcolor{lightgray}0& \cellcolor{lightgray}q111& \cellcolor{lightgray}0& 0& \cellcolor{lightgray}0& \cellcolor{lightgray}0& \cellcolor{lightgray}0& \cellcolor{lightgray}0& 0& \cellcolor{lightgray}0& \cellcolor{lightgray}-q414& \cellcolor{lightgray}0& \cellcolor{lightgray}0& q411\\
\baa \otimes (\bcc \wedge \bgg)& 0& 0& 0& 0& 0& -q213& q212& 0& 0& -q313& q312& 0& 0& 0& 0& 0\\
 \baa \otimes (\bcc \wedge \btt) &0& -q114& 0& 0& 0& 0& 0& q212& 0& 0& 0& q312& 0& -q414& 0& 0\\
 \baa \otimes (\bgg \wedge \btt) &0& 0& -q114& 0& 0& 0& 0& q213& 0& 0& 0& q313& 0& 0& -q414& 0\\
\bcc \otimes (\baa \wedge \bcc) (*) &-q122& \cellcolor{lightgray}0& \cellcolor{lightgray}0& \cellcolor{lightgray}0& \cellcolor{lightgray}0& q221& \cellcolor{lightgray}0& \cellcolor{lightgray}0& \cellcolor{lightgray}0& \cellcolor{lightgray}q321& 0& \cellcolor{lightgray}0& \cellcolor{lightgray}-q422& \cellcolor{lightgray}0& \cellcolor{lightgray}0& 0\\
\bcc \otimes (\baa \wedge \bgg) &-q123& 0& 0& 0& 0& 0& q221& 0& 0& 0& q321& 0& -q423& 0& 0& 0\\
\bcc \otimes (\baa \wedge \btt)& 0& 0& 0& 0& -q224& 0& 0& q221& -q324& 0& 0& q321& 0& 0& 0& 0\\
\bcc \otimes (\bcc \wedge \bgg)(*)& 0& \cellcolor{lightgray}-q123& \cellcolor{lightgray}q122& \cellcolor{lightgray}0& \cellcolor{lightgray}0& 0& \cellcolor{lightgray}0& \cellcolor{lightgray}0& \cellcolor{lightgray}0& \cellcolor{lightgray}0& 0& \cellcolor{lightgray}0& \cellcolor{lightgray}0& \cellcolor{lightgray}-q423& \cellcolor{lightgray}q422& 0\\
\bcc \otimes (\bcc \wedge \btt) (*) & 0& \cellcolor{lightgray}0& \cellcolor{lightgray}0& \cellcolor{lightgray}q122& \cellcolor{lightgray}0& -q224& \cellcolor{lightgray}0& \cellcolor{lightgray}0& \cellcolor{lightgray}0& \cellcolor{lightgray}-q324& 0& \cellcolor{lightgray}0& \cellcolor{lightgray}0& \cellcolor{lightgray}0& \cellcolor{lightgray}0& q422\\
\bcc \otimes (\bgg \wedge \btt)& 0& 0& 0& q123& 0& 0& -q224& 0& 0& 0& -q324& 0& 0& 0& 0& q423\\
\hline
%\\
\bgg \otimes (\baa \wedge \bcc) & -q132& 0& 0& 0& 0& q231& 0& 0& 0& q331& 0& 0& -q432& 0& 0& 0\\
\bgg \otimes (\baa \wedge \bgg)(*) & -q133& \cellcolor{lightgray}0& \cellcolor{lightgray}0& \cellcolor{lightgray}0& \cellcolor{lightgray}0& 0& \cellcolor{lightgray}q231& \cellcolor{lightgray}0& \cellcolor{lightgray}0& \cellcolor{lightgray}0& q331& \cellcolor{lightgray}0& \cellcolor{lightgray}-q433& \cellcolor{lightgray}0& \cellcolor{lightgray}0& 0\\
\bgg \otimes (\baa \wedge \btt) & 0& 0& 0& 0& -q234& 0& 0& q230& -q334& 0& 0& q330& 0& 0& 0& 0\\
\bgg \otimes (\bcc \wedge \bgg)(*)&  0& \cellcolor{lightgray}-q133& \cellcolor{lightgray}q132& \cellcolor{lightgray}0& \cellcolor{lightgray}0& 0& \cellcolor{lightgray}0& \cellcolor{lightgray}0& \cellcolor{lightgray}0& \cellcolor{lightgray}0& 0& \cellcolor{lightgray}0& \cellcolor{lightgray}0& \cellcolor{lightgray}-q433& \cellcolor{lightgray}q432& 0\\
 \bgg \otimes (\bcc \wedge \btt) &0& 0& 0& q132& 0& -q234& 0& 0& 0& -q334& 0& 0& 0& 0& 0& q432\\
 \bgg \otimes (\bgg \wedge \btt)(*) &0& \cellcolor{lightgray}0& \cellcolor{lightgray}0& \cellcolor{lightgray}q133& \cellcolor{lightgray}0& 0& \cellcolor{lightgray}-q234& \cellcolor{lightgray}0& \cellcolor{lightgray}0& \cellcolor{lightgray}0& -q334& \cellcolor{lightgray}0& \cellcolor{lightgray}0& \cellcolor{lightgray}0& \cellcolor{lightgray}0& q433\\
\btt \otimes (\baa \wedge \bcc) &0& q141& 0& 0& -q242& 0& 0& 0& -q342& 0& 0& 0& 0& q441& 0& 0\\
\btt \otimes (\baa \wedge \bgg) & 0& 0& q141& 0& -q243& 0& 0& 0& -q343& 0& 0& 0& 0& 0& q441& 0\\
\btt \otimes (\baa \wedge \btt) (*) &-q144& \cellcolor{lightgray}0& \cellcolor{lightgray}0& \cellcolor{lightgray}q141& \cellcolor{lightgray}0& 0& \cellcolor{lightgray}0& \cellcolor{lightgray}0& \cellcolor{lightgray}0& \cellcolor{lightgray}0& 0& \cellcolor{lightgray}0& \cellcolor{lightgray}-q444& \cellcolor{lightgray}0& \cellcolor{lightgray}0& q441\\
\btt \otimes (\bcc \wedge \bgg) & 0& 0& 0& 0& 0& -q243& q242& 0& 0& -q343& q342& 0& 0& 0& 0& 0\\
 \btt \otimes (\bcc \wedge \btt) (*) &0& \cellcolor{lightgray}-q144& \cellcolor{lightgray}0& \cellcolor{lightgray}0& \cellcolor{lightgray}0& 0& \cellcolor{lightgray}0& \cellcolor{lightgray}q242& \cellcolor{lightgray}0& \cellcolor{lightgray}0& 0& \cellcolor{lightgray}q342& \cellcolor{lightgray}0& \cellcolor{lightgray}-q444& \cellcolor{lightgray}0& 0\\
 \btt \otimes (\bgg \wedge \btt) (*) & 0& \cellcolor{lightgray}0& \cellcolor{lightgray}-q144& \cellcolor{lightgray}0& \cellcolor{lightgray}0& 0& \cellcolor{lightgray}0& \cellcolor{lightgray}q243& \cellcolor{lightgray}0& \cellcolor{lightgray}0& 0& \cellcolor{lightgray}q343& \cellcolor{lightgray}0& \cellcolor{lightgray}0& \cellcolor{lightgray} -q444& 0\\
\end{array}
$$

\caption{\label{matrixSSM} Matrix representation of $f(\tau)$ in the Fourier basis for SSM used in section \ref{sec:ssm}. Light gray cells represent the entries lying in distinguished rows and columns.}
\end{table}
\end{landscape}

%\section{Acknowledgements}

%\newpage

\bibliographystyle{plain}
%\bibliography{biblio}

\end{document}